\documentstyle[12pt]{article}  
\def\sq{\hbox {\rlap{$\sqcap$}$\sqcup$}}
\def\sq{\hbox {\rlap{$\sqcap$}$\sqcup$}}
\def\R{ {\rm R \kern -.31cm I \kern .15cm}}
\def\C{ {\rm C \kern -.15cm \vrule width.5pt \kern .12cm}}
\def\Z{ {\rm Z \kern -.27cm \angle \kern .02cm}}
\def\N{ {\rm N \kern -.26cm \vrule width.4pt \kern .10cm}}
\def\1{{\rm 1\mskip-4.5mu l} }
\def\lsim{\raise0.3ex\hbox{$<$\kern-0.75em\raise-1.1ex\hbox{$\sim$}}}
\def\gsim{\raise0.3ex\hbox{$>$\kern-0.75em\raise-1.1ex\hbox{$\sim$}}}
\def\noi{\noindent}

\def\beq{\begin{equation}}   \def\eeq{\end{equation}}
\def\bea{\begin{eqnarray}}  \def\eea{\end{eqnarray}}
\def\nn{\nonumber}
\def\noi{\noindent}
\def\beeq{\begin{eqnarray}} \def\eeeq{\end{eqnarray}}
\newcommand\mysection{\setcounter{equation}{0}\section}

\newcounter{hran}

\begin{document} 
\centerline{\large\bf Scattering Theory for the Schr\"odinger Equation} 
 \vskip 3 truemm \centerline{\large\bf in some External Time Dependent Magnetic Fields\footnote{Work supported in part by NATO Collaborative Linkage Grant 979341}} 

\vskip 0.5 truecm

\centerline{\bf J. Ginibre}
\centerline{Laboratoire de Physique Th\'eorique\footnote{Unit\'e Mixte de
Recherche (CNRS) UMR 8627}}  \centerline{Universit\'e de Paris XI, B\^atiment
210, F-91405 ORSAY Cedex, France}
\vskip 3 truemm
\centerline{\bf G. Velo}
\centerline{Dipartimento di Fisica, Universit\`a di Bologna}  \centerline{and INFN, Sezione di
Bologna, Italy}

\vskip 1 truecm

\begin{abstract}
We study the theory of scattering for a Schr\"odinger equation in an
external time dependent magnetic field in the Coulomb gauge, in space dimension 3. The
magnetic vector potential is assumed to satisfy decay properties in time that
are typical of solutions of the free wave equation, and even in some
cases to be actually a solution of that equation. That problem appears
as an intermediate step in the theory of scattering for the
Maxwell-Schr\"odinger (MS) system. We prove in particular the existence
of wave operators and their asymptotic completeness in spaces of
relatively low regularity. We also prove their existence or at least
asymptotic results going in that direction in spaces of higher
regularity. The latter results are relevant for the MS system. As a
preliminary step, we study the Cauchy problem for the original equation
by energy methods, using as far as possible time derivatives instead of
space derivatives. \end{abstract}

\vskip 3 truecm
\noi AMS Classification : Primary 35P25. Secondary 35B40, 35Q40, 81U99.  \par \vskip 2 truemm

\noi Key words : Schr\"odinger equation, time dependent magnetic fields, Scattering\break \noindent theory. 
\vskip 1 truecm

\noindent LPT Orsay 04-02\par
\noindent January 2004\par

\newpage
\pagestyle{plain}
\baselineskip 18pt

\mysection{Introduction}
\hspace*{\parindent}
This paper is devoted to the theory of scattering and in particular to
the construction of the wave operators for a Schr\"odinger equation
minimally coupled to an external time dependent magnetic field in the
Coulomb gauge, namely
\beq \label{1.1e} i \partial_t u = - (1/2) \Delta_A u \eeq

\noi in space dimension 3. Here $u$ is a complex function defined in
space time ${I\hskip-1truemm R}^{3+1}$,
\beq \label{1.2e} \Delta_A \equiv \nabla_A^2 \equiv (\nabla - i A)^2
\eeq

\noi and the magnetic potential $A$ is an ${I\hskip-1truemm R}^3$
vector valued function defined in ${I\hskip-1truemm R}^{3+1}$ and
satisfying the condition $\nabla\cdot A = 0$, which is the Coulomb gauge
condition. The magnetic potential will be assumed to be sufficiently
smooth and to satisfy a number of decay estimates in time that are
satisfied by sufficiently regular solutions of the free wave equation
$\sq A = 0$, where $\sq$ is the d'Alembertian operator. At some places
$A$ will even be assumed to be a solution of that equation. \par

The present problem arises, actually is an intermediate step, in the
theory of scattering for the Maxwell-Schr\"odinger (MS) system. In the
Coulomb gauge\break\noindent $\nabla\cdot A = 0$, that system takes the form 
\beq
\label{1.3e} \left \{ \begin{array}{l}  i\partial_t u = - (1/2)
\Delta_A u + A_0 u \\ \\ \Delta A_0 = - |u|^2 \\ \\ \sq A + \nabla
\left ( \partial_t A_0 \right ) = {\rm Im} \ \bar{u} \nabla_Au \equiv
J(u, A) \ .   \end{array} \right . \eeq

\noi Using the second equation in (\ref{1.3e}) to eliminate $A_0$, one
can recast that system into the formally equivalent form
\beq
\label{1.4e} \left \{ \begin{array}{l}  i\partial_t u = - (1/2)
\Delta_A u + g(u)\  u \\ \\  \sq A = P \ {\rm Im}\ \bar{u}  \nabla_A u = P\ J(u, A)    \end{array} \right . \eeq

\noi where $g(u)$ is the Coulomb interaction term
$$g(u) = (4 \pi |x|)^{-1} \ * \ |u|^2$$

\noi and $P = \1 - \nabla \Delta^{-1} \nabla$ is the projector on
divergence free vector fields. \par

The theory of scattering for the MS system (\ref{1.4e}) has been
studied and in particular the construction of modified wave operators
has been performed in \cite{13r} \cite{16r} in the case of small
asymptotic states and solutions and in \cite{7r} for asymptotic states
and solutions of arbitrary size for the Schr\"odinger function, but
only in the special case of vanishing asymptotic magnetic field, namely
in the case where the asymptotic state for the magnetic potential is
zero. More precisely the method used in \cite{7r} starts with the
replacement of the Maxwell equation by the associated integral equation
with infinite initial time, namely 
\beq \label{1.5e} A = A_0 + A_1 =
\dot{K}(t) A_+ + K(t) \dot{A}_+ - \int_t^{\infty} dt' K (t - t') PJ
(u,A) (t') \eeq

\noi where 
\beq \label{1.6e} K(t) = \omega^{-1} \sin \omega t \ , \
\dot{K}(t) = \cos \omega t \ , \ \omega = (- \Delta)^{1/2} \ . \eeq

\noi The case treated in \cite{7r} is that where the asymptotic state
$(A_+, \dot{A}_+)$ for the magnetic potential is zero. As an
intermediate step towards the treatment of the general case, it is
useful to consider the complementary case where only the free part
$A_0$ is kept in the equation (\ref{1.5e}). If in addition one omits
the now well controlled Hartree interaction $g(u)$, one is led to
study the theory of scattering for (\ref{1.1e}) where now
$\sq A = 0$. This is the purpose of the present paper.\par

The variables $(u, A)$ are not convenient to study the asymptotic
behaviour in time and the theory of scattering for (\ref{1.1e}). The
present paper will rely in an essential way on the use of new variables
defined as follows. The unitary group which solves the free
Schr\"odinger equation can be written as \beq \label{1.7e} U(t) = \exp
\left ( i(t/2)\Delta \right ) = M(t)\ D(t)\ F \ M(t) \eeq

\noi where $M(t)$ is the operator of multiplication by the function 
\beq
\label{1.8e}
M(t) = \exp \left ( i x^2/2t\right ) \ ,
\eeq

\noi $F$ is the Fourier transform and $D(t)$ is the dilation operator defined by 
\beq
\label{1.9e}
D(t) = (it)^{-3/2} \ D_0(t) \ ,
\eeq
\beq
\label{1.10e}
\left ( D_0(t) f\right ) (x) = f(x/t) \ .
\eeq

\noi We replace the variables $(u, A)$ by new variables $(w, B)$
defined by 
\beq \label{1.11e} \left \{ \begin{array}{l}  u(t) = M(t) \
D(t) \ \overline{w(1/t)}\\ \\  A(t) = - t^{-1} \ D_0 (t) \ B(1/t) \ .  
\end{array} \right . \eeq

\noi Actually the change of variables from $u$ to $w$ is simply the
pseudoconformal inversion. In terms of the new variables $(w, B)$ the
original equation, namely (\ref{1.1e}), is easily seen to become 
\beq
\label{1.12e} i \partial_t w = - (1/2) \Delta_B w - \hbox{\it \v B} w \eeq

\noi where 
\beq \label{1.13e} \hbox{\it \v B(t)} = t^{-1} (x \cdot B(t)) \
. \eeq

\noi The study of the asymptotic behaviour in time and of the theory of
scattering for (\ref{1.1e}) is then reduced to the study of
(\ref{1.12e}) near $t= 0$ and will be performed by studying that
equation in the interval $[0, 1]$. \par

In \cite{7r} we performed a different change of variables and used new
variables called $(w, B)$ in \cite{7r} and which we now denote $(w_*,
B_*)$. They are defined by 
\beq \label{1.14e} \left \{ \begin{array}{l}
 u(t) = M(t) \ D(t) \ w_*(t)\\ \\  A(t) = t^{-1} \ D_0 (t) \ B_*(t)  
\end{array} \right . \eeq

\noi so that 
\beq \label{1.15e} w(t) = \overline{w_*(1/t)}\quad , \quad
B(t) = - B_*(1/t) \ . \eeq

In terms of the variables $(w_*, B_*)$, (\ref{1.1e}) becomes \beq
\label{1.16e} i \partial_t w_* = - \left ( 2t^2 \right )^{-1}
\Delta_{B_*} w_* - \hbox{\it \v B}_* w_* \eeq

\noi where 
\beq \label{1.17e} \hbox{\it \v B}_*(t) = t^{-1} \left ( x
\cdot B_*(t) \right ) = - t^{-1} \left ( x \cdot B(1/t)\right ) \ .
\eeq

\noi The study of the asymptotic behaviour in time and the theory of
scattering for (\ref{1.1e}) is then reduced to the same problem for
(\ref{1.16e}). The change of variables from $u$ to $w$ or $w_*$ can be
rewritten in a slightly different way by introducing 
\beq \label{1.18e}
\widetilde{u}(t) = U(-t) \ u(t) \quad , \quad \widetilde{w}(t) = U(-t)
\ w(t) \ . \eeq

\noi Then 
\beq \label{1.19e} F\widetilde{u}(t) =
\overline{\widetilde{w}(1/t)} = U(1/t) w_*(t) \ . \eeq

Scattering theory is the asymptotic study at infinity in time of an
evolution, in the present case that defined by (\ref{1.1e}) or
(\ref{1.16e}), by comparison with a simpler reference evolution. In the
present case the latter is inspired by the free Schr\"odinger equation
and will be chosen to drive $\widetilde{u}$ or $w_*$ to a limit as $t
\to \infty$, or equivalently to drive $w$ to a limit as $t \to 0$. The
simplest evolution of this type is the free Schr\"odinger evolution
itself for $u$, which is equivalent to taking $\widetilde{u}$ constant
or equivalently $\widetilde{w}$ constant. One is then led to study the
convergence properties of $\widetilde{u}$, $\widetilde{w}$ or possibly
of $w$, $w_*$ in a suitable function space $X$ in the appropriate time
limit.\par

The theory of scattering thereby obtained will in general depend
significantly on the choice of the space $X$. Let $X$ and $Y$ be two
Banach spaces with $Y \subset X$. When comparing the theories of
scattering in $X$ and $Y$, one may encounter the following situations.
On the one hand, one can in some cases use the theory of scattering in
$Y$ together with uniform bounds on the evolution in $X$ to construct
the theory of scattering in $X$. This possibility will be exploited in Section 4
below with $X = L^2$ and $Y = H^2$. On the other hand, when restricting
the theory of scattering from $X$ to $Y$, one may eliminate some
(insufficiently regular) solutions of the original equation, thereby
restricting the possible set of asymptotic behaviours. However since
one is then interested in convergence properties in a stronger sense,
namely in the norm of $Y$ instead of that of $X$, it may also happen
that the asymptotic behaviours that were sufficiently accurate in $X$
norm are no longer so in $Y$ norm and have to be replaced by more
accurate ones. This last possibility will appear in Sections 6 and 7
below when going from $X = H^2$ to $Y = H^k$ for $3 \leq k \leq 4$.\par

We now comment on the methods used in the present paper. In order to
treat the full nonlinear MS system (\ref{1.4e}) we need to use spaces
of sufficiently regular functions, more precisely in \cite{7r} we
needed to take $w$ or $w_*$ in $H^k$ for $k > 5/2$. Furthermore we have
to use methods that are sufficiently simple and robust to accomodate
the non linearities. This is the case of the energy methods used in
\cite{7r} and we shall confine our attention to such methods in this
paper. The Schr\"odinger equation with time dependent magnetic fields
has been studied in the literature by more sophisticated methods using
the abstract theory of evolution equations or semi classical
approximations. We refer to \cite{17r} and references therein quoted.
Those methods do not seem to be readily extendable to the MS system,
and we shall not use them in this paper. \par

We now turn to an important feature of the MS system, namely the fact
that it couples two equations with different scaling properties. In
fact one time derivative is homogeneous to two space derivatives for
the Schr\"odinger equation but only to one space derivative for the
wave equation. When treating the MS system with $u \in H^{2k}$ by using
space derivatives, one needs to consider $\Delta^ku$ and
correspondingly $\partial_x^{\alpha}A$ for multiindices $\alpha$ with
$|\alpha| = 2k$. If instead one uses time derivatives, one needs to
consider $\partial_t^k u$ and correspondingly $\partial_t^kA$ which is
homogeneous to $\partial_x^{\alpha}	A$ with $|\alpha | = k$ only.
When considering time decay properties, this fact makes little
difference in the case of the variables $(u, A)$ because for solutions
of the free wave equation $\sq A = 0$, one can ensure that
$\partial_t^j \partial_x^{\alpha} A$ has the same time decay as $A$ for
any $j$ and $\alpha$ by considering sufficiently regular solutions. However the same phenomenon
occurs when studying the MS system in the transformed variables $(w,
B)$. In that case it follows from (\ref{1.11e}) that
derivatives $\partial_t^j \partial_x^{\alpha}$ applied to $B$ generate
a factor $t^{-j-|\alpha|}$ no matter how regular $A$ is, which is a
disaster as regards the behaviour at $t \to 0$. It is therefore
important to apply as few derivatives as possible to $B$ and for that
purpose it is advantageous to use time derivatives instead of space
derivatives of $w$ in order to control the regularity in spaces
$H^{2k}$. In this paper we shall therefore use time derivatives as much
as possible in the treatment of (\ref{1.12e}). On the other hand, since
as mentioned above we need simple methods that can be extended to the
full MS system, we shall use only time (and also space) derivatives of
integer order. \par

We now turn to a description of the contents of this paper
and to a statement of a representative sample of the results. After a preliminary section
containing notation and estimates of a general nature (Section 2), we
begin with the study of the Cauchy problem at finite times for a class
of equations (see (\ref{3.1e}) (\ref{3.2e}) below) that generalize
slightly both (\ref{1.1e}) and (\ref{1.12e}) (Section 3). The main
result is that the Cauchy problem is well posed in $L^2$ (Proposition
3.1), in $H^2$ (Proposition 3.2), in $H^3$ (Proposition 3.3) and at the
level of $H^4$ (Proposition 3.4) under assumptions on $A$ (or $B$) of a
general nature. Let $v$ denote the unknown function. Following the
remarks made above, the method is based on energy estimates of $v$, of
$\partial_t v$, of $\nabla \partial_t v$ and of $\partial_t^2 v$
respectively, and the number of derivatives on $A$ or $B$ in the
assumptions is kept to a minimum, namely 0 or 1, 1, 2 and 2
respectively. \par

We next turn to the theory of scattering for (\ref{1.1e})
and (\ref{1.16e}) via the study of (\ref{1.12e}) in $[0, 1]$, in the
more particular case where $A$ satisfies some of the natural decay
properties associated with the free wave equation. In Section 4, we
exploit Propositions 3.1 and 3.2 to prove the existence and the
asymptotic completeness of the wave operators for the equation
(\ref{1.1e}) with $u$ or rather $\widetilde{u}$ in $L^2$ and in $FH^2$,
as compared with the free Schr\"odinger evolution, or equivalently by
(\ref{1.19e}), for the equation (\ref{1.16e}) with $w_*$ in $L^2$ and
in $H^2$, as compared with the constant evolution. The main result can
be stated as follows (see Propositions 4.1 and 4.2 below for more
details).\\

\noi{\bf Proposition 1.1.} {\it Let $A$ satisfy 
\beq \label{1.20e} \parallel P^j\ \partial_x^{\alpha} A
\parallel_r \ \vee \ \parallel P^j(x\cdot A) \parallel_r \ \leq C \ t^{-1 + 2/r} \eeq

\noi where $P = t \partial_t + x \cdot \nabla$, for $0 \leq j + |\alpha | \leq 1$, $2 \leq r \leq \infty$ and for all
$t \in [1, \infty )$. \par

(1) Let $X = L^2$ or $FH^2$. Then for any $u_+ \in X$, there exists a
unique solution $u$ of (\ref{1.1e}) such
that $\widetilde{u} \in {\cal C}([1, \infty ), X)$ and such
that
\beq
\label{1.21e}
\parallel \widetilde{u}(t) - u_+; X \parallel \ \to 0 \quad \hbox{\it when $t \to \infty$} \ .
\eeq

\noi Conversely for any solution $u$ of (\ref{1.1e}) such
that $\widetilde{u} \in {\cal C}([1, \infty ), X)$, there exists $u_+ \in X$ such that (\ref{1.21e}) holds.

(2) Let $X = L^2$ or $H^2$. Then for any $w_+ \in X$, there exists a
unique solution $w_* \in {\cal C}([1, \infty ), X)$ of (\ref{1.16e})
such that 
\beq
\label{1.22e}
\parallel w_*(t) - w_+;X\parallel \ \to 0\quad \hbox{\it when $t \to \infty$} \ .
\eeq

\noi Conversely for any solution $w_* \in {\cal C}([1, \infty ), X)$ of
(\ref{1.16e}), there exists $w_+ \in X$ such that (\ref{1.22e})
holds.} \\

Note that the time decay that occurs in (\ref{1.20e}) is the optimal
time decay that can be obtained for the relevant norms of $A$ if $A$ is
a solution of the wave equation $\sq A = 0$ satisfying the Coulomb
gauge condition $\nabla\cdot A = 0$. In that case, that decay can be
easily ensured by making appropriate assumptions on the Cauchy data for
$A$ (see Section 2 below, especially Lemma 2.4). \par

If one makes stronger assumptions on the asymptotic state $u_+$ or $w_+$, one
obtains stronger convergence properties than (\ref{1.21e})
(\ref{1.22e}) for the solutions $u$ or $w_*$ constructed in Proposition
1.1 as $t \to \infty$. The following typical result is extracted from a
special case of Proposition 4.3 below, to which we refer for a slightly
more general result. \\

\noi {\bf Proposition 1.2.} {\it Let $A$ satisfy the assumptions of
Proposition 1.1.\par

(1) Let $u_+ \in FH^2$ and let $u$ be the solution of (\ref{1.1e}) with
$\widetilde{u} \in {\cal C}([1, \infty ), FH^2)$ obtained in
Proposition 1.1, part (1). Then the following estimate holds for all $t
\geq 1$~:
\beq
\label{1.23e}
\parallel \widetilde{u}(t) - u_+ \parallel_2 \ \leq \ C\ t^{-1} \ .
\eeq

\noi Let in addition $u_+ \in FH^3$. Then the following estimates hold for all $t \geq 1$~:
\beq
\label{1.24e}
\parallel \widetilde{u}(t) - u_+ \parallel_2 \ \leq \ C\ t^{-3/2} \ ,
\eeq

\beq
\label{1.25e}
t^2 \parallel \partial_t \left ( \widetilde{u}(t) - u_+\right ) \parallel_2 \ \vee  \ \parallel x^2\left ( \widetilde{u}(t) - u_+\right ) \parallel_2\ \leq\ C\ t^{-1/2} \ .
\eeq

(2) Let $w_+ \in H^2$ and let $w_* \in {\cal C} ([1, \infty ), H^2)$ be
the solution of (\ref{1.16e}) obtained in Proposition 1.1, part (2).
Then the following estimate holds for all $t \geq 1$~:
\beq
\label{1.26e}
\parallel w_*(t) - w_+ \parallel_2 \ \leq \ C\ t^{-1} \ .
\eeq

\noi Let in addition $w_+ \in H^3$. Then the following estimates hold for all $t\geq 1$~:}
\beq
\label{1.27e}
\parallel w_*(t) - U^*(1/t) w_+ \parallel_2 \ \leq \ C\ t^{-3/2} \ ,
\eeq
\beq
\label{1.28e}
t^2 \parallel \partial_t \left ( w_*(t) - U^* (1/t) w_+\right ) \parallel_2 \ \vee  \ \parallel \Delta \left ( w_*(t) - U^*(1/t) w_+\right ) \parallel_2\ \leq\ C\ t^{-1/2} \ .
\eeq

We next study the $H^k$ regularity of the $L^2$  wave operators for
$u$ constructed in Section 4, as given in particular by Proposition 1.1
part (1) above with $X = L^2$. As mentioned above, time decay is not
impaired by derivatives for sufficiently regular solutions of the free
wave equation. Making decay assumptions of that type on $A$ at a
sufficient level of regularity, we prove in Section 5 that the $L^2$
wave operators for $u$ essentially preserve $H^k$ regularity for
arbitrarily high $k$. The following typical result is a special case of
Proposition 5.2, to which we refer for a slightly more general
statement.\\

\noi {\bf Proposition 1.3.} {\it Let $j \geq 0$ be an integer and let $A$ satisfy
\beq
\label{1.29e}
\parallel \partial_t^l \
\partial_x^{\alpha} \ A \parallel_r \ \vee \ \parallel \partial_t^l \
\partial_x^{\alpha} (x\cdot A) \parallel_r \ \leq C\ t^{-1+2/r} 
\eeq
\noi for all $r$, $2 \leq r \leq \infty$, for $0 \leq l + |\alpha |/2
\leq j$ and for all $t \geq 1$. Let $u_+$, $x u_+ \in H^{2j} \cap
H_{1}^{2j}$. Then
there exists a unique solution $u \in \displaystyle{\mathrel{\mathop {\rm \cap}_{0
\leq l \leq j}}}{\cal C}^{j-l} ([1, \infty ) , H^{2l})$ of (\ref{1.1e})
satisfying the estimates
\beq
\label{1.30e}
\parallel \partial_t^{j-l} \Delta^l(u(t) - U(t)
u_+)\parallel_2 \ \leq C\ t^{-3/2} 
\eeq
\noi for $0 \leq l \leq j$ and for all $t
\geq 1$. The solution is actually unique in ${\cal C}([1, \infty ),
L^2)$ under the condition (\ref{1.30e}) for $j = 0$.}\\

\noi (See Section 2 below for the definition of $H_1^{2j}$). At the
level of regularity of $H^k$ for $u$, however, we lose control of
asymptotic completeness since the method does not even allow to prove
that generic solutions of (\ref{1.1e}) remain bounded in $H^k$ for $k >
0$ under assumptions on $A$ of the same type (see Proposition 5.1 for a
partial result in that direction). \par

We next turn to the theory of scattering at
the level of regularity of $H^k$ with $k \geq 3$ for $w_*$ or $w$,
which is of primary interest for subsequent application to the MS
system, by studying again (\ref{1.12e}) at that level of regularity. Now
however Propositions 3.3 and 3.4 do not apply with initial time $t_0 =
0$ because the relevant norms of $B$ blow up too fast as $t \to 0$. As
a consequence we cannot prove that the generic solutions of
(\ref{1.12e}) (resp. (\ref{1.16e})) remain bounded as $t \to 0$ (resp.
($t \to \infty$)) and we lose control of asymptotic completeness.
Furthermore, in order to construct the wave operators, which in that
case amounts to solving the Cauchy problem for (\ref{1.12e}) with
initial time zero, we have to resort to the same indirect method that
was used in \cite{7r}. We give ourselves a presumed asymptotic
behaviour of the solution $w$ of (\ref{1.12e}) near zero in the form of a function $W$ with
prescribed value at $t = 0$ and we look for $w$ in the form $w= W + q$ with $q$ tending to zero as $t \to
0$. The evolution equation for $q$ is obviously 
\beq
\label{1.31e}
i\partial_t q = - (1/2) \Delta_B q - \hbox{\it \v B}q - R(W)
\eeq

\noi where
\beq
\label{1.32e}
R(W) = i \partial_t W + (1/2) \Delta_BW + \hbox{\it \v B}W \ .
\eeq

\noi We take $t_0 > 0$, we apply Proposition 3.3 or 3.4 to define a
solution $q_{t_0}$ of (\ref{1.31e}) in $(0, 1]$ with suitably small
initial condition $q_{t_0}(t_0) = q_0$ at $t_0$, and we take the limit
of $q_{t_0}$ as $t_0 \to 0$. We carry out this program in two steps. In
Section 6, we construct solutions $q$ of (\ref{1.31e}) tending to zero
at $t \to 0$ under general assumptions on $W$, taking mainly the form
of time decay estimates of $R(W)$. This is done in $H^3$ in Proposition
6.1 and at the level of $H^4$ in  Proposition 6.2. Actually in the
latter case, although $\partial_t^2 q$ tends to zero in $L^2$ as $t \to
0$, $q$ itself tends to zero in $H^k$ only for some $k$ with $3 < k <
4$, depending on the decay assumptions on $R(W)$, but in general not in
$H^4$. In Section 7 we construct asymptotic functions $W$ satisfying
the assumptions required in Section 6, restricting our attention to the
case where $A$ is a solution of the free wave equation. The simplest
way to do that, in keeping with the fact that we compare (\ref{1.1e})
with the free Schr\"odinger equation, consists in taking $W(t) = U(t)
\overline{w_+}$ where $w_+ = Fu_+$, the Fourier transform of the
Schr\"odinger asymptotic state $u_+$. One then encounters a standard
difficulty in that problem, namely the difference of propagation
properties of the wave equation and of the Schr\"odinger equation
\cite{5r} \cite{16r}. That difficulty can be circumvented by imposing a
support condition on $w_+$ saying in effect that $w_+$ vanishes in a
neighborhood of the unit sphere, so that $u_+$ generates a solution of
the free Schr\"odinger equation which is asymptotically small in a
neighborhood of the light cone. In Section 7 we first produce an
appropriate $W$ under such a support condition. In Proposition 7.1 we
reduce the problem to a joint condition on $B$ and on the support of
$w_+$ and in Proposition 7.2 we prove that such a condition can be
ensured under a suitable support condition on $w_+$ and suitable decay
assumptions at infinity in space on the asymptotic state $(A_+,
\dot{A}_+)$ of $A$ (see (\ref{1.5e})). Collecting the implications of
Section 6 and of that first part of Section 7 on the theory of
scattering for (\ref{1.1e}) at the level of regularity of $FH^4$ for
$\widetilde{u}$, we construct solutions of (\ref{1.1e}) with prescribed
asymptotic state $u_+$ at that level of regularity. The main result is
stated in Section 8 as Proposition 8.1 to which we refer for a full
mathematical statement. Since that statement is too complicated to be
presented at the level of this introduction, we give here only a
heuristic preview thereof, stripped from most technicalities.\\

\noi {\bf Proposition 1.4} {\it Let $A$ be a solution of the free wave
equation satisfying suitable conditions of regularity and of decay at
infinity in space. Let $u_+ \in FH^5$ and let $w_+ \equiv Fu_+ (\in
H^5)$ satisfy the support condition
\beq
\label{1.33e}
{\rm Supp} \ w_+ \subset \left \{ x: |\ |x|-1| \geq \eta \right \}
\eeq

\noi for some $\eta$, $0 < \eta < 1$. Then there exists a unique
solution $u$ of (\ref{1.1e}) such that $\widetilde{u} \in {\cal C} ([1,
\infty ), FH^4) \cap {\cal C}^1([1, \infty ), FH^2) \cap {\cal C}^2([1,
\infty ),L^2)$ and such that $\widetilde{u} - u_+$ tends to zero as $t
\to \infty$ in suitable norms (related to the previous space) with
power law decay in time, with exponents depending on the assumptions
made on $A$.}\\

We next try to eliminate the support condition (\ref{1.33e}). This has
been done in \cite{12r} in the case of the Wave-Schr\"odinger system
(see also \cite{6r}) and we try to apply the same method in the present
case. We are only partly successful, namely we succeed in controlling
the $B\cdot\nabla W$ and $\hbox{\it \v B}W$ terms in $R(W)$, but not
the $B^2W$ term. The results appear in Propositions 7.3 to 7.5 as
regards the construction of the asymptotic $W$, but the situation is
too intricate to be described here whereas the result is hopefully not
final and we refer to the discussion in Section 7 for details.
Collecting the implications of Section 6 and of that second part of
Section 7, we finally construct solutions of (\ref{1.1e}) at the level
of regularity of $FH^4$ for $\widetilde{u}$ with prescribed asymptotic
behaviour in time given by an asymptotic function $\widetilde{u}_a$
constructed from the asymptotic $W$ mentioned above according to
(\ref{8.7e}). The result is stated in Section 8 as Proposition 8.2 to
which we refer for a full mathematical statement. A heuristic overview
of that proposition can be obtained from Proposition 1.4 above by
having the support condition (\ref{1.33e}) removed and by having $u_+$
replaced by the asymptotic function $\widetilde{u}_a$ mentioned above
in the final decay estimates. In addition, the assumptions on $A$ are
slightly different and include some moment conditions.

\mysection{Notation and preliminary estimates} \hspace*{\parindent} In
this section we introduce some notation and collect a number of
estimates which will be used freely throughout this paper. We denote by
$\parallel \cdot \parallel_r$ the norm in $L^r \equiv
L^r({I\hskip-1truemm R}^3)$, $1 \leq r \leq \infty$. For any
nonnegative integer $k$ and for $1 \leq r \leq \infty$, we denote by
$H_r^k$ the Sobolev spaces
$$H_r^k = \left \{ u : \parallel u ; H_r^k\parallel \ = \sum_{\alpha :
0 \leq |\alpha | \leq k} \parallel \partial_x^{\alpha} u \parallel_r \
< \infty \right \}$$

\noi where $\alpha$ is a multiindex. For $1 < r < \infty$, those spaces
can be defined equivalently (with equivalent norms) by

$$H_r^k = \left \{ u : \parallel u; H_r^k\parallel\ = \ \parallel
<\omega >^k u \parallel_r \ < \infty \right \}$$

\noi where $\omega = (- \Delta)^{1/2}$ and $<\cdot > = (1 + |\cdot
|^2)^{1/2}$. The latter definition extends immediately to any $k \in
{I\hskip-1truemm R}$ and we shall occasionally use such spaces. The
subscript $r$ in $H_r^k$ will be omitted in the case $r = 2$. For any
interval $I$ and for any Banach space $X$ we denote by ${\cal C}(I, X)$
(resp. ${\cal C}_w (I, X))$ the space of strongly (resp. weakly)
continuous functions from $I$ to $X$. For any positive integer $k$, we
denote by ${\cal C}^k (I, X)$ the space of $k$ times differentiable
functions from $I$ to $X$. For any $r$, $1 \leq r \leq \infty$, we
denote by $L^r(I, X)$ (resp. $L_{loc}^r(I, X))$ the space of $L^r$
integrable (resp. locally $L^r$ integrable) functions from $I$ to $X$ if
$r < \infty$ and the space of measurable essentially bounded (resp.
locally essentially bounded) functions from $I$ to $X$ if $r = \infty$.
For $I$ an open interval we denote by $D'(I, X)$ the space of vector
valued distributions from $I$ to $X$ \cite{9r}. For any integer $k$, $0
< k < \infty$, for any $r$, $1 \leq r \leq \infty$, we denote by
$H_r^k(I, X)$ (resp. $H_{r,loc}^k(I, X)$) the space of functions from $I$
to $X$ whose derivatives up to order $k$, taken in $D'(I_0, X)$ with
$I_0$ the interior of $I$, are in $L^r(I, X)$ (resp. $L_{loc}^r(I, X)$). In
the same vein we say that an evolution equation (like (\ref{1.1e})) has
a solution in $I$ with values in $X$ if the equation is satisfied in
$D'(I_0, X)$. \par

We shall use extensively the following Sobolev inequalities, stated
here in ${I\hskip-1truemm R}^n$, but to be used only for $n = 3$. \\

\noi {\bf Lemma 2.1.} {\it Let $1 < q$, $r < \infty$, $1 < p \leq
\infty$ and $0 \leq j < n$. If $p = \infty$, assume that $k - j > n/r$.
Let $\sigma$ satisfy $j/k \leq \sigma \leq 1$ and 

\beq \label{2.1e} n/p
- j = (1 - \sigma )n/q + \sigma (n/r - k) \ . \eeq

\noi Then the following inequality holds~:} 
\beq \label{2.2e} \parallel
\omega^j u \parallel_p \ \leq C \parallel u \parallel_q^{1 - \sigma} \
\parallel \omega^ku \parallel_r^{\sigma} \ .  \eeq

\indent The proof follows from the Hardy-Littlewood-Sobolev inequality
\cite{8r} (from the Young inequality if $p = \infty$), from
Paley-Littlewood theory and interpolation.\\

Occasionally a special case of (\ref{2.2e}) will be used with the
ordinary derivative $\nabla$ replaced by the covariant derivative
$\nabla_A = \nabla - iA$, where $A$ is a real vector-valued function,
namely 
\beq \label{2.3e} \parallel u \parallel_p ÷ \leq C \parallel u
\parallel_q^{1 - \sigma}\ \parallel \nabla_A u\parallel_r^{\sigma} \eeq

\noi which holds under the assumptions of Lemma 2.1 with $j = 0$, $k =
1$. The proof of (\ref{2.3e}) is an immediate consequence of
(\ref{2.2e}) with $j = 0$, $k = 1$ applied to $|u|$ and of the
inequality $|\nabla |u|| \leq |\nabla_Au|$.\par

We shall also make use of the following two lemmas.\\

\noi {\bf Lemma 2.2.} {\it Let $0 \leq \alpha_j < 1$, $a_j \in
{I\hskip-1truemm R}^+$, $1 \leq j \leq n$ and let $y \in
{I\hskip-1truemm R}^+$ satisfy
 \beq \label{2.4e} y \leq \sum_{1 \leq j
\leq n} a_j y^{\alpha_j} \ .\eeq 

\noi Then 
\beq \label{2.5e} y \leq C
\sum_{1 \leq j \leq n} (a_j)^{1/(1 - \alpha_j)} \eeq

\noi where} 
$$C = \mathrel{\mathop {\rm Max}_{1 \leq j \leq n}}\
n^{1/(1- \alpha_j)} \ .$$

\noi {\bf Proof.} 
$$\left \{ y : y \leq \sum_{1 \leq j \leq n} a_j
y^{\alpha_j} \right \} \subset \mathrel{\mathop {\cup}_{1 \leq j \leq
n}}\ \left \{ y : y \leq n a_j y^{\alpha_j} \right \} =$$
$$\mathrel{\mathop {\cup}_{1 \leq j \leq n}}\ \left \{ y : y \leq ( n
a_j)^{1/(1 - \alpha_j)} \right \} \subset \left \{ y : y \leq \sum_{1
\leq j \leq n} (n a_j)^{1/(1 - \alpha_j)} \right \}$$

\noi which implies (\ref{2.5e}). \par\nobreak \hfill $\sq$\\

\noi {\bf Lemma 2.3.} {\it Let $\alpha_j$ satisfy $0 \leq \alpha_j <
\alpha_1 < 1$, $2 \leq j \leq n$. Let $I$ be an interval, let $a_j \in
{\cal C}(I,{I\hskip-1truemm R}^+)$, $0 \leq j \leq n$ and let $y \in
{\cal C} (I, {I\hskip-1truemm R}^+)$ be absolutely continuous and satisfy 
\beq \label{2.6e} \left |
\partial_t y \right | \leq a_0y + \sum_{1 \leq j \leq n} a_j
y^{\alpha_j} \eeq 

\noi for all $t \in I$. Let $t_0 \in I$, $y(t_0) =
y_0$ and define 
$$A_0(t) = \left | \int_{t_0}^t dt' \ a_0(t')\right |
\quad , \quad A_j(t) = (1 - \alpha_j) \left | \int_{t_0}^t dt' \ a_j
(t')\right | \ , \ 1 \leq j \leq n \ .$$

\noi Then the following inequality holds 
\beq \label{2.7e} y(t) \leq
\exp \left ( A_0(t)\right ) \left \{ y_0^{1 - \alpha_1} + \sum_{1 \leq
j \leq n} A_j(t)^{(1- \alpha_1)/(1 - \alpha_j)}\right \}^{1/(1-
\alpha_1)} \eeq

\noi for all $t\in I$.} \\

\noi {\bf Proof.} We consider only the case $t \geq t_0$. By
exponentiating the term with $a_0$ in (\ref{2.6e}) and by passing to
the rescaled variable $$y' \equiv \exp (- A_0) y$$

\noi we can rewrite (\ref{2.6e}) and the initial condition at $t_0$ in
the form 
$$\left \{ \begin{array}{l} \partial_t y' \leq \displaystyle{\sum_{1 \leq j
\leq n}} a_j \exp \left ( (\alpha_j - 1) A_0 \right ) y'^{\alpha_j}\\ \\
y'(t_0) = y_0\end{array}\right .$$

\noi or equivalently 
$$\left \{ \begin{array}{l} \partial_t y' \leq
\displaystyle{\sum_{1 \leq j \leq n}} (1 - \alpha_j)^{-1} \left ( \partial_t
A'_j\right ) y'^{\alpha_j} \\ \\ y'(t_0) = y_0\end{array}\right .$$

\noi where 
\beq \label{2.8e} A'_j(t) \equiv (1 - \alpha_j) \int_{t_0}^t
dt' \ a_j(t') \exp \left ( (\alpha_j - 1) A_0(t')\right ) \ . \eeq

\noi On the other hand it is well known that $y'(t) \leq z(t)$ where
\beq \label{2.9e} \left \{ \begin{array}{l} \partial_t z = \displaystyle{\sum_{1
\leq j \leq n}} (1 - \alpha_j)^{-1} \left ( \partial_t A'_j\right )
z^{\alpha_j} \\ \\ z(t_0) = y_0\ . \end{array}\right . \eeq

\noi From (\ref{2.9e}) we obtain 
$$\partial_t z\geq (1 - \alpha_j)^{-1}
\left ( \partial_t A'_j \right ) z^{\alpha_j} \ , \qquad 1 \leq j \leq n$$

\noi so that
\beq \label{2.10e} 
z \geq (A'_j)^{1/(1-\alpha_j)} \ , \qquad 1 \leq j \leq n \ .
\eeq

\noi Inserting (\ref{2.10e}) into the differential equation in
(\ref{2.9e}) yields 
$$\partial_t\left ( z^{1 - \alpha_1}\right ) \leq
\sum_{1 \leq j \leq n} (1 - \alpha_1) (1 - \alpha_j)^{-1} \left (
\partial_t A'_j\right ) A{'}_j^{(\alpha_j - \alpha_1)/(1 - \alpha_j)}$$

\noi which implies by integration 
\beq \label{2.11e} z(t)^{1 -
\alpha_1} \leq y_0^{1-\alpha_1} + \sum_{1 \leq j \leq n} \left (
A'_j(t)\right )^{(1- \alpha_1)/(1 - \alpha_j)} \ . \eeq

\noi Now (\ref{2.7e}) follows from $y'(t) \leq z(t)$ and from
(\ref{2.11e}) by coming back to the original variable $y$ and by
replacing the exponential inside the integral of (\ref{2.8e}) by
1.\par\nobreak \hfill $\sq$\\

We now collect some properties of the solutions of the wave equation
$\sq A = 0$. The general solution can be written as 
\beq \label{2.12e}
A(t) = \cos \omega t \ A_+ + \omega^{-1} \sin \omega t \ \dot{A}_+ \ .
\eeq

\noi If $A$ is vector valued and divergence free $(\nabla \cdot A = 0$)
then also $\sq (x \cdot A) = 0$ and $x \cdot A$ can be written as 
\beq
\label{2.13e} x \cdot A(t) = \cos \omega t \left ( x \cdot A_+ \right )
+ \omega^{-1} \sin \omega t \left ( x \cdot \dot{A}_+\right ) \ . \eeq

\noi We shall need the dilation generator 
\beq \label{2.14e} P = t
\partial_t + x \cdot \nabla \ . \eeq

\noi The operator $P$ satisfies the following commutation relations
\bea \label{2.15e} &&\left [ P, \exp (i\omega t) \right ] = 0 \\ &&P
D_0(t) = D_0(t) t \partial_t \\ \label{2.16e} &&P \omega^{-j} =
\omega^{-j} (P + j) \\ \label{2.17e} &&P \partial^{\alpha} =
\partial^{\alpha} (P - |\alpha |) \label{2.18e} \eea

\noi for any integer $j$ and for any multiindex $\alpha$. \par

If $\sq A = 0$ then also $\sq PA = 0$ and $PA$ can be written as 
\beq
\label{2.19e} (PA)(t) = \cos \omega t \left ( x\cdot \nabla A_+\right )
+ \omega^{-1} \sin \omega t \left ( (1 + x \cdot \nabla) \dot{A}_+
\right ) \ . \eeq

\noi When changing variables from $A$ to $B$ according to
(\ref{1.11e}) we obtain 
\beq \label{2.20e} x \cdot A(t) = - t^{-1}
D_0(t) \hbox{\it \v B}(1/t) \eeq

\noi and 
\beq \label{2.21e} \left ( (P + 1)^j A \right ) (t) =
(-1)^{j-1} \ t^{-1} D_0(t) \left ( (t \partial _t)^j B \right ) (1/t)
\eeq

\noi for any integer $j \geq 0$. \par

We finally collect some estimates of divergence free vector solutions
of the wave equation.\\

\noi {\bf Lemma 2.4.} {\it Let $j \geq 0$ be an integer and let
$\alpha$ be a multiindex. Assume that $(A_+, \dot{A}_+)$ satisfy the
conditions
\bea \label{2.22e} &{\cal A} \in L^2 \quad , &\quad \nabla^2{\cal A} \in L^1 \\
&\omega^{-1} \dot{\cal A}\in L^2 \quad , &\quad \nabla \dot{\cal A} \in L^1
\label{2.23e} \eea

\noi for 
\bea \label{2.24e} &{\cal A} = (x \cdot \nabla )^{j'}
\partial_x^{\alpha} A_+ \quad , &\quad {\cal A} = (x \cdot \nabla)^{j'}
\partial_x^{\alpha} \left ( x \cdot A_+ \right ) \\ &\dot{\cal A} = (x \cdot
\nabla )^{j'} \partial_x^{\alpha} \dot{A}_+ \quad , &\quad \dot{\cal A} = (x
\cdot \nabla)^{j'} \partial_x^{\alpha} \left ( x \cdot \dot{A}_+ \right
) \eea

\noi for $0 \leq j' \leq j$. Then $A$ satisfies the following estimates~:
 \beq
\label{2.26e} \parallel (P+1)^{j'} \ \partial_x^{\alpha} A(t)
\parallel_r \ \vee \ \parallel (P+1)^{j'} \partial_x^{\alpha} (x \cdot
A(t)) \parallel_r \ \leq b \ t^{-1+2/r} \eeq

\noi for $0 \leq j' \leq j$, for $2 \leq r \leq \infty$ and for all $t
> 0$.\par 

Let $B$ and $\hbox{\it \v B}$ be defined by (\ref{1.11e}) and
(\ref{1.13e}). Then $B$ and $\hbox{\it \v B}$ satisfy the following
estimates~:}

\beq \label{2.27e} \parallel \partial_t^j \ \partial_x^{\alpha} \
B(t)\parallel_r\ \vee \ \parallel \partial_t^j \ \partial_x^{\alpha} \
\hbox{\it \v B}(t) \parallel_r\  \leq b \ t^{-j-|\alpha | + 1/r} \ . \eeq

\noi {\bf Proof.} The estimate (\ref{2.26e}) is standard for $j = 0$
\cite{14r}. For $j \not= 0$ it is a consequence of the case $j = 0$ and
of the commutation relations satisfied by $P$. The estimate
(\ref{2.27e}) follows from (\ref{2.26e}), (\ref{1.11e}), (\ref{2.20e})
and (\ref{2.21e}). \par\nobreak \hfill $\sq$

\mysection{The Cauchy problem at finite time} \hspace*{\parindent} In
this section we study the Cauchy problem for the equations (\ref{1.1e})
and (\ref{1.12e}) and for related nonautonomous equations that will
appear in Sections 5 and 6, for finite initial time. We shall write
those equations in the general form \beq \label{3.1e} i \partial_tv =
Kv + f \eeq

\noi where \beq \label{3.2e} K = - (1/2) \Delta_A + V \ . \eeq

\noi The equations (\ref{1.1e}) and (\ref{1.12e}) are obtained by
replacing $(v, A, V, f)$ by $(u, A, 0, 0)$ and by $(w, B, -\hbox{\it \v B}, 0)$
respectively. We shall need a parabolically regularized version of
(\ref{3.1e}) of the form 
\beq \label{3.3e} i\partial_t v = K_{\eta} v +
f \eeq

\noi where \beq \label{3.4e} K_{\eta} = - (1/2) (1 - i \eta ) \Delta_A
+ V \eeq

\noi for some $\eta > 0$. Furthermore we shall regularize $A$, $V$ and
$f$.\par

We shall use extensively the conservation law and/or estimates of the
$L^2$ norm which are formally associated with (\ref{3.1e}) and
(\ref{3.3e}) and which actually hold for sufficiently regular $A$, $V$
and $f$. We state them in the following two Lemmas.\\

\noi {\bf Lemma 3.1.} {\it Let $I$ be an interval, let $A \in
L_{loc}^2(I, L^{\infty})$, $\nabla A \in L_{loc}^1(I, L^{\infty})$, $V
\in L_{loc}^1 (I, L^{\infty})$ and $f \in L_{loc}^1(I, L^2)$. Let $v
\in (L_{loc}^{\infty} \cap {\cal C}_w) (I, L^2)$.\par

(1) Let $\eta > 0$ and let $v$ satisfy (\ref{3.3e}) in $I$. Then $v \in
L_{loc}^2(I,H^1)$ and $v$ satisfies the following inequality for all $t_1,
t_2 \in I$, $t_1 \leq t_2$~: 
$$\parallel v(t_2)\parallel_2^2 \ - \ \parallel
v(t_1)\parallel_2^2 \ + \ \eta \int_{t_1}^{t_2}  dt \parallel \nabla_A v(t) \parallel_2^2$$ 
\beq \label{3.5e} \leq \int_{t_1}^{t_2} dt
\ 2 \ {\rm Im} \ <v,f>(t) \ . \eeq

(2) Let $v$ satisfy (\ref{3.1e}) in $I$. Then $v \in {\cal C}(I, L^2)$
and $v$ satisfies the following equality for all $t_1$, $t_2$ in $I$~:}
\beq \label{3.6e} \parallel v(t_2)\parallel_2^2 \ - \ \parallel
v(t_1)\parallel_2^2 \ = \int_{t_1}^{t_2}  dt \ 2 \ {\rm Im}\ <v,f>(t) \ .
\eeq

There is a large number of results of the type of Lemma 3.1 in the
literature. The proof relies on commutator estimates of vector fields
with standard mollifiers. For instance for the proof of (\ref{3.6e}),
one starts from the regularized identity \bea \label{3.7e} &&\parallel
\varphi * v(t_2)\parallel_2^2 \ - \ \parallel \varphi *
v(t_1)\parallel_2^2 \ = 2 \ {\rm Re} \int_{t_1}^{t_2} dt <\varphi * v,
[\varphi *, A \cdot \nabla ]v>(t) \nn \\ &&+ \ 2 \ {\rm Im}
\int_{t_1}^{t_2} dt \left ( <\varphi * v, [\varphi *, A^2/2 + V ]v> +
<\varphi * v, \varphi * f> \right ) (t)\eea

\noi where $\varphi$ is a standard mollifier in the space variables.
The only delicate term in the elimination of the mollifier is that
containing the commutator $[\varphi *, A \cdot \nabla]$. We refer to
\cite{2r} \cite{3r} for detailed proofs. \par

If $v$ is more regular, similar results hold under weaker assumptions
on $A$, $V$, $f$. We state only one such result related to
(\ref{3.1e}). \\

\noi {\bf Lemma 3.2.} {\it Let $I$ be an interval. Let $A \in
L_{loc}^2(I, L^3 + L^{\infty})$, $V \in L_{loc}^1(I, L^{3/2} +
L^{\infty})$ and $f \in L_{loc}^1(I, H^{-1})$. Let $v \in
(L_{loc}^{\infty} \cap {\cal C}_w)(I, H^1)$ satisfy (\ref{3.1e}) in $I$.
Then $v$ satisfies (\ref{3.6e}) for all $t_1, t_2 \in I$.}\\

The proof of Lemma 3.2 is an elementary variant of that of Lemma 3.1.
One starts again from (\ref{3.7e}), where now the elimination of the
mollifier is elementary under the assumptions made on $A$, $V$ and $f$. \\

\noi {\bf Remark 3.1.} The assumption $v \in (L_{loc}^{\infty} \cap
{\cal C}_w)(I,L^2)$ in Lemma 3.1 is partly redundant. On the one hand,
if one is given a solution $v \in L_{loc}^{\infty}(I, L^2)$ of
(\ref{3.1e}) or (\ref{3.3e}) in $I$, that solution has a representative
$v' \in {\cal C}(I, H^{-2})$. By a standard compactness argument $v'
\in {\cal C}_w(I, L^2)$ and therefore $v'$ satisfies the assumptions of
the Lemma (see for instance \cite{9r}). On the other hand, if $v \in {\cal C}_w(I,L^2)$ then by the
uniform boundedness principle, $v \in L_{loc}^{\infty}(I,L^2)$. In
practice, we shall use Lemma 3.1 mostly for $\eta > 0$ and solutions $v
\in {\cal C}(I, L^2)$ or for weak$*$ limits of such solutions in
$L^{\infty}(I, L^2)$ which will be extended to ${\cal C}_w(I, L^2)$ as
explained above.\\

We now begin the study of the Cauchy problem for the equation
(\ref{3.1e}). As mentioned in the introduction, we shall study that
problem successively at the level of regularity of $L^2$, $H^2$, $H^3$
and $H^4$ using time derivatives as much as possible to control the
regularity. \par

We first state the result at the level of $L^2$.\\

\noi {\bf Proposition 3.1.} {\it Let $I$ be an interval, let $A \in
L_{loc}^2(I, L^4 + L^{\infty})$, $V \in L_{loc}^1(I, L^2 + L^{\infty})$
and $f \in L_{loc}^1(I, L^2)$. Let $t_0 \in I$ and $v_0 \in L^2$. Then
\par

(1) There exists a solution $v \in (L_{loc}^{\infty} \cap {\cal
C}_w)(I, L^2)$ of (\ref{3.1e}) in $I$ with $v(t_0) = v_0$. That
solution satisfies 
\beq \label{3.8e} \parallel v(t) \parallel_2^2 -
\parallel v_0 \parallel_2^2 \ \leq \int_{t_0}^t dt'\ 2\ {\rm Im}
<v,f>(t') \eeq

\noi for all $t \in I$.\par

(2) Let in addition $A \in L_{loc}^2(I, L^{\infty})$, $V \in
L_{loc}^1(I, L^{\infty})$ and $\nabla A \in L_{loc}^1(I, L^{\infty})$.
Then the previous solution $v$ is unique in $(L_{loc}^{\infty} \cap
{\cal C}_w)(I, L^2)$, $v \in {\cal C}(I, L^2)$ and $v$ satisfies
(\ref{3.6e}) for all $t_1, t_2 \in I$.}\\

\noi {\bf Proof.} \par

\noi \underbar{Part (1)}. The proof proceeds by a parabolic
regularization and a limiting procedure. We consider separately the
cases $t \geq t_0$ and $t \leq t_0$ and we begin with $t \geq t_0$. We
replace (\ref{3.1e}) by (\ref{3.3e}) with $0 < \eta \leq 1$, where in
addition we regularize $A$ in space and time and $V$ in space by the
use of standard mollifiers parametrized by $\eta$, so that the
regularized $A$ and $V$ belong to ${\cal C}^{\infty}(I, H_{\infty}^N)$
and to $L^1(I, H_{\infty}^N)$ respectively for any $N \geq 0$. We
recast the Cauchy problem for the regularized equation in the form of
the integral equation \beq \label{3.9e} v(t) = (\phi (v))(t) \equiv
U_{\eta}(t-t_0) v_0 - i \int_{t_0}^t dt' \ U_{\eta} (t - t')
F_{\eta}(t') \eeq

\noi where 
\bea \label{3.10e} &&U_{\eta}(t) = \exp ( i(t/2) (1 - i
\eta ) \Delta ) \ , \\ &&F_{\eta} = (1 - i \eta ) \left ( iA \cdot
\nabla v + (1/2) A^2 v \right ) + Vv + f \ . \label{3.11e} \eea

We first solve (\ref{3.9e}) locally in time by contraction in ${\cal
C}([t_0, t_0 + T], L^2)$ for some $T > 0$. The semigroup $U_{\eta}$
satisfies the estimate 
\beq \label{3.12e} \parallel U_{\eta} (t) \
\partial_x^{\alpha} ; {\cal B}(L^2)\parallel \ \leq C_{\alpha} (\eta
t)^{-|\alpha |/2} \qquad \hbox{for $t > 0$} \eeq

\noi so that 
$$\parallel U_{\eta} (t- t')\
F_{\eta}(t')\parallel_2 ÷ \leq C \left \{ (\eta (t-t'))^{-1/2}
\parallel A(t')\parallel_{\infty} \ + \ \parallel
A(t')\parallel_{\infty}^2 \ + \ \parallel V(t')\parallel_{\infty}
\right \}$$
\beq \label{3.13e} 
 \times \ \parallel v(t')\parallel_2 ÷ + \ \parallel
f(t') \parallel_2 \ . \eeq

\noi The RHS of (\ref{3.13e}) is in $L^1$ in the variable $t'$. It then
follows from (\ref{3.13e}) that (\ref{3.9e}) can be solved by
contraction in ${\cal C}([t_0, t_0 + T], L^2)$ for $T$ sufficiently
small. By a standard argument using the fact that the equation is
linear, one can extend the solution to $I_+ = I \cap \{t:t\geq t_0\}$.
Let $v_{\eta}$ be that solution.\par

We next take the limit where $\eta$ tends to zero. For that purpose we
use the fact that $v_{\eta}$ satisfies
 \beq \label{3.14e} \parallel
v_{\eta}(t)\parallel_2^2 \ \leq \ \parallel v_0\parallel_2^2 \ +
\int_{t_0}^t dt'\ 2 \ {\rm Im} <v_{\eta}, f> (t') \eeq

\noi for all $t \in I_+$ by Lemma 3.1, part (1), so that \beq
\label{3.15e} \parallel v_{\eta}(t)\parallel_2 \ \leq \ \parallel v_0
\parallel_2 \ + \int_{t_0}^t dt'\parallel f(t')\parallel_2 \ , \eeq

\noi and $v_{\eta}$ is uniformly bounded with respect to $\eta$ in
$L_{loc}^{\infty}(I_+, L^2)$. We take the limit $\eta \to 0$ by a
compactness argument. We can extract a subsequence of $v_{\eta}$ which
converges to some $v \in L_{loc}^{\infty}(I_+, L^2)$ in the
weak$*$ sense in $L^{\infty}(J, L^2)$ for all $J \subset\subset I_+$.
One can see that $v$ satisfies the equation (\ref{3.1e}). In particular
$v \in {\cal C}(I_+, H^{-2})$ and $v$ can therefore be taken in ${\cal
C}_w (I_+, L^2)$ (see Remark 2.1). Furthermore $v_{\eta}$
converges to $v$ weakly in $L^2$ pointwise in $t$. This can be seen
easily from the identity 
\beq \label{3.16e} v_{\eta}(t) = \theta^{-1}
\int_t^{t+ \theta} dt'\ v_{\eta}(t') - \theta^{-1} \int_t^{t+\theta}
dt' (t + \theta - t') \partial_t \ v_{\eta}(t') \eeq

\noi and from the previous convergence of $v_{\eta}$ to $v$. This allows
to prove that $v$ satisfies the initial condition $v(t_0) = v_0$ and to
take the limit $\eta \to 0$ in (\ref{3.14e}), thereby obtaining
(\ref{3.8e}). This completes the proof for $t \geq t_0$. A similar
proof applies to the case $t \leq t_0$.\\

\noindent \underbar{Part (2)}. It follows from Lemma 3.1, part (2) that
$v$ satisfies (\ref{3.6e}) and that $v \in {\cal C}(I, L^2)$.
Uniqueness follows from (\ref{3.6e}) with $f = 0$ applied to the
difference of two solutions.\par\nobreak \hfill $\sq$\par

We now turn to the study of the Cauchy problem for the equation
(\ref{3.1e}) at the level of regularity of $H^2$. As mentioned
previously, we shall control that regularity through the use of the
time derivative $\partial_tv$. That derivative satisfies the equation
obtained by taking the time derivative of (\ref{3.1e}), namely \beq
\label{3.17e} i \partial_t \ \partial_tv = K \partial_t v + f_1 \eeq

\noi \beq \label{3.18e} f_1 = \left ( \partial_t K \right ) v +
\partial_t f = i \left ( \partial_tA\right ) \cdot \nabla_Av + \left (
\partial_tV\right ) v + \partial_t f \ . \eeq

We state the result as Proposition 3.2 below. That proposition is a
minor extension of Lemma 3.1 in \cite{10r} since we have in addition a
linear potential $V$ and an inhomogeneous term $f$. In \cite{10r} the
authors use the covariant space derivatives, which is equivalent to
using the time derivative for $V = f = 0$, but not for nonzero $V$ and
$f$.\\

\noi {\bf Proposition 3.2.} {\it Let $I$ be an interval, let $A \in
L_{loc}^{\infty}(I, L^6 + L^{\infty}) \cap {\cal C}(I, L^3 +
L^{\infty})$, $\partial_t A \in L_{loc}^1(I, L^3 + L^{\infty})$, $V \in
{\cal C}(I, L^2 + L^{\infty})$, $\partial_t V \in L_{loc}^1(I, L^2 +
L^{\infty})$, $f \in {\cal C}(I, L^2)$ and $\partial_tf \in L^1(I,
L^2)$. Let $t_0 \in I$ and $v_0 \in H^2$. Then\par

(1) There exists a unique solution $v \in {\cal C}(I, H^2) \cap {\cal
C}^1(I, L^2)$ of (\ref{3.1e}) in $I$ with $v(t_0) = v_0$. That solution
satisfies (\ref{3.6e}) for all $t_1, t_2 \in I$. That solution is
actually unique in $(L_{loc}^{\infty} \cap {\cal C}_w)(I, H^1)$. \par

(2) Let in addition $A \in L_{loc}^2(I, L^{\infty})$, $\nabla A \in
L_{loc}^1(I, L^{\infty})$ and $V \in L_{loc}^1(I, L^{\infty})$. Then
the previous solution $v$ is actually unique in $(L_{loc}^{\infty} \cap
{\cal C}_w)(I, L^2)$ and $v$ satisfies 
\beq \label{3.19e} \parallel
\partial_t v(t_2) \parallel_2^2 \ - \ \parallel \partial_t
v(t_1)\parallel_2^2 \ = \int_{t_1}^{t_2} dt\ 2 \ {\rm Im} <\partial_t v,
f_1> (t) \eeq

\noi for all $t_1, t_2 \in I$, where $f_1$ is defined by
(\ref{3.18e}).}\\

\noi {\bf Remark 3.2.} Assumptions of the type $\partial_t A \in
L_{loc}^1(I, X)$ almost imply that $A \in {\cal C}(I, X)$. The latter
condition serves simply to exclude that $A$ contains a constant term in
time which does not belong to $X$. A similar remark will apply to more
complicated assumptions of the same type made in Propositions 3.3 and
3.4 below.\\

\noi {\bf Proof.} \underbar{Part (1)}. The proof proceeds by a
parabolic regularization and a limiting procedure as in the case of
Proposition 3.1. We consider first the case $t \geq t_0$. We replace
(\ref{3.1e}) by (\ref{3.3e}) with $0 \leq \eta \leq 1$, where in
addition we regularize $A$, $V$ and $f$ by standard mollifiers in space
parametrized by $\eta$, in such a way that the regularization decreases the
relevant $L^r$ norms. We shall not indicate the regularization in the
notation for $A$ and $V$. As regards $f$, we shall in general not
indicate it either, except in cases of doubt where we shall use the
notation $f_{\eta}$. The regularized $A$, $V$ and $f$ satisfy
conditions obtained from the assumptions of the proposition by
replacing $L^r + L^{\infty}$ by $H_{\infty}^N$ for $A$ and $V$ and $L^2$
by $H^N$ for $f$, for arbitrary $N \geq 0$. We consider again the
integral equation, namely (\ref{3.9e}). We shall need the time
derivative of $\phi (v)$. Integration by parts yields
\beq
\label{3.20e}
\left ( i \partial_t \phi (v)\right ) (t) = U_{\eta} \left ( t-t_0 \right ) \left ( K_{\eta} v_0 + f \right ) (t_0) + \int_{t_0}^t dt'\ U_{\eta} (t- t') \partial_t  F_{\eta}(t')
\eeq

\noi while (\ref{3.11e}) yields
\beq
\label{3.21e}
\partial_t F_{\eta} = (1 - i \eta) \left ( i A \cdot \nabla \partial_t v + (1/2) A^2 \partial_t v \right ) + V \partial_t v + f_{1 \eta}
\eeq

\noi where \beq \label{3.22e} f_{1 \eta} = \left ( \partial_t K_{\eta}
\right ) v + \partial_tf = i (1 - i \eta ) \left ( \partial_t A \right
) \cdot \nabla_A v + \left ( \partial_t V\right ) v + \partial_t f \ .
\eeq

\noi We first solve (\ref{3.9e}) locally in time by contraction in
${\cal C}([t_0,t_0 + T], H^2) \cap$\break \noindent  ${\cal C}^1([t_0,t_0 + T], L^2)$ for
some $T > 0$. For that purpose we estimate 
\bea \label{3.23e} \parallel
U_{\eta}(t-t')F_{\eta}(t')\parallel_2 &\leq & \Big \{ \parallel
A\parallel_{\infty} \ \parallel \nabla v\parallel_2 \ + \left (
\parallel A \parallel_{\infty}^2 \ + \ \parallel V \parallel_{\infty}\right ) \parallel v \parallel_2 \ \nn \\
&&+ \ \parallel f \parallel_2 \Big \} (t') \eea 
$$\parallel
U_{\eta}(t-t')\Delta F_{\eta}(t')\parallel_2 \ \leq C (\eta
(t-t'))^{-1/2} \parallel \nabla \left ( F_{\eta} - f\right ) (t')
\parallel_2 \ + \ \parallel \Delta f(t')\parallel_2$$ 
$$\leq C
(\eta (t- t'))^{-1/2} \Big \{ \parallel A \parallel_{\infty} \
\parallel \Delta v \parallel_2 \ + \left ( \parallel \nabla A
\parallel_{\infty}\ + \ \parallel A \parallel_{\infty}^2 \ + \
\parallel V \parallel_{\infty} \right ) \parallel \nabla v \parallel_2$$
\beq \label{3.24e} 
+ \left ( \parallel A \parallel_{\infty} \ \parallel \nabla A
\parallel_{\infty} \ + \ \parallel \nabla V\parallel_{\infty} \right )
\parallel v \parallel_2 \Big \} (t') + \ \parallel \Delta
f(t')\parallel_2 \eeq 
\bea \label{3.25e} &&\parallel
U_{\eta}(t-t')\partial_t F_{\eta}(t')\parallel_2 \ \leq C (\eta
(t-t'))^{-1/2} \parallel A(t')\parallel_{\infty} \ \parallel \partial_t
v(t') \parallel_2\nn \\ &&+ \ C \Big \{ \left ( \parallel A
\parallel_{\infty}^2 \ + \ \parallel V \parallel_{\infty} \right )
\parallel \partial_tv\parallel_2\ + \ \parallel \partial_t A
\parallel_{\infty} \ \parallel \nabla v \parallel_2   \nn \\ &&+ \left (
\parallel A \parallel_{\infty} \ \parallel \partial_t A
\parallel_{\infty} \ + \ \parallel \partial_t V\parallel_{\infty}
\right ) \parallel v \parallel_2 \Big \} (t') + \ \parallel \partial_t
f(t')\parallel_2 \ . \eea

\noi The RHS of (\ref{3.23e})-(\ref{3.25e}) are in $L^1$ of the variable
$t'$. By the same argument as in Proposition 3.1, one obtains a
solution $v_{\eta} \in {\cal C}(I_+, H^2) \cap {\cal C}^1(I_+, L^2)$
with $I_+ = I \cap \{t:t \geq t_0 \}$. \par

We next take the limit where $\eta$ tends to zero and for that purpose
we need estimates of $v_{\eta}$ uniform in $\eta$ in the relevant
space. We remark that $v_{\eta}$ satisfies the equation (compare
(\ref{3.17e}) (\ref{3.18e}) with (\ref{3.22e})) 
\beq
\label{3.26e}
i \partial_t \ \partial_t  v_{\eta} = K_{\eta} \partial_t v_{\eta} + f_{1 \eta} \ .
\eeq

\noi Let $y_0 = \parallel v_{\eta}\parallel_2$ and $y_1 = \parallel
\partial_t v_{\eta}\parallel_2$ . By Lemma 3.1, part (1), $\partial_t
v_{\eta}$ satisfies 
\beq \label{3.27e} y_1(t_2)^2 - y_1(t_1)^2 \leq \int_{t_1}^{t_2} dt \ 2\ 
{\rm Im} <\partial_t v_{\eta}, f_{1 \eta}>(t) \eeq

\noi for all $t_1, t_2 \in I_+$, $t_1 \leq t_2$, so that 
\beq
\label{3.28e}
\partial_t y_1 \leq \ \parallel f_{1 \eta}\parallel_2 \ .
\eeq

\noi We have already estimated $y_0$ by (\ref{3.15e}) uniformly in
$\eta$. For brevity we continue the estimates by omitting the index
$\eta$ in $v_{\eta}$. Furthermore, we keep only the most dangerous,
namely the most singular parts of $A$ and $V$ by keeping only the $L^p$
component in all the $L^p + L^{\infty}$ spaces that occur in the
assumptions of the proposition. The more regular $L^{\infty}$
components yield contributions that can be estimated similarly and
involve lower norms of $v$. From (\ref{3.22e}) we obtain
\bea
\label{3.29e} \parallel f_{1 \eta}\parallel_2 & \leq &(1 + \eta^2)^{1/2}
\parallel \partial_t A \parallel_3 \left ( \parallel \nabla
v\parallel_6 \ + \ \parallel A\parallel_6 \ \parallel v
\parallel_{\infty} \right ) + \ \parallel \partial_t V \parallel_2\ 
\parallel v \parallel_{\infty} \nn \\
&& + \ \parallel \partial_t f \parallel_2
\ . \eea

\noi On the other hand from (\ref{3.3e}) we obtain by direct estimation 
\beq \label{3.30e} \parallel \Delta v \parallel_2 \ \leq 2 \left \{ y_1
+ \ \parallel A \parallel_6\ \parallel \nabla v \parallel_3 \ + \
\parallel A \parallel_6^2 \ \parallel v \parallel_6 \ + \ \parallel V
\parallel_2 \ \parallel v \parallel_{\infty} \ + \ \parallel f
\parallel_2 \right \} \ . \eeq

\noi By Lemma 2.1, we estimate
\beq \label{3.31e} \left \{ \begin{array}{l} \parallel \nabla v
\parallel_3 \ \vee \ \parallel v \parallel_{\infty} \ \leq C \parallel
v \parallel_2^{1/4} \parallel \Delta v \parallel_2^{3/4} \ , \\ \\
\parallel v \parallel_6 \ \leq C \parallel v \parallel_2^{1/2}\ 
\parallel \Delta v \parallel_2^{1/2}  .\end{array}\right . \eeq

\noi Substituting (\ref{3.31e}) into (\ref{3.30e}) and using Lemma 2.2, we obtain
\beq \label{3.32e} \parallel \Delta v \parallel_2 \ \leq C \left ( y_1
+ \ \parallel f \parallel_2 \ + \ m^4y_0 \right ) \eeq

\noi where
\beq
\label{3.34e}
m = \ \parallel A \parallel_6 \ + \ \parallel V \parallel_2 \ .
\eeq

\noi From (\ref{3.29e}) and Lemma 2.1, we obtain 
$$\parallel f_{1\eta} \parallel_2 \ \leq C \left \{ \parallel
\partial_tA\parallel_3 \ \parallel \Delta v\parallel_2 \ + \left (
\parallel \partial_t A \parallel_3 \ \parallel A \parallel_6 \ + \
\parallel \partial_t V \parallel_2 \right ) y_0^{1/4} \parallel \Delta
v \parallel_2^{3/4} \right \}$$ \beq \label{3.35e} +\ \parallel
\partial_t f \parallel_2 \ \equiv M_0 \left ( \parallel \Delta
v\parallel_2\ , y_0\right ) + \ \parallel \partial_t f\parallel_2 \eeq

\noi so that by (\ref{3.28e}) (\ref{3.32e}) 
\beq \label{3.36e} \partial_t y_1 \leq \ \parallel f_{1\eta}\parallel_2
\ \leq M\left ( y_1 \ + \ \parallel f\parallel_2\ , y_0 \right ) + \
\parallel \partial_t f \parallel_2 \eeq

\noi where
\bea \label{3.37e} M(z, y_0) &=& C \Big \{ \parallel \partial_t A
\parallel_3 \left ( z + m^4y_0\right ) + \left ( \parallel \partial_t A
\parallel_3 \ \parallel A \parallel_6 \ + \ \parallel \partial_t V
\parallel_2 \right ) \nn \\ &&\times \left ( y_0^{1/4} z^{3/4} + m^3
y_0 \right ) \Big \} \ . \eea

\noi It then follows from Lemma 2.3 that $y_1$ is estimated by
\bea \label{3.38e}
&&y_1(t) \leq \exp \left ( C \int_{t_0}^t dt'
\parallel \partial_t A(t')\parallel_3 \right ) \Big \{ y_1(t_0)^{1/4} +
\int_{t_0}^t dt' \left ( \parallel \partial_t V \parallel_2
\bar{y}_0^{1/4} \right ) (t') \nn \\ &&+ \Big \{ \int_{t_0}^tdt'\Big (
\parallel \partial_t A \parallel_3 (\parallel f\parallel_2 + \ m^4
\bar{y}_0 ) + \ \parallel \partial_t V \parallel_2 \left ( \parallel f
\parallel_2^{3/4}\bar{y}_0^{1/4} + m^3 \bar{y}_0 \right ) \nn \\ &&+
\parallel \partial_t f \parallel_2 \Big ) \Big \}^{1/4}\Big \}^4\eea

\noi where $\bar{y}_0$ is an estimate of $y_0$ uniform in $\eta$ as
obtained previously from (\ref{3.15e}). Substituting (\ref{3.38e}) into
(\ref{3.32e}) finally yields an a priori estimate of $\parallel \Delta
v\parallel_2$. The estimates thereby obtained are uniform in $\eta$. In
fact 
$$y_1(t_0) = \ \parallel (K_{\eta} v_0 +
f)(t_0)\parallel_2 $$
\beq \label{3.39e} \leq \ \parallel \Delta v_0 \parallel_2 \ +
2 \parallel A\parallel_3 \ \parallel \nabla v_0 \parallel_6 \ + \left (
\parallel A \parallel_4^2 \ + \ \parallel V \parallel_2 \right )
\parallel v_0\parallel_{\infty} \ + \ \parallel f\parallel_2 \ . \eeq

\noi The norms of $A$, $V$ and $f$ that occur in (\ref{3.32e}),
(\ref{3.38e}) and (\ref{3.39e}) are controlled by the assumptions of
the Proposition, so that the norms of the regularized quantities are
bounded uniformly with respect to the regularization, since the
regularisation is taken such as to decrease those norms.\par

We can now take the limit $\eta \to 0$. The solution $v_{\eta}$ is
uniformly bounded in $L_{loc}^{\infty}(I_+, H^2) \cap H_{\infty loc}^1
(I_+, L^2)$. By compactness  we can extract a subsequence which
converges in the weak$*$ sense to some $v \in L_{loc}^{\infty}(I_+,
H^2) \cap H_{\infty loc}^1 (I_+, L^2)$. One can see that $v$ satisfies
(\ref{3.1e}) and therefore can be chosen in ${\cal C}_w(I_+, H^2) \cap
{\cal C}_w^1(I_+, L^2)$. Furthermore $v_{\eta}$ converges pointwise to
$v$ weakly in $H^2$ and $\partial_t v_{\eta}$ converges pointwise to
$\partial_t v$ weakly in $L^2$. Together with the fact that $(K_{\eta}
v_0 + f_{\eta})(t_0)$ converges to $(Kv_0 + f) (t_0)$ strongly in
$L^2$, this allows to prove that $v$ satisfies the initial conditions
$v(t_0) = v_0$ and $i\partial_tv(t_0) = (Kv_0 + f)(t_0)$. By Lemma 3.2,
the solution $v$ satisfies (\ref{3.6e}) and therefore is unique in
$(L_{loc}^{\infty} \cap {\cal C}_w)(I_+, H^1)$.\par

We now turn to the strong continuity of $\partial_tv$ in $L^2$ and of
$v$ in $H^2$. When $\eta \to 0$, $y_1(t_0)$ converges to its value for
$\eta = 0$ so that in the limit $\eta \to 0$ the RHS of (\ref{3.38e})
is bounded by its value for $\eta = 0$. On the other hand $y_1(t)$ is
non increasing in that limit under pointwise weak convergence of
$\partial_tv_{\eta}$ in $L^2$. Therefore (\ref{3.38e}) also holds in
the limit $\eta \to 0$. Since $\partial_tv$ is weakly continuous in
$L^2$ and since the RHS of (\ref{3.38e}) tends to $y_1(t_0)$ when $t$
decreases to $t_0$, $\partial_t v$ is strongly continuous from the
right at $t_0$. \par

This completes the proof for $t \geq t_0$, except for strong
continuity. A similar proof yields the corresponding results for $t
\leq t_0$. In particular it yields strong continuity of $\partial_tv$
from the left at $t_0$, which together with the previous result yields
strong continuity of $\partial_tv$ at $t_0$. Strong continuity of
$\partial_tv$ for any $t \in I$ now follows from strong continuity at
$t_0$ and from uniqueness by varying $t_0$ for a given solution $v$.
\par

Finally strong continuity of $v$ in $H^2$ follows from the strong
continuity of $\partial_tv$ and from (\ref{3.1e}) under the available
continuity assumptions on $A$, $V$ and $f$.\\

\noindent \underbar{Part (2)}. Uniqueness is a rewriting of the
uniqueness part of Proposition 3.1, part (2). The equality
(\ref{3.19e}) follows from Lemma 3.1, part (2) applied to $\partial_tv$
as a solution of (\ref{3.17e}).\par\nobreak \hfill $\sq$ \par

We now turn to the study of the Cauchy problem for (\ref{3.1e}) at the
level of regularity of $H^3$. In addition to the space derivatives of
$v$ and to the time derivative $\partial_tv$ which satisfies
(\ref{3.17e}) (\ref{3.18e}), we shall use the mixed space time
derivative $\nabla_A \partial_t v$. That derivative satisfies an
evolution equation obtained by applying $\nabla_A$ to (\ref{3.17e}),
namely

\beq \label{3.40e} i \partial_t \nabla_A \partial_tv = - (1/2) \nabla_A
\Delta_A \partial_t v + V \nabla_A \partial_tv + g \eeq

\noi where
\bea \label{3.41e} g&=& \left ( \partial_t A + \nabla V\right )
\partial_t v + \nabla_A f_1 \nn \\ &=& \left ( \partial_t A + \nabla
V\right ) \partial_t v + i ( \partial_t A) \cdot \nabla_A \otimes
\nabla_A v + i ( \nabla \partial_t A ) \cdot
\nabla_A v \nn \\ &&+ \left ( \partial_t V\right ) \nabla_A v + \left
( \nabla \partial_t V\right ) v + \nabla_A \partial_t f \ . \eea

The result can be stated as follows. \\

\noi {\bf Proposition 3.3.} {\it Let $I$ be an interval. Let $A$, $V$
and $f$ satisfy $A \in {\cal C}(I, L^s + L^{\infty})$, $\partial_t A
\in L_{loc}^1(I, L^s + L^{\infty})$, $\nabla A \in {\cal C}(I, L^r +
L^{\infty})$, $\partial_t\nabla A \in L_{loc}^1(I, L^r + L^{\infty})$
for some $r$, $s$, with $1/r + 1/s = 1/2$, $2 \leq r \leq 3$, $V \in
{\cal C}(I, L^6 + L^{\infty})$, $\partial_t V \in L_{loc}^1(I, L^6 +
L^{\infty})$, $\nabla V \in {\cal C}(I, L^2 + L^{\infty})$, $\partial_t
\nabla V \in L_{loc}^1 (I, L^2 + L^{\infty})$, $f \in {\cal C}(I,H^1)$,
$\partial_t f \in L_{loc}^1 (I, H^1)$. Let $t_0 \in I$ and $v_0 \in
H^3$. Then \par

(1) There exists a unique solution $v \in {\cal C} (I, H^3) \cap {\cal
C}^1(I, H^1)$ of (\ref{3.1e}) in $I$ with $v(t_0) = v_0$. That solution
satisfies (\ref{3.6e}) and (\ref{3.19e}) for all $t_1,t_2 \in I$. That
solution is actually unique in $(L_{loc}^{\infty} \cap {\cal C}_w)(I,
H^1)$. \par

(2) Let in addition $A \in L_{loc}^2(I, L^{\infty})$ , $\nabla A \in
L_{loc}^1(I, L^{\infty})$ and $V \in L_{loc}^1(I, L^{\infty})$. Then the
previous solution $v$ is actually unique in $(L_{loc}^{\infty} \cap
{\cal C}_w)(I, L^2)$ and $v$ satisfies
\beq \label{3.42e} \parallel \nabla_A \partial_tv(t_2)\parallel_2^2 \ -
\ \parallel \nabla_A \partial_t v(t_1)\parallel_2^2 \ =
\int_{t_1}^{t_2} dt \ 2 \ {\rm Im} < \nabla_A \partial_t v, g>(t) \eeq

\noi for all $t_1, t_2 \in I$, where $g$ is defined by
(\ref{3.41e}).}\\

\noi {\bf Remark 3.3.} We recall that assumptions on the time
derivative of a function and on that function are related as explained
in Remark 2.2. On the other hand the assumptions on $(A, \partial_t A)$
follow from those on $(\nabla A, \partial_t \nabla A)$ by Lemma 2.1 for
$r \geq 12/5$, $s \leq 12$. We have written both assumptions explicitly
in order to avoid that restriction. Similarly the assumptions on $(V,
\partial_t V)$ follow from those on $(\nabla V, \partial_t \nabla
V)$.\\

\noi {\bf Proof.} The proof proceeds by a parabolic regularization and
a limiting procedure as in the case of Proposition 3.2. We consider
first the case $t \geq t_0$. We replace (\ref{3.1e}) by (\ref{3.3e})
with $0 < \eta \leq 1$ and with $A$, $V$ and $f$ regularized in space
as in Proposition 3.2. We use again the integral equation for $v$,
namely (\ref{3.9e}), and we solve that equation locally in time by
contraction in ${\cal C}([t_0, t_0 + T], H^2) \cap {\cal
C}^1([t_0,t_0+T], H^1)$ for some $T >0$. For that purpose we use again
the estimates (\ref{3.23e}) (\ref{3.24e}) (\ref{3.25e}), supplemented
by an additional estimate for $U_{\eta}(t-t') \nabla \partial_t
F_{\eta}(t')$, namely 
\bea \label{3.43e} &&\parallel U_{\eta}(t-t') \nabla \partial_t
F_{\eta} (t')\parallel_2 \ \leq C ( \eta (t-t'))^{-1/2} \parallel
A(t')\parallel_{\infty} \ \parallel \nabla \partial_t v(t')\parallel_2 \nn\\ 
&&+ \ C \Big \{ \left ( \parallel \nabla A \parallel_{\infty} \ + \
\parallel A \parallel_{\infty}^2 \ + \ \parallel V \parallel_{\infty}
\right ) \parallel \nabla \partial_t v \parallel_2 \nn \\ &&+ \parallel
\partial_t A \parallel_{\infty} \ \parallel \Delta v \parallel_2 \ +
\left ( \parallel A \parallel_{\infty} \ \parallel \nabla A
\parallel_{\infty} \ + \ \parallel \nabla V \parallel_{\infty}\right )
\parallel \partial_t v \parallel_2 \nn \\ &&+\left ( \parallel \nabla
\partial_t A \parallel_{\infty} \ + \ \parallel A \parallel_{\infty}\ 
\parallel \partial_t A \parallel_{\infty}\ + \ \parallel \partial_t V
\parallel_{\infty}\right ) \parallel \nabla v \parallel_{2} \nn \\
&&+\left ( \parallel \nabla \partial_t A^2 \parallel_{\infty} \ + \
\parallel \nabla \partial_t V \parallel_{\infty} \right ) \parallel v
\parallel_{2}\Big \} (t')  + \ \parallel \nabla \partial_t f(t') 
\parallel_{2} \ . \eea

\noi The RHS of (\ref{3.43e}) as well as those of
(\ref{3.23e})-(\ref{3.25e}) is in $L^1$ of the variable $t'$. By the
same argument as in Propositions 3.1 and 3.2, one obtains a solution
$v_{\eta} \in {\cal C}(I_+, H^2) \cap {\cal C}^1(I_+, H^1)$, with $I_+
= I \cap \{t:t\geq t_0\}$. \par

We next take the limit where $\eta$ tends to zero, and for that purpose
we need estimates of $v_{\eta}$ uniform in $\eta$ in the relevant
space. We shall need the evolution equation for $\nabla_A \partial_t v_{\eta}$, namely
\beq \label{3.44e} i \partial_t \nabla_A \partial_t v_{\eta} = - (1/2)
(1 - i \eta ) \nabla_A \Delta_A \partial_t v_{\eta} + V \nabla_A
\partial_t v_{\eta} + g_{\eta} \eeq

\noi where (see (\ref{3.41e}))
\bea \label{3.45e} g_{\eta}&=& \left ( \partial_t A + \nabla V \right )
\partial_t v_{\eta} + i (1 - i \eta ) \left ( (\partial_t A)\cdot
\nabla_A \otimes \nabla_A v_{\eta} + (\nabla \partial_t
A)\cdot \nabla_A v_{\eta} \right ) \nn \\ &&+ \left (
\partial_t V \right ) \nabla_A v_{\eta} + \left ( \nabla \partial_t V
\right ) v_{\eta} + \nabla_A \partial_t f\ . \eea

\noi Let $y_0 = \parallel v_{\eta}\parallel_2$ , $y_1 = \parallel
\partial_t v_{\eta}\parallel_2$ and $y = \parallel \nabla_A \partial_t
v_{\eta} \parallel_2$ . By a minor variation of Lemma 3.1, part (1),
$\nabla_A \partial_t v_{\eta}$ satisfies
\beq \label{3.46e} y(t_2)^2 - y(t_1)^2 \leq \int_{t_1}^{t_2} dt\ 2 \ {\rm
Im} <\nabla_A \partial_t v_{\eta}, g_{\eta}> (t) \eeq

\noi for all $t_1, t_2 \in I_+$, $t_1 \leq t_2$, so that 
\beq
\label{3.47e}
\partial_t y \leq \ \parallel g_{\eta}\parallel_2 \ .
\eeq

\noi We have already estimated $y_0$ by (\ref{3.15e}), $y_1$ by
(\ref{3.38e}) and $\parallel \Delta v_{\eta}\parallel_2$ by
(\ref{3.32e}) uniformly in $\eta$. We now estimate $y$ and at the same
time $\parallel \Delta \nabla v_{\eta}\parallel_2$ (which is not part of the
norm of the space of resolution). As in the proof of Proposition 3.2,
we omit the index $\eta$ in $v_{\eta}$ and we keep only the most
dangerous parts of $A$ and $V$ by dropping the $L^{\infty}$ components
allowed by the assumptions on $A$, $V$. From (\ref{3.45e})
(\ref{3.47e}) we obtain 
$$\partial_t y \leq \left ( \parallel
\partial_t A \parallel_3 \ + \ \parallel\nabla V \parallel_3 \right )
\parallel \partial_t v \parallel_6 \ + (1 + \eta^2)^{1/2} \Big (
\parallel \partial_t A \parallel_s \ \parallel\nabla_A^2 v \parallel_r$$
\beq \label{3.48e} 
+ \ \parallel\nabla \partial_t A \parallel_r \ \parallel \nabla_A v
\parallel_s \Big ) + \ \parallel \nabla
\partial_t V \parallel_2 \ \parallel v \parallel_{\infty}\ + \ \parallel
\partial_t V \parallel_6 \ \parallel \nabla_A v \parallel_3\ + \
\parallel \nabla_A \partial_t f \parallel_2 \ . \eeq

\noi We then estimate
\beq \label{3.49e} \parallel \nabla_A^2 v \parallel_r \ \leq \
\parallel\nabla^2 v \parallel_r \ + \ 2\parallel A \parallel_s \ \parallel \nabla v \parallel_{r_1} \ + \
\parallel A \parallel_s^2 \ \parallel v \parallel_{r_2}\ + \
\parallel\nabla A \parallel_r \ \parallel v \parallel_{\infty} \eeq
 \bea
 \label{3.50e}
&&\parallel \nabla_A v \parallel_s \ \leq \ \parallel\nabla v
\parallel_s \ + \ \parallel A \parallel_s \ \parallel v
\parallel_{\infty}\\  &&\nn \\ &&\parallel \nabla_A v
\parallel_3 \ \leq \ \parallel \nabla v \parallel_3 \ + \ \parallel A
\parallel_s \ \parallel v \parallel_{r_3} \label{3.51e} \eea

\noi where $1/r_1 = 1/2 - 2/s$, $1/r_2 = 1/2 - 3/s$, $1/r_3 = 1/3 -
1/s$, so that $r_1 \leq r_3 \leq 6 \leq s$. We next estimate by
(\ref{2.3e}) and Lemma 2.1
 \bea \label{3.52e} &&\parallel \partial_t v
\parallel_6  \ \leq  C \parallel \nabla_A \partial_t v \parallel_2 \ =
Cy\\  \nn \\ 
\label{3.53e}
&&\parallel\nabla^2v \parallel_r \ \leq C
\parallel\nabla \Delta v \parallel_2^{\delta} \ \parallel \Delta v
\parallel_2^{1 - \delta}\\  && \nn \\ &&\parallel\nabla v
\parallel_s \ \leq C \parallel\nabla \Delta v \parallel_2^{1/2 -
\delta} \ \parallel \Delta v \parallel_2^{1/2 + \delta} \label{3.54e} \eea

\noi where $0 \leq \delta = \delta (r) \equiv 3/2 - 3/r \leq 1/2$. Furthermore
\bea \label{3.55e} \parallel\nabla \Delta v \parallel_2 & \leq &
\parallel\nabla \Delta_A v \parallel_2 \ + 2\parallel A \parallel_s\ 
\parallel\nabla^2 v \parallel_r\ + 2 \parallel\nabla A\parallel_r\ 
\parallel\nabla v \parallel_s \nn \\ &&+ \ 2 \parallel \nabla A
\parallel_r \ \parallel A \parallel_s \ \parallel v \parallel_{\infty} \ +
\ \parallel A \parallel_s^2 \ \parallel \nabla v \parallel_{r_1}  \eea

\noi so that by (\ref{3.53e}) (\ref{3.54e}) and Lemma 2.2
\bea \label{3.56e} \parallel\nabla \Delta v \parallel_2 & \leq& C \Big
\{  \parallel\nabla \Delta_A v \parallel_2 \ + \left ( \parallel A
\parallel_s^{1/(1 - \delta)} \ + \ \parallel\nabla A
\parallel_r^{1/(1/2 + \delta )}\right )  \parallel\Delta v\parallel_2 
\nn \\ &&+ \ 2 \parallel \nabla A \parallel_r \ \parallel A \parallel_s\ 
\parallel v \parallel_{\infty} \ + \ \parallel A \parallel_s^2\ 
\parallel \nabla v \parallel_{r_1}\Big \} \ . \eea

\noi On the other hand, by a direct estimate of (\ref{3.3e}), we obtain
\bea \label{3.57e} &&\parallel\nabla \Delta_A v \parallel_2 \ \leq 2
\Big \{  y +\ C\parallel A  \parallel_s \ y_1^{1-\delta} y^{\delta} + \
\parallel \nabla V \parallel_2 \ \parallel v \parallel_{\infty }\ + \
\parallel V\parallel_6 \ \parallel \nabla v \parallel_3\nn \\ &&+ \ \parallel
\nabla f \parallel_2 \Big \} \eea

\noi where we have used (\ref{3.52e}). \par

We substitute (\ref{3.57e}) into (\ref{3.56e}), we substitute the
result into (\ref{3.53e}) (\ref{3.54e}),  we substitute the result into
(\ref{3.49e}) (\ref{3.50e}), and we substitute the result and
(\ref{3.51e}) (\ref{3.52e}) into (\ref{3.48e}). Using the fact that the
remaining norms of $v$ in (\ref{3.49e})-(\ref{3.51e}) and in
(\ref{3.56e}) (\ref{3.57e}) are controlled by $\parallel v;
H^2\parallel$ and using (\ref{3.32e}), we finally obtain an estimate
of the form 
\beq \label{3.58e} \partial_t y \leq N (y, y_1, y_0) \eeq

\noi where $N$ depends in addition on $A$, $V$ and $f$ through the norms
associated with the assumptions of the proposition, $N$ is homogeneous
of degree 1 in $y$, $y_0$, $y_1$ and $f$, and $N$ as a function of $y$
is the sum of a finite number of powers between 0 and 1. The estimate
(\ref{3.58e}) plays the same role in the proof of this proposition as
(\ref{3.36e}) in the proof of Proposition 3.2. Using the fact that
$y_0$ and $y_1$ have already been estimated uniformly in $\eta$ in the
proof of Proposition 3.2 and applying Lemma 2.3, we obtain an estimate
of the form 
\beq \label{3.59e} y(t) \leq P\left ( y(t_0), t\right ) \eeq

\noi where $P(z, t)$ is uniform in $\eta$, increasing in $z$ and
continuous and increasing in $t$ for $t \geq t_0$ with $P(z, t_0) = z$.
The estimate (\ref{3.59e}) plays the same role in the proof of this
proposition as (\ref{3.38e}) in the proof of Proposition 3.2. Now
$y(t_0)$ is estimated uniformly in $\eta$  for $v_0 \in H^3$. This
follows from the estimate  
\bea \label{3.60e} 
&&y(t_0) = \ \parallel\nabla_A (K_{\eta} v_0 + f) (t_0)
\parallel_2 \ \leq C \Big \{ \parallel\nabla \Delta v_0 \parallel_2 \ +
\ \parallel A \parallel_s \ \parallel \nabla^2 v_0 \parallel_r \nn \\
&&+\ \parallel\nabla A \parallel_r \ \parallel \nabla v_0\parallel_s \ + \
\parallel A \parallel_s^2 \ \parallel\nabla v_0\parallel_{r_1} \ + \
\parallel V \parallel_6 \ \parallel \nabla v_0 \parallel_3\nn \\
&&+ \left ( \parallel
A\parallel_s \ \parallel \nabla A \parallel_r\ + \ \parallel \nabla V
\parallel_2 \right ) \parallel v_0 \parallel_{\infty} \ + \ \parallel A
\parallel_s^3 \ \parallel v_0\parallel_{r_2} \nn \\
&& + \ \parallel A \parallel_s\ 
\parallel V \parallel_6 \ \parallel v_0\parallel_{r_3}\ + \ \parallel \nabla f \parallel_2 \ + \ \parallel A \parallel_s \ \parallel
f \parallel_r \Big \} \ . \eea

\noi The estimates (\ref{3.59e}) (\ref{3.60e}) provide an estimate of
$y$ uniform in $\eta$. Together with the estimates of $\partial_t v$ in
$L^2$ and of $v$ in $H^2$ that follow from (\ref{3.32e}) (\ref{3.38e})
and with (\ref{3.56e}) (\ref{3.57e}), they provide an a priori estimate
of $v$ in $L_{loc}^{\infty}(I_+, H^3) \cap H_{\infty loc}^1(I_+, H^1)$,
uniformly in $\eta$. \par

We can now take the limit $\eta \to 0$. The end of the proof is the
same as in Proposition 3.2 and will be omitted. \\

\noi \underbar{Part (2)} is proved in the same way as Part (2) of
Proposition 3.2.\par\nobreak \hfill $\sq$ \par

We now turn to the study of the Cauchy problem for (\ref{3.1e}) at the
level of regularity of $H^4$. We shall control the regularity of $v$ at
that level through the use of the second time derivative
$\partial_t^2v$. That derivative satisfies an evolution equation
obtained by applying $\partial_t^2$ to (\ref{3.1e}), namely
\beq
\label{3.61e}
i \partial_t \ \partial_t^2 v = K \partial_t^2 v + f_2
\eeq

\noi where
\beq
\label{3.62e}
f_2 = 2 \left ( \partial_t K \right ) \partial_t v + f_3
\eeq
\beq \label{3.63e} f_3 = \left ( \partial_t^2 K \right ) v +
\partial_t^2 f  = i \left ( \partial_t^2 A \right ) \cdot \nabla_A v
+ \left ( (\partial_t A)^2 + \partial_t^2 V  \right ) v + \partial_t^2
f \eeq

\noi and $\partial_t K$ can be read from (\ref{3.18e}). \par

In the present case however there arises a difficulty with the initial
condition. We shall perform the same regularization as in the proof of
Proposition 3.2, replacing (\ref{3.1e}) by (\ref{3.3e}) and in addition
regularizing $A$, $V$ and $f$ in space. We shall again use the integral
equation associated with (\ref{3.3e}), namely (\ref{3.9e}), and we shall need
the second time derivative of $\phi (v)$. By integration by parts,
that derivative is seen to be
$$ \left ( \partial_t^2 \phi (v)\right ) (t) = - U_{\eta}(t-t_0) \left
( K_{\eta}(K_{\eta} v_0 + f_{\eta}) + i \left ( \partial_t K_{\eta}
\right ) v_0 + i \partial_t f_{\eta} \right ) (t_0)$$
\beq
\label{3.64e}
- i \int_{t_0}^t dt'\ U_{\eta} (t-t') \partial_t^2 F_{\eta} (t') \ .
\eeq

\noi For smooth regularized $A$ and $V$, the domain of $K_{\eta}^2$ is
$H^4$, thereby suggesting to take $v_0 \in H^4$. However we shall make
only weak assumptions on the space regularity of $A$ and $V$, namely
assumptions of the $L^r$ type, but no assumptions on the space
derivatives. Under such assumptions, whereas the domain of $K$ for
$\eta = 0$ is easily seen to remain $H^2$, the domain of $K^2$ can be
very complicated and completely different from $H^4$. As a consequence,
for fixed $v_0 \in H^4$, $K_{\eta}^2v_0$ may very well blow up in the
limit $\eta \to 0$. This has two consequences. First the initial
condition $v_0$ should be chosen in a way adapted to $K$, namely such
that $K(Kv_0 + f) \in L^2$ for $t = t_0$ and $\eta = 0$. Such a $v_0$
will in general not be in $H^4$. Second we need to regularize the
initial $v_0$ to some $v_{0\eta}$ for $\eta > 0$. This is most simply
done by imposing the condition that $K_{\eta}(K_{\eta} v_{0\eta} +
f_{\eta})$ be independent of $\eta$. Actually for technical reasons it
is convenient to replace $K_{\eta}$ in that condition by $(\rho +
K_{\eta})$ for some sufficiently large positive $\rho$, such that $-
\rho$ belongs to the resolvent set of $K_{\eta}$ and that $(\rho +
K_{\eta})$ be invertible for all $\eta$. Thus we choose $v_0$ and we
regularize it in such a way that at $t = t_0$ 
\beq \label{3.65e} \left ( \rho + K_{\eta}\right )^2 v_{0\eta} + \left
( \rho + K_{\eta}\right ) f_{\eta} = (\rho + K)^2 v_0 + (\rho + K)f = z
\in L^2 \ , \eeq

\noi namely such that the LHS of (\ref{3.65e}) be independent of $\eta$
and be a fixed $z \in L^2$. In other words $v_0$ should be chosen as 
\beq
\label{3.66e}
v_0 = - (\rho + K)^{-1} f + (\rho + K)^{-2} z\eeq

\noi namely $v_0$ should be a given vector of ${\cal D}(K)$ modulo
${\cal D}(K^2)$, and $v_0$ should be regularized to
\beq \label{3.67e} v_{0\eta} = - \left ( \rho + K_{\eta}\right )^{-1}
f_{\eta} + \left ( \rho + K_{\eta}\right )^{-2} z\ . \eeq

As in Propositions 3.2 and 3.3, we shall need that all the initial
conditions for the relevant norms converge when $\eta \to 0$ namely
that $(\phi (v))t_0 = v_{0\eta}$ and $i(\partial_t \phi (v))(t_0) =
(K_{\eta} v_{0\eta} + f_{\eta})(t_0)$ converge in $H^2$ and that
\beq
\label{3.68e} - \left ( \partial_t^2 \phi (v)\right ) (t_0) = \left (
K_{\eta} (K_{\eta} v_{0\eta} + f_{\eta}) + i \left ( \partial_t
K_{\eta}\right ) v_{0\eta} + i \partial_t f_{\eta} \right ) (t_0) \eeq

\noi converge in $L^2$ as $\eta \to 0$. This will require some
information on the operator $K_{\eta}$ and in particular the
convergence of $K_{\eta}$ to $K$ in the strong resolvent sense
\cite{11r}.\\

\noi {\bf Lemma 3.3.} {\it Let $I$ be an interval, let $A \in {\cal
C}(I, L^6 + L^{\infty}) \cap {\cal C}^1(I, L^3 + L^{\infty})$, $V \in
{\cal C}^1(I, L^2 + L^{\infty})$. Let $0 \leq \eta \leq 1$. Let $K$ be
defined by (\ref{3.2e}) and let $K_{\eta}$ be defined by (\ref{3.4e})
with $A$ and $V$ regularized as in the proof of Proposition 3.2. Then,
for fixed $t = t_0 \in I$\par

(1) There exists $\rho > 0$ independent of $\eta$ such that $\rho + K$
and $\rho + K_{\eta}$ have bounded inverses from $L^2$ to $H^2$, with
$(\rho + K_{\eta})^{-1}$ uniformly bounded in $\eta$ as an operator
from $L^2$ to $H^2$.\par

(2) When $\eta \to 0$, $(\rho + K_{\eta})^{-1}$ converges strongly to
$(\rho + K)^{-1}$ and $(\partial_t K_{\eta})(\rho + K_{\eta})^{-1}$
converges strongly to $(\partial_t K)(\rho + K)^{-1}$ in $L^2$.}\\

\noi {\bf Proof.} \underbar{Part (1)} follows by standard arguments
from the fact that the $A$ and $V$ dependent parts in $K_{\eta}$ are a
Kato small perturbation of the Laplacian uniformly with respect to
$\eta \in [0, 1]$. In fact let $A = A_6 + A_{\infty}$ with $A_6 \in
L^6$ and $A_{\infty} \in L^{\infty}$, and similarly $V = V_2 +
V_{\infty}$. Then by Lemma 2.1
\bea \label{3.69e} &&\parallel iA \cdot \nabla v + (A^2/2)v \parallel_2
\ \leq \ \parallel A_6\parallel_6 \ \parallel\nabla v\parallel_3 \ + \
\parallel A_6 \parallel_6^2 \ \parallel v \parallel_6 \nn \\ &&+ \
\parallel A_{\infty}\parallel_{\infty} \ \parallel\nabla v \parallel_2 \
+ \ \parallel A_{\infty} \parallel_{\infty}^2 \ \parallel v \parallel_2
\nn \\ &&\leq C \left ( \parallel A_6\parallel_6 \ \parallel v
\parallel_2^{1/4} \ \parallel \Delta v \parallel_2^{3/4} \ + \ \left (
\parallel A_6 \parallel_6^2 \ + \ \parallel
A_{\infty}\parallel_{\infty} \right ) \parallel v \parallel_2^{1/2}\ 
\parallel \Delta v\parallel_2^{1/2}\right ) \nn \\ &&+ \
\parallel A_{\infty}\parallel_{\infty}^2 \ \parallel v \parallel_2 \nn \\
&& \leq \mu \parallel \Delta v \parallel_2 \ + C \left \{ \left (
\mu^{-3} + \mu^{-1}\right ) \parallel A_6\parallel_6^4 \ + \left (
\mu^{-1} + 1 \right ) \parallel A_{\infty}\parallel_{\infty}^2  \right
\} \parallel v \parallel_2 \eea

\noi for all $\mu > 0$, and similarly
\beq \label{3.70e} \parallel Vv \parallel_2 \ \leq \mu \parallel \Delta
v \parallel_2 \ + \left ( C \mu^{-3} \parallel V_2 \parallel_2^4 \ + \
\parallel V_{\infty} \parallel_{\infty} \right ) \parallel v \parallel_2
\ . \eeq

\noi The uniformity in $\eta$ follows from the fact that the
regularization does not increase the norms of $A$ and $V$ that appear
in (\ref{3.69e}) (\ref{3.70e}). \\

\noi \underbar{Part (2)}. We first remark that $K_{\eta} - K$ and
$\partial_t K_{\eta} - \partial_t K$ converge strongly to zero from
$H^2$ to $L^2$ when $\eta \to 0$. For $K_{\eta} - K$ and for the $A_6$
and $V_2$ components, this follows from the estimates (\ref{3.69e})
(\ref{3.70e}) applied to the differences and from the fact that the
regularization tends strongly to the identity in $L^r$ for $1 \leq r <
\infty$. For the $L^{\infty}$ components, the result follows from the
pointwise almost everywhere convergence of the regularized quantities
to the unregularized ones and from the dominated convergence theorem.
For $\partial_t K_{\eta} - \partial_tK$ the same argument applies with
the only difference that now $\partial_t A$ is decomposed as
$(\partial_t A)_3 + (\partial_t A)_{\infty}$ and that the contribution
of $(\partial_t A)_3$ is estimated by 
\bea \label{3.71e} \parallel \left ( \partial_t A\right )_3 \cdot \nabla_A v
\parallel_2  &\leq & \parallel \left ( \partial_t A\right )_3 \parallel_3 \left (
\parallel \nabla v \parallel_6 \ + \ \parallel A_6\parallel_6 \ \parallel
v \parallel_{\infty} \ + \ \parallel A_{\infty} \parallel_{\infty}\ 
\parallel v \parallel_6 \right ) \nn \\ &\leq& C \parallel \left ( \partial_t
A\right )_3 \parallel_3 \Big ( \parallel \Delta v \parallel_2 \ + \ \parallel
A_6\parallel_6 \ \parallel v \parallel_2^{1/4} \ \parallel \Delta v \parallel_2^{3/4} \nn \\
&&+ \ \parallel A_{\infty}
\parallel_{\infty} \ \parallel v \parallel_2^{1/2}  \ \parallel \Delta v
\parallel_2^{1/2} \Big) \ . \eea

\noi We now turn to the proof of strong resolvent convergence. Let $v \in L^2$. Then
\beq \label{3.72e} \left ( ( \rho + K)^{-1} - (\rho +
K_{\eta})^{-1}\right ) v = \left ( \rho + K_{\eta}\right )^{-1} \left (
K_{\eta} - K\right ) (\rho + K)^{-1} v \eeq

\noi which tends to zero strongly in $L^2$ by the previous convergence
of $K_{\eta} - K$ since $(\rho + K)^{-1} v$ is a fixed vector in $H^2$
and $(\rho + K_{\eta})^{-1}$ is uniformly bounded in $\eta$ as an
operator in $L^2$. Finally 
\bea
\label{3.73e}
&&\left \{ \left ( \partial_t K_{\eta}\right ) \left ( \rho + K_{\eta}\right )^{-1} - \left ( \partial_t K\right ) (\rho + K)^{-1} \right \} v = \left ( \partial_t K_{\eta} - \partial_t K \right ) (\rho + K)^{-1} v \nn \\
&&- \left ( \partial_t K_{\eta} \right ) \left ( \rho + K_{\eta} \right )^{-1} \left ( K_{\eta} - K\right ) (\rho + K)^{-1} v \ .
\eea

\noi The first term in the RHS converges to zero strongly in $L^2$ by
the previous strong convergence of $\partial_t K_{\eta}$ to $\partial_t
K$ as an operator from $H^2$ to $L^2$. The second term converges to
zero by the previous convergence of $K_{\eta}$ to $K$ and from the fact
that $(\partial_t K_{\eta})(\rho + K_{\eta})^{-1}$ is uniformly bounded
in $\eta$ as an operator in $L^2$. \par\nobreak \hfill $\sq$ \par

\noi {\bf Remark 3.4.} If it were not for the fact that the regularization
does not converge to the identity strongly in $L^{\infty}$, and in
particular if $A$, $V$ and $\partial_tA$, $\partial_t V$ did not have
$L^{\infty}$ components, we would obtain norm convergence instead of
strong convergence in Part (2) of the Lemma. \\

We can now state the result for the Cauchy problem for (\ref{3.1e}) at
the level of regularity of $H^4$.\\

\noi {\bf Proposition 3.4.} {\it Let $I$ be an interval. Let $A$, $V$
and $f$ satisfy $A \in {\cal C}(I, L^6 + L^{\infty}) \cap {\cal
C}^1(I, L^3 + L^{\infty})$, $\partial_t A \in L_{loc}^2(I, L^4 +
L^{\infty})$, $\partial_t^2A \in L_{loc}^1(I, L^3 + L^{\infty})$, $V \in
{\cal C}^1(I, L^2 + L^{\infty})$, $\partial_t^2V \in L_{loc}^1(I, L^2 +
L^{\infty})$, $f \in {\cal C}^1(I, L^2)$, $\partial_t^2f \in
L_{loc}^1(I, L^2)$. Let $t_0 \in I$ and let $v_0 \in H^2$ be such that
$v_0 + (\rho + K(t_0))^{-1} f(t_0) \in {\cal D}(K(t_0)^2)$ for some
$\rho > 0$ sufficiently large. Then\par

(1) There exists a unique solution $v \in {\cal C}^1(I, H^2) \cap {\cal
C}^2(I, L^2)$ of (\ref{3.1e}) in $I$ with $v(t_0) = v_0$. That solution
satisfies (\ref{3.6e}) and (\ref{3.19e}), namely
\beq \label{3.74e} \parallel \partial_t^j v(t_2) \parallel_2^2 \ - \
\parallel \partial_t^jv(t_1) \parallel_2^2 \ = \int_{t_1}^{t_2} dt \ 2\ 
{\rm Im} < \partial_t^jv, f_j>(t) \eeq

\noi for $j = 0,1$ and for all $t_1,t_2 \in I$, where $f_0 = f$ and
$f_1$ is defined by (\ref{3.18e}). That solution is actually unique in
$(L_{loc}^{\infty} \cap {\cal C}_w)(I,H^1)$. Furthermore $iK\partial_t
v = K(Kv + f) \in {\cal C}(I, L^2)$. \par

(2) Let in addition $A \in L_{loc}^2(I, L^{\infty})$, $V \in
L_{loc}^1(I, L^{\infty})$ and $\nabla A \in L_{loc}^1(I, L^{\infty})$.
Then the previous solution $v$ is actually unique in $(L^{\infty} \cap
{\cal C}_w)(I, L^2)$ and $v$ satisfies (\ref{3.74e}) for $j = 2$ and
for all $t_1,t_2 \in I$, where $f_2$ is defined by (\ref{3.62e})
(\ref{3.63e}).}\\

\noi {\bf Remark 3.5.} If $A$, $V$ and $f$ are sufficiently regular in
the space variable, the condition $K(Kv + f) \in {\cal C}(I, L^2)$ is
equivalent to the condition $v \in {\cal C}(I, H^4)$.\\

\noi {\bf Proof.} The proof proceeds by a parabolic regularization and
a limiting procedure as in the case of Proposition 3.2. We consider
first the case $t \geq t_0$. We replace (\ref{3.1e}) by (\ref{3.3e})
with $0 < \eta \leq 1$ and with $A$, $V$ and $f$ regularized as in
Proposition 3.2. We use again the integral equation for $v$, namely
(\ref{3.9e}), now however with the initial data $v_0$ regularized to
$v_{0\eta}$ according to (\ref{3.67e}) with $z$ defined by
(\ref{3.66e}), or equivalently by (\ref{3.65e}), and we solve that
equation locally in time by contraction in ${\cal C}([t_0, t_0 + T],
H^2) \cap {\cal C}^1 ([t_0, t_0 + T], H^1) \cap {\cal C}^2([t_0, t_0 +
T], L^2)$ for some $T > 0$. For that purpose we use again the estimates
(\ref{3.23e}) (\ref{3.24e}) (\ref{3.25e}) and (\ref{3.43e}),
supplemented by an additional estimate for $U_{\eta} (t-t')
\partial_t^2 F_{\eta}(t')$, namely
\bea \label{3.75e} &&\parallel U_{\eta}(t-t') \partial_t^2 F_{\eta}
(t')\parallel_2 \ \leq C ( \eta (t-t'))^{-1/2} \parallel
A(t')\parallel_{\infty} \ \parallel \partial_t^2 v(t')\parallel_2 \nn\\
&&+ C \Big \{ \parallel \partial_t A \parallel_{\infty} \ \parallel
\nabla \partial_t v \parallel_2\  + \ \parallel \partial_t^2 A
\parallel_{\infty} \ \parallel \nabla v \parallel_2 \nn \\ &&+ \left (
\parallel A \parallel_{\infty}^2 \ + \ \parallel V\parallel_{\infty}
\right ) \parallel \partial_t^2 v\parallel_2 \ + \left ( \parallel A
\parallel_{\infty} \ \parallel \partial_t A \parallel_{\infty}\ + \
\parallel \partial_t V \parallel_{\infty}\right ) \parallel \partial_t
v \parallel_2 \nn \\ &&+ \left ( \parallel \partial_t^2 A
\parallel_{\infty} \ \parallel A \parallel_{\infty} \ + \ \parallel
\partial_t A \parallel_{\infty}^2  \ + \ \parallel \partial_t^2 V
\parallel_{\infty}\right ) \parallel v \parallel_2 \Big \} + \
\parallel \partial_t^2 f \parallel_{2}  \ . \eea 

\noi The RHS of (\ref{3.75e}) as well as those of
(\ref{3.23e})-(\ref{3.25e}), (\ref{3.43e}) is in $L^1$ of the
variable $t'$. On the other hand, the choice (\ref{3.67e}) of the
regularized $(\phi (v))(t_0)$ yields (see (\ref{3.68e})) 
\beq
\label{3.76e}
(\phi (v))(t_0) = v_{0\eta} \in {\cal D}(K_{\eta}) = H^2 \ ,
\eeq
\beq \label{3.77e} \left ( i \partial_t \phi (v)\right ) (t_0) = \left
( K_{\eta} v_{0\eta} + f_{\eta} \right ) = - \rho v_{0\eta} + \left (
\rho + K_{\eta}\right )^{-1} z \in {\cal D}(K_{\eta}) = H^2 \ , \eeq
\bea \label{3.78e}
&&- \left  (  \partial_t^2 \phi (v)\right ) (t_0) = K_{\eta} \left (
K_{\eta} v_{0\eta} + f_{\eta} \right ) + i \left ( \partial_t K_{\eta}
\right ) v_{0\eta} + i \partial_t f_{\eta}\nn \\
&&= - \rho
K_{\eta} v_{0\eta} + K_{\eta} \left ( \rho + K_{\eta}\right )^{-1} z +
i \left ( \partial_t K_{\eta} \right ) v_{0 \eta} + i \partial_t f_{\eta}
\in L^2 \ . \eea

\noi Using those properties and the previous estimates, one obtains a
solution $v_{\eta} \in {\cal C}(I_+, H^2) \cap {\cal C}^1(I_+, H^1)
\cap {\cal C}^2(I_+, L^2)$ of (\ref{3.1e}) with $I_+ = I \cap \{t:t
\geq t_0\}$ by the same argument as in the proof of Proposition 3.1.
Furthermore it follows from (\ref{3.3e}), more precisely from
(\ref{3.26e}), that $\Delta \partial_t v_{\eta} \in {\cal C} (I_+,
L^2)$ so that $v_{\eta} \in {\cal C}^1(I_+, H^2)$. \par

We next take the limit where $\eta$ tends to zero and for that purpose
we need estimates of $v_{\eta}$ uniform in $\eta$ in the relevant
spaces. We remark that $v_{\eta}$ satisfies the equation (compare with
(\ref{3.61e})-(\ref{3.63e}))
\beq
\label{3.79e}
i \partial_t \ \partial_t^2 v_{\eta} = K_{\eta} \ \partial_t^2 v_{\eta} + f_{2 \eta}
\eeq

\noi where
\beq
\label{3.80e}
f_{2\eta} = 2\left ( \partial_t K_{\eta}\right ) \partial_t v_{\eta} + f_{3 \eta}
\eeq
\bea \label{3.81e} f_{3\eta} &=& \left ( \partial_t^2 K_{\eta}\right )
v_{\eta} + \partial_t^2 f_{\eta} \nn \\ &=&i (1 - i\eta) \left ( \partial_t^2
A\right ) \cdot \nabla_A  v_{\eta} + \left ( (1 - i \eta ) \left ( \partial_t
A\right )^2 + \partial_t^2 V \right ) v_{\eta} + \partial_t^2 f_{\eta}
\eea

\noi and $\partial_t K_{\eta}$ can be read from (\ref{3.22e}). \par 

Let
$y_j = \parallel \partial_t^j v_{\eta} \parallel_2$ , $j = 0,1,2$. By
Lemma 3.1, part (1), $\partial_t^2 v_{\eta}$ satisfies
\beq \label{3.82e} y_2(t_2)^2 - y_2(t_1)^2 \leq \int_{t_1}^{t_2} dt \ 2\ 
{\rm Im} < \partial_t^2 v_{\eta} , f_{2 \eta}>(t) \eeq

\noi for all $t_1, t_2 \in I_+$, $t_1 \leq t_2$, so that 
\beq
\label{3.83e}
\partial_t y_2 \leq \ \parallel f_{2\eta}\parallel_2 \ .
\eeq

\noi We have already estimated $y_0$ by (\ref{3.15e}) and $y_1$ in the
proof of Proposition 3.2. Now however the initial values of $y_0$ and
$y_1$ are
$$y_0 (t_0) = \parallel v_{0\eta}\parallel_2 \quad , \quad y_1(t_0) = \
\parallel \left ( K_{\eta} v_{0\eta} + f_{\eta}\right ) (t_0) \parallel_2$$

\noi by (\ref{3.76e}) (\ref{3.77e}). It follows from (\ref{3.15e}) and
from Lemma 3.3 that $y_0(t_0)$ and $y_1(t_0)$ are uniformly bounded in
$\eta$, so that the estimates of $y_0$ and $y_1$ are also uniform in
$\eta$. (See especially (\ref{3.36e}) (\ref{3.38e})). We next estimate
$y_2$, omitting again the index $\eta$ in $v_{\eta}$ for brevity and keeping only the
most dangerous parts of $A$ and $V$ by dropping again the $L^{\infty}$
components allowed by the assumptions on $A$ and $V$. From
(\ref{3.80e}) (\ref{3.81e}), by exactly the same estimates as in the
proof of Proposition 3.2 (see (\ref{3.29e}) (\ref{3.35e})), we obtain 
\beq \label{3.84e} \parallel f_{2\eta}\parallel \ \leq 2 M_0 \left (
\parallel \Delta \partial_tv \parallel_2\ , y_1 \right ) + \ \parallel
f_{3\eta}\parallel_2 \eeq

\noi where $M_0$ is defined in (\ref{3.35e}). On the other hand from
(\ref{3.26e}) and again by exactly the same estimates as in the proof
of Proposition 3.2 (see (\ref{3.30e}) (\ref{3.32e})) we obtain
\beq \label{3.85e} \parallel \Delta \partial_t v\parallel_2 \ \leq C
\left ( y_2 \ + \ \parallel f_{1\eta}\parallel_2\ + \ m^4 y_1 \right ) \eeq

\noi with the same constant $C$ as in (\ref{3.32e}). From
(\ref{3.83e})-(\ref{3.85e}) we then obtain
\beq \label{3.86e} \partial_t y_2 \leq \ \parallel
f_{2\eta}\parallel_2\ \leq 2 M \left ( y_2 \ + \ \parallel
f_{1\eta}\parallel_2\ ,  y_1 \right ) + \parallel f_{3\eta}\parallel_2
\eeq

\noi where $M$ is defined by (\ref{3.37e}). The $L^2$ norm of
$f_{1\eta}$ is already estimated by (\ref{3.35e}) and it remains only
to estimate $f_{3\eta}$. We obtain 
\bea \label{3.87e} \parallel f_{3\eta}\parallel_2 & \leq & \parallel
\partial_t^2 A\parallel_3 \left ( \parallel \nabla  v\parallel_6 \ + \
\parallel A\parallel_6 \ \parallel v\parallel_{\infty}\right )\nn \\ &+&
\left ( \parallel \partial_t A\parallel_4^2 \ + \ \parallel
\partial_t^2 V\parallel_2 \right ) \parallel v\parallel_{\infty} \ + \
\parallel \partial_t^2 f\parallel_2 \nn \\ &\leq & M_1 \left ( y_1\  + \
\parallel f\parallel_2\ , y_0\right ) + \ \parallel \partial_t^2
f\parallel_2\eea

\noi by Lemma 2.1 and (\ref{3.32e}), where
\bea \label{3.88e} &&M_1 (z, y_0) = C \Big \{ \parallel \partial_t^2
A\parallel_3 \left ( z + m^4y_0\right ) \nn \\ &&+ \left ( \parallel
\partial_t^2 A\parallel_3 \ \parallel A\parallel_6 \ + \ \parallel
\partial_t A\parallel_4^2 \ + \ \parallel \partial_t^2 V\parallel_2
\right ) \left ( y_0^{1/4} z^{3/4} + m^3 y_0 \right ) \Big \} \ . \eea

\noi Note in particular that the assumptions on the time derivatives of
$A$, $V$ and $f$ made in the proposition are taylored to ensure that
the estimate (\ref{3.87e}) is integrable in time uniformly with respect
to the regularization under the already available estimates on $y_0$,
$y_1$. From (\ref{3.86e}) (\ref{3.87e}), from the previous estimates of
$y_0$ and $y_1$ and from Lemma 2.1, it follows that $y_2$ satisfies an
estimate of the form 
\beq
\label{3.89e}
y_2 (t) \leq P_2\left ( y_2(t_0), t\right )
\eeq

\noi where $P_2(z,t)$ is uniform in $\eta$, increasing in $z$,
continuous and increasing in $t$ for $t \geq t_0$, and satisfies $P(z,
t_0) = z$. Actually that estimate is obtained from (\ref{3.38e}) by
replacing $y_0$, $y_1$, $f$, $\partial_tf$ by $y_1$, $y_2$,
$f_{1\eta}$, $f_{3\eta}$ and using the available estimates for $y_1$,
$f_{1\eta}$ and $f_{3\eta}$. Substituting (\ref{3.89e}) into
(\ref{3.85e}) then yields an estimate of $\parallel \Delta
\partial_tv\parallel_2$. Furthermore the initial value of $y_2$, namely 
\beq \label{3.90e} y_2(t_0) = \ \parallel \left ( K_{\eta} \left ( K_{\eta}  \ v_{0\eta}
+ f_{\eta}\right ) + i \left ( \partial_t K_{\eta}\right )
v_{0\eta} + i \partial_t f_{\eta}\right ) (t_0) \parallel_2 \eeq

\noi is uniformly bounded in $\eta$ by (\ref{3.67e}) (\ref{3.78e}) and
Lemma 3.3. Therefore the estimates (\ref{3.89e}) of $y_2$ and
(\ref{3.85e}) of $\parallel \Delta \partial_t v\parallel_2$ are also
uniform in $\eta$, so that $v$ is estimated in $H_{\infty loc}^1(I_+,
H^2) \cap H_{\infty loc}^2 (I_+, L^2)$ uniformly in $\eta$.\par

We can now take the limit $\eta \to 0$. For that purpose, we need the
strong convergence of the initial conditions as $\eta \to 0$, more
precisely the convergence of $(\phi (v))(t_0)$ and $(\partial_t \phi
(v))(t_0)$ in $H^2$ and the convergence of $(\partial_t^2 \phi
(v))(t_0)$ in $L^2$. Those convergences follow from
(\ref{3.76e})-(\ref{3.78e}), from (\ref{3.67e}) and from Lemma 3.3.
With that information available, the end of the proof is the same as
that of Proposition 3.2 and will be omitted.\\

\noi \underbar{Part (2)} is proved in the same way as Part (2) of
Proposition 3.2. \par\nobreak \hfill $\sq$

\mysection{Scattering theory at the level of L$^{\bf 2}$ and H$^{\bf 2}$ for~w}
\hspace*{\parindent}
In this section we begin the study of scattering theory for
(\ref{1.1e}) with a potential $A$ satisfying conditions of the type
(\ref{2.26e}), or equivalently for (\ref{1.16e}) with a potential $B$
satisfying conditions of the type (\ref{2.27e}). Here we study that
theory for (\ref{1.1e}) in the spaces $L^2$ and $FH^2$, or equivalently
for (\ref{1.16e}) in the spaces $L^2$ and $H^2$. This will be done by
studying (\ref{1.12e}) in $L^2$ and $H^2$ and will rely on
Propositions 3.1 and 3.2. \par

The main result of this section has been stated as Proposition 1.1 in the introduction and is repeated here as the following proposition.\\ 

\noi{\bf Proposition 4.1.} {\it Let $A$ satisfy 
$$\parallel P^j \partial_x^{\alpha} A
\parallel_r \ \vee \ \parallel P^j (x\cdot A) \parallel_r \ \leq C \ t^{-1 + 2/r} \eqno(1.20)\equiv (4.1)$$

\noi where $P = t \partial_t + x \cdot \nabla$, for $0 \leq j + |\alpha | \leq 1$, $2 \leq r \leq \infty$ and for all
$t \in [1, \infty )$. \par

(1) Let $X = L^2$ or $FH^2$. Then for any $u_+ \in X$, there exists a
unique solution $u$ of (\ref{1.1e}) such that $\widetilde{u} \in {\cal C}([1, \infty ), X)$ and such
that
$$\parallel \widetilde{u}(t) - u_+; X \parallel \ \to 0 \quad \hbox{\it when $t \to \infty$} \ .
\eqno(1.21)\equiv (4.2)$$

\noi Conversely for any solution $u$ of (\ref{1.1e}) such that
$\widetilde{u} \in {\cal C}([1, \infty ), X)$, there exists $u_+ \in X$
such that (4.2) holds.

(2) Let $X = L^2$ or $H^2$. Then for any $w_+ \in X$, there exists a
unique solution $w_* \in {\cal C}([1, \infty ), X)$ of (\ref{1.16e})
such that 
$$\parallel w_*(t) - w_+;X\parallel \ \to 0\quad \hbox{\it when $t \to \infty$} \ .
\eqno(1.22)\equiv (4.3)$$

\noi Conversely for any solution $w_* \in {\cal C}([1, \infty ), X)$ of
(\ref{1.16e}), there exists $w_+ \in X$ such that (4.3)
holds.} \\

\noi {\bf Proof.} By (\ref{1.19e}), Parts (1) and (2) are equivalent,
with $w_+ = Fu_+$. We concentrate on Part (2). By (\ref{1.11e})
(\ref{2.20e}) (\ref{2.21e}), the assumption (4.1) on $A$ can be
rewritten in terms of $B$ as
$$\parallel\partial_t^j \  \partial_x^{\alpha} B
\parallel_r \ \leq C\ t^{-j-|\alpha | + 1/r}\ , \ \parallel  \partial_t^j \hbox{\it \v B} \parallel_r \ \leq C \ t^{-j + 1/r}\ .$$

For $X = H^2$, Part (2) is then obtained as an immediate consequence of
Proposition 3.2, part (1) applied with $(v, A, V, f)$ replaced by $(w,
B, - \hbox{\it \v B}, 0)$ with $w(t) = \overline{w_*(1/t)}$. In fact, from the previous assumption on $B$, 
$$\parallel B(t)\parallel_6 \ \leq C\ t^{1/6} \ , \ \parallel
\partial_t B(t) \parallel_3 \ \leq C\ t^{-2/3}\ , \ \parallel
\partial_t \hbox{\it \v B}(t)\parallel_2\ \leq C\ t^{-1/2}$$

\noi so that the assumptions of Proposition 3.2, part
(1) are satisfied. \par

For $X = L^2$, the situation is slightly more delicate. For the same
choice of $(v, A, V, f)$, the assumptions of Proposition 3.1, part (1)
are satisfied in $I = [0, 1]$, but the assumptions of Proposition 3.1, part (2) are
satisfied only in $(0, 1]$ because
$$\parallel \nabla B\parallel_{\infty} \ \leq C\ t^{-1}$$

\noi is not integrable at $t= 0$. We shall therefore combine
Proposition 3.1 in $(0, 1]$ with an approximation argument using
Proposition 3.2, part (1). Let first $t_0 \in [0, 1]$ and $w_0 \in L^2$.
We approximate $w_0$ by a sequence $\{w_{0j}\}$ in $H^2$ converging
strongly to $w_0$ in $L^2$. By Proposition 3.2, part (1), each $w_{0j}$
generates a solution $w_j \in {\cal C}([0,1], H^2)$ of (\ref{1.12e})
with $w_j (t_0) = w_{0j}$. Furthermore $L^2$ norm conservation holds for
those solutions so that
$$\parallel w_j(t) - w_l(t)\parallel_2 \ = \ \parallel w_{0j} -
w_{0l}\parallel_2 \qquad \hbox{for all $t \in [0, 1]$} \ .$$

\noi Therefore $w_j$ converges in norm in $L^{\infty}([0,1], L^2)$ to a
solution $w \in {\cal C}([0, 1], L^2)$, with constant $L^2$ norm. Using
that result with $t_0 = 0$ and $w_0 = \overline{w_+}$ yields the existence
part of the first statement of Part (2). We next prove uniqueness. Let
$w_1, w_2 \in {\cal C}([0,1],L^2)$ be two solutions of (\ref{1.12e})
in $[0, 1]$ with $w_1(0) = w_2(0) = \overline{w_+}$. By Proposition 3.1
applied in $(0,1]$, $w_1 - w_2$ satisfies $L^2$ norm conservation,
namely 
$$\parallel w_1(t) - w_2(t) \parallel_2 \ = C \quad \hbox{for $t \in (0,1]$}\ .$$ 

\noi Since $w_1$ and $w_2$ have the same strong $L^2$ limit as $t \to
0$, the last constant is zero and therefore $w_1 = w_2$. \par

We finally prove the second statement of Part (2).\par

Let $w \in {\cal C}((0, 1], L^2)$ be a solution of (\ref{1.12e}). By
Proposition 3.1 applied in $(0, 1]$, $w$ is uniquely determined in
${\cal C}((0, 1], L^2)$ (actually in $({\cal C}_w \cap
L_{loc}^{\infty})((0, 1], L^2))$ by its value $w_0 = w(t_0)$ for some
$t_0 > 0$. By the previous $H^2$ approximation method, we can construct
a solution $w' \in {\cal C}([0, 1];L^2)$ of (\ref{1.12e}) with $w'(t_0)
= w_0$. By uniqueness, $w = w'$ in $(0,1]$ and therefore $w$ has a
strong limit in $L^2$ as $t \to 0$. This proves the second statement of
Part (2).\par\nobreak \hfill $\sq$ \par

In the language of scattering theory, Proposition 4.1 essentially
expresses the existence of the wave operators and their asymptotic
completeness for the equation (\ref{1.1e}) in $L^2$ and in $FH^2$, as
compared with the free Schr\"odinger evolution, and for the equation
(\ref{1.16e}) in $L^2$ and in $H^2$, as compared with the constant
evolution. \par

In Proposition 4.1, we have obtained the existence and asymptotic
completeness of the wave operators by using Propositions 3.1 and 3.2. If
one is only interested in the existence of the wave operators in $L^2$,
one can avoid using Proposition 3.2 and use only Proposition 3.1. As a
consequence no assumption is needed on the time derivative of $B$ or equivalently on $PA$. \\

\noi {\bf Proposition 4.2.} {\it Let $A$ satisfy 
$$\parallel \partial_x^{\alpha} A \parallel_r \ \vee \ \parallel x \cdot A\parallel_r \
\leq C\ t^{-1+2/r} \eqno(4.4)$$

\noi for $0 \leq |\alpha | \leq 1$, $2 \leq r \leq \infty$ and for all $t \in [1, \infty )$. Then\par

(1) For any $u_+ \in L^2$, there exists a unique solution $u \in {\cal
C}([1, \infty ), L^2)$ of (\ref{1.1e}) such that
$$
\parallel \widetilde{u}(t) - u_+ \parallel_2 \ \to 0 \qquad \hbox{\it when $t \to \infty $} \ .
\eqno(4.5)$$

(2) For any $w_+ \in L^2$, there exists a unique solution $w_* \in
{\cal C}([ 1 , \infty ), L^2)$ of (\ref{1.16e}) such that}
$$\parallel w_*(t) - w_+ \parallel_2 \ \to 0 \qquad \hbox{\it when $t \to \infty $} \ .
\eqno(4.6)$$

\noi {\bf Proof.} Again Parts (1) and (2) are equivalent by
(\ref{1.19e}) with $w_+ = Fu_+$ and can be rephrased in an obvious way
in terms of $w$ and of (\ref{1.12e}). The assumption (4.4) on $A$ can be rewritten in terms of $B$ as
$$\parallel \partial_x^{\alpha} B \parallel_r \ \leq \ C\ t^{-|\alpha|+1/r}\ , \quad \parallel \hbox{\it \v B} \parallel_r \
\leq C\ t^{1/r} \ .$$

\noi Under that assumption Proposition 3.1 holds in $(0, 1]$. Let $w_+ \in H^2$. Let $t_0 >
0$. By Proposition 3.1, there exists a (unique) solution $w_{t_0} \in
{\cal C}((0, 1], L^2)$ of (\ref{1.12e}) such that $w_{t_0}(t_0) =
\overline{w_+}$. Furthermore for $t_0 > t_1 > 0$, $w_{t_0} - w_{t_1}$
satisfies $L^2$ norm conservation in $(0, 1]$, while $w_{t_0} -
\overline{w_+}$ satisfies the equation
$$i \partial_t \left ( w_{t_0} - \overline{w_+}\right ) =
\left ( - (1/2) \Delta_B - \hbox{\it \v B} \right ) \left ( w_{t_0} -
\overline{w_+}\right ) - R(\overline{w_+})\eqno(4.7)$$

\noi where $R(\cdot )$ is defined by (\ref{1.32e}), so that 
$$R(\overline{w_+}) = \left ( (1/2) \Delta_B + \hbox{\it \v B}\right ) \overline{w_+} \ .
\eqno(4.8)$$

\noi Therefore, for all $t \in (0, 1]$
$$\parallel w_{t_0}(t) - w_{t_1}(t) \parallel_2 \ = \ \parallel
w_{t_0}(t_1) - \overline{w_+} \parallel_2 \ \leq \int_{t_1}^{t_0} dt \parallel
R(\overline{w_+})\parallel_2 \ .$$

\noi Now
$$\parallel R(\overline{w_+})\parallel\ \leq \ \parallel
\Delta w_+ \parallel_2 \ + \ \parallel B\parallel_3 \ \parallel
\nabla w_+ \parallel_6\ + \left ( \parallel B \parallel_4^2 \ + \
\parallel \hbox{\it \v B} \parallel_2 \right ) \parallel w_+
\parallel_{\infty} \ \leq C \eqno(4.9)$$

\noi so that 
$$\parallel w_{t_0}(t) - w_{t_1}(t) \parallel_2\ \leq
C |t_1 - t_0| \qquad \hbox{for all $t \in (0, 1]$} \ . \eqno(4.10)$$

\noi It follows from (4.10) that $w_{t_0}$ converges in norm in
$L^{\infty}((0, 1], L^2)$ to a limit $w \in {\cal C}((0,1],L^2)$ which
is also a solution of (\ref{1.12e}) when $t_0 \to 0$. Furthermore by taking the limit
$t_1 \to 0$ in (4.10), we obtain

$$\parallel w_{t_0}(t) - w(t) \parallel_2\ \leq C\ t_0
\qquad \hbox{for all $t \in (0, 1]$} \eqno(4.11)$$

\noi and in particular
$$\parallel w(t_0) - \overline{w_+} \parallel_2\ \leq C\ t_0
\eqno(4.12)$$

\noi so that $w$ can be extended to ${\cal C}([0, 1], L^2)$ with $w(0)
= \overline{w_+}$. \par

This proves the existence part in the proposition in the special case
where $w_+ \in H^2$. The proof for general $w_+ \in L^2$ follows
therefrom by the same approximation argument as in Proposition 4.1.
\par

The proof of the uniqueness part is the same as in Proposition 4.1
since it does not use Proposition 3.2. \par\nobreak \hfill $\sq$\par

In the framework of Proposition 4.1, if we assume additional regularity
properties of $w_+$ and $u_+$, we obtain stronger convergence
properties than (4.2) (4.3), in the form of time decay as powers of
$t$. We have stated typical results of this type, expressed in terms of
$u$ and of $w_*$, as Proposition 1.2 in the introduction. That
proposition follows as a special case (namely with $m=1$) of the
following proposition, where the corresponding results are stated in
terms of $w$.\\

\noi {\bf Proposition 4.3.} {\it Let $A$ satisfy the assumptions of
Proposition 4.1. Let $w_+ \in H^2$ and let $w \in {\cal C}([0, 1],
H^2)$ be the solution of (\ref{1.12e}) with $w(0) = \overline{w_+}$ obtained
by Proposition 3.2. Then the following estimates
hold~:\par

(1)
$$\parallel w(t) - \overline{w_+} \parallel_2 \ \leq C \ t \ .
\eqno(4.13)$$

(2) Let in addition $w_+ \in H^{2+m}$ for $m \geq 0$. Then
$$\parallel w(t) - U(t) \overline{w_+} \parallel_2 \ \leq
\left \{ \begin{array}{ll} C\ t^{(4 + m)/3} &\hbox{\it for $m < 1/2$}
\\ &\\ C\ t^{(3 - \varepsilon)/2} &\hbox{\it for $m = 1/2$} \\ &\\ C\
t^{3/2} &\hbox{\it for $m > 1/2$} \ . \end{array} \right . \eqno(4.14)$$

(3) Let $w_+ \in H^3$. Then}
$$\parallel \partial_t (w(t) - U(t) \overline{w_+} )
\parallel_2 \ \vee\ \parallel \Delta (w(t) - U(t)
\overline{w_+})\parallel_2\ \leq C\ t^{1/2} \ . \eqno(4.15)$$

\noi {\bf Proof of Proposition 1.2.} Part (2) of that proposition
follows from the special case $m = 1$ of Proposition 4.3 and from
(\ref{1.15e}). Part (1) follows from Part (2) and from (\ref{1.19e}).
\\

\noi {\bf Proof of Proposition 4.3.} \underbar{Part (1)}. In the same way as in the proof
of Proposition 4.2, we estimate
$$\left | \partial_t \parallel w(t) -
\overline{w_+}\parallel_2 \right | \ \leq \ \parallel
R(\overline{w_+})\parallel_2\ \leq C \eqno(4.16)$$

\noi from which (4.13) follows by integration.\\

\noi \underbar{Part (2)}. We estimate similarly
$$\left | \partial_t \parallel w(t) - U(t)
\overline{w_+}\parallel_2 \right | \ \leq \ \parallel
R\left ( U(t)\overline{w_+}\right )\parallel_2
\eqno(4.17)$$

\noi where $R(\cdot )$ is defined by (\ref{1.32e}). We compute
$$R\left ( U(t)\overline{w_+}\right ) = - i B \cdot \nabla
U(t)\overline{w_+} - \left ( B^2/2 - \hbox{\it \v B}\right ) U(t) \overline{w_+}
\eqno(4.18)$$

\noi and we estimate
$$ \parallel B \cdot  \nabla U(t) \overline{w_+} \parallel_2 \ \leq
\left \{ \begin{array}{ll} C \parallel B \parallel_{3/(1+m)}\ \parallel
\omega^{m+2} w_+ \parallel_2 \ \leq C\ t^{(1 + m)/3} &\hbox{for $m
< 1/2$} \\ &\\ C \parallel B \parallel_{2/(1- \varepsilon )} \
\parallel w_+;H^{3/2}\parallel\ \leq C\ t^{(1 - \varepsilon)/2}
&\hbox{for $m = 1/2$} \\ &\\ C \parallel B \parallel_2 \ \parallel
w_+;H^{2+m}\parallel\ \leq C\ t^{1/2} &\hbox{for $m > 1/2$} \ ,
\end{array} \right . \eqno(4.19)$$
$$\parallel ( B^2/2 - \hbox{\it \v B}) U(t) \overline{w_+}
\parallel_2 \ \leq C \left ( \parallel B \parallel_4^2 \ + \
\parallel\hbox{\it \v B}\parallel_2 \right ) \parallel w_+;H^2\parallel\
\leq C\ t^{1/2} \eqno(4.20)$$

\noi by (4.1) and Lemma 2.1. The result follows from
(4.17) (4.19) (4.20) by integration on time.\\

\noi \underbar{Part (3)}. We use the estimates in the proof of
Proposition 3.2 applied with $(v, A, V, f)$ replaced by $(w - U(t)
\overline{w_+}, B, - \hbox{\it \v B}, - R(U(t)\overline{w_+}))$. In addition to
the estimate
$$\parallel R(U(t) \overline{w_+}) \parallel_2 \ \leq C \
t^{1/2} \eqno(4.21)$$

\noi which follows from (4.19) (4.20), we need the
estimate
$$\parallel \partial_t R(U(t) \overline{w_+}) \parallel_2
\ \leq  \ \parallel \partial_t B\parallel_2 \ \parallel\nabla U(t)
\overline{w_+} \parallel_{\infty} \ + \ \parallel B \parallel_{\infty} \
\parallel\nabla \Delta w_+ \parallel_2 $$ $$+ \ \parallel B
\partial_t B - \partial_t \hbox{\it \v B} \parallel_2 \ \parallel U(t)
\overline{w_+} \parallel_{\infty} \ + \left ( \parallel B
\parallel_{\infty}^2 \ + \ \parallel \hbox{\it \v B}\parallel_{\infty}
\right ) \parallel\Delta w_+ \parallel_2$$ $$\leq C\parallel
w_+;H^3\parallel\ t^{-1/2} \leq C\ t^{-1/2}\eqno(4.22)$$

\noi by (4.1) and Lemma 2.1.\par

From (4.21) it follows immediately by integration that
$$\parallel w - U(t) \overline{w_+} \parallel_2\ \leq C \ t^{3/2} \ .
\eqno(4.23)$$

\noi From the estimate (\ref{3.38e}), it follows by substituting
(4.21)-(4.23) that

$$\parallel \partial_t \left ( w - U(t)
\overline{w_+}\right )\parallel_2 \ \leq C \ t^{1/2} \eqno(4.24)$$

\noi which yields the estimate of the first term in (4.15).
Substituting (4.24) (4.21) (4.23) into
(\ref{3.32e}) yields the estimate of the second term in (4.15).\par\nobreak \hfill $\sq$ 

\mysection{H$^{\bf k}$  regularity of the wave operators for u}
\hspace*{\parindent} In this section we study the theory of scattering
for (\ref{1.1e}) in spaces $H^k$ for $k > 0$ for sufficiently smooth
$A$ satisfying conditions of the type (\ref{2.26e}). We have already
obtained $L^2$ wave operators satisfying asymptotic completeness in
$L^2$ for (\ref{1.1e}) in Proposition 4.1, part (1), and the problem is
that of additional regularity for those wave operators. We restrict our
attention to the case where $k$ is an even integer. For brevity in all
this section we shall take for granted the existence of solutions of
(\ref{1.1e}) at the required level of regularity and we shall
concentrate on the derivation of higher norm estimates. The existence
results follow from Propositions 3.1-3.4 for $k \leq 4$ and can be proved in
the same way for higher $k$. We first derive bounds for higher norms of
generic solutions of (\ref{1.1e}). Those bounds are unfortunately not
uniform in $t$. \\

\noi {\bf Proposition 5.1.} {\it Let $j \geq 1$ be an integer and let
$A$ satisfy the estimates
\beq \label{5.1e} \parallel \partial_t^l A \parallel_{\infty} \ \leq C\
t^{-1} \qquad \hbox{\it for $0 \leq l \leq j$} \ . \eeq

\noi Let $u \in {\cal C}^j([1, \infty ), L^2) \cap {\cal C}^{j-1}([1,
\infty ) ,H^2)$ be a solution of (\ref{1.1e}). Then $u$ satisfies the
estimates
\beq \label{5.2e} \parallel \partial_t^j u \parallel_2 \ \vee \
\parallel \partial_t^{j-1} \Delta u \parallel_2 \ \leq C\left ( 1 +
\ell n \ t\right )^{2j} \ .\eeq

\noi Let in addition $A$ satisfy the estimates \beq \label{5.3e}
\parallel \partial_t^l \ \partial_x^{\alpha} \ A\parallel_{\infty} \ \leq C
\ t^{-1} \quad \hbox{\it for $0 \leq |\alpha |/2 + l \leq j - 1$} \eeq

\noi and let $u \in \displaystyle{\mathrel{\mathop {\cap}_{0 \leq l
\leq j}}} {\cal C}^{j-l}([1, \infty ), H^{2l})$. Then $u$ satisfies the
estimates}
\beq \label{5.4e} \parallel \partial_t^{j-l} \Delta^l u \parallel_2\
\leq C (1 + \ell n \ t)^{2j} \quad \hbox{for $0 \leq l \leq j$} \ . \eeq

\noi {\bf Proof.} The proof proceeds by induction on $j$ and possibly
$l$, starting from (\ref{1.1e}). Since for the proof of Proposition
5.2 below, we shall need a similar induction for a slightly more
general equation, we shall replace (\ref{1.1e}) by the more
general inhomogeneous equation
\beq
\label{5.5e}
i \partial_t u = - (1/2) \Delta_A u + f \ .
\eeq

\noi For the needs of the present proof we shall take $f = 0$ at the
end. Taking the $j$-th time derivative of (\ref{5.5e}), we obtain 
\bea \label{5.6e} &&i \partial_t^{j+1} u = - (1/2) \Delta_A \partial_t^j
u + \sum_{0 \leq l < j} i\ C_j^l \left \{ \partial_t^{j-l-1}
((\partial_t A)\cdot \nabla_A) \right \} \partial_t^l u + \partial_t^j f\nn
\\ &&= - (1/2)\Delta_A \partial_t^j u + \sum_{0 \leq l < j} C_j^l \Big
\{ i \left ( \partial_t^{j-l} A\right ) \cdot \nabla_A \partial_t^l u \nn \\ &&+ \sum_{0 < m
< j-l} C_{j-l-1}^m \left ( \partial_t^{j-l-m} A \right ) \left (
\partial_t^m A\right ) \partial_t^l u\Big \}  + \partial_t^j f \eea

\noi where we have used the relation (see (\ref{3.18e}))
\beq
\label{5.7e}
- (1/2) \left ( \partial_t \Delta_A \right ) = i \left ( \partial_t A \right ) \cdot\nabla_A \ .
\eeq

\noi We now define $$y_l = \ \parallel \partial_t^{l} u\parallel_2\
\hbox{for $0 \leq l \leq j$}\ ,\ z_l = \ \parallel \Delta_A
\partial_t^{l-1} u \parallel_2\ \hbox{for $1 \leq l \leq j$, $z_0 = 0$} \ .$$

\noi In the same way as in Section 3 and by using (\ref{5.1e}), we
estimate 
\beq \label{5.8e} \left | \partial_t y_j \right | \leq C \
t^{-1} \sum_{0\leq l < j} \left ( \parallel \nabla_A \partial_t^l u
\parallel_2\ + t^{-1} y_l \right ) + \ \parallel \partial_t^j f
\parallel_2\ . \eeq

\noi Now
\beq \label{5.9e} \parallel \nabla_A \partial_t^l u \parallel_2 \ \leq
\left (y_l z_{l + 1} \right )^{1/2} \leq (1/2) \left ( y_l + z_{l+1}
\right ) \ . \eeq

\noi We substitute (\ref{5.9e}) into (\ref{5.8e}), using the middle
bound of (\ref{5.9e}) for $l = j - 1$ and the last bound of
(\ref{5.9e}) for $l < j - 1$. We obtain 
\beq \label{5.10e} \left | \partial_t y_j \right | \leq C \ t^{-1}
\left \{ \left ( y_{j-1} z_j \right )^{1/2} + \sum_{0\leq l < j} \left
( y_l + z_l \right ) \right \} + \ \parallel \partial_t^j f \parallel_2 \ . \eeq

\noi On the other hand by a similar estimate, we obtain

\beq \label{5.11e} z_j \leq 2y_j + C\ t^{-1} \sum_{0 \leq l < j-1}
\left ( \parallel \nabla_A \partial_t^l u \parallel_2 \ + t^{-1} y_l
\right ) + \ 2\parallel\partial_t^{j-1} f \parallel_2 \eeq

\noi so that by (\ref{5.9e}) again
\beq
\label{5.12e}
z_j \leq 2y_j + C\ t^{-1} \sum_{0 \leq l < j}
\left ( y_l + z_l \right ) + \ 2\parallel\partial_t^{j-1} f \parallel_2 \ .
\eeq

\noi We now prove the estimate (\ref{5.2e}) by induction on $j$ through
the use of (\ref{5.10e}) (\ref{5.12e}) with $f = 0$. The starting point
is $y_0 = C$ by $L^2$ norm conservation. We next assume that 
$$y_l \vee z_l \leq C (1 + \ell n \ t)^{2l}\qquad \hbox{for $0 \leq l < j$}\ .$$

\noi Substituting that assumption into (\ref{5.10e}) (\ref{5.12e}) yields 
\beq \label{5.13e} \left | \partial_t y_j \right | \leq C \ t^{-1}
\left \{ ( 1 + \ell n\ t )^{j-1} \ z_j^{1/2} + (1 + \ell n \ t )^{2(j-1)} \right \}\eeq
\beq
\label{5.14e}
z_j \leq 2y_j + C\ t^{-1} (1 + \ell n t)^{2(j-1)}\eeq

\noi so that 
\beq \label{5.15e} \left | \partial_t y_j \right | \leq C \ t^{-1}
\left \{ ( 1 + \ell n\ t )^{j-1} \ y_j^{1/2} + (1 + \ell n \ t )^{2(j-1)} \right \}\ .\eeq

\noi Integrating over time by Lemma 2.3 yields the first estimate of
(\ref{5.2e}), from which the second estimate follows by (\ref{5.14e})
and (\ref{5.1e}).\par

We now turn to the proof of (\ref{5.4e}). For that purpose we define
$$z_j^l = \ \parallel \partial_t^{j-l} \Delta^l u \parallel_2 \qquad
\hbox{for $0 \leq l \leq j$}$$

\noi so that $z_j^0 = y_j$ is already estimated by (\ref{5.2e}). Now for $1 \leq l \leq j$, 
\bea \label{5.16e} &&z_j^l = \ \parallel \partial_t^{j-l} \Delta^{l-1}
\left ( (- 2 i \partial_t + 2iA\cdot \nabla + A^2 ) u + 2f \right
)\parallel_2 \nn \\ &&\leq 2 z_j^{l-1} + C \sum_{0\leq m \leq j-l} \
\sum_{|\alpha + \beta | = 2(l-1)} \Big \{ \parallel \partial_t^m
\partial_x^{\alpha} A \parallel_{\infty} \ \parallel \partial_t^{j-l-m}
\partial_x^{\beta} \nabla u \parallel_2 \nn \\ &&+ \ \parallel
\partial_t^m \partial_x^{\alpha} (A^2) \parallel_{\infty} \ \parallel
\partial_t^{j-l-m} \partial_x^{\beta} u \parallel_2 \Big \} + \ 2
\parallel \partial_t^{j-l} \Delta^{l-1}f \parallel_2 \nn \\ &&\leq 2
z_j^{l-1} + C \ t^{-1} \left \{ \left ( z_j^l z_{j-1}^{l-1}\right
)^{1/2} + \sum_{m \leq k \leq j-1} z_k^m \right \} + \ 2 \parallel
\partial_t^{j-l} \Delta^{l-1} f \parallel_2 \eea

\noi by (\ref{5.3e}) and
\beq \label{5.17e} \parallel \partial_t^{j-l} \Delta^{l-1}\nabla u
\parallel\ \leq \left ( z_j^l z_{j-1}^{l-1}\right )^{1/2} \leq {1 \over
2} \left ( z_j^l + z_{j-1}^{l-1}\right ) \ . \eeq

We now prove the result by induction on $j$ and for each $j$ by
induction on $l$ starting from $z_0^0 = y_0 = C$. Thus we assume that
$z_k^m \leq C (1 + \ell n\ t)^{2k}$ for $0 \leq m \leq k \leq j-1$ and
for $k = j$, $0 \leq m \leq l-1$. The estimate (\ref{5.4e}) follows
immediately from (\ref{5.16e}) with $f = 0$ and from that
assumption.\par\nobreak \hfill $\sq$ \par

Although we are unable to prove the $H^k$ boundedness of general $H^k$
solutions of (\ref{1.1e}), we can nevertheless construct densely defined $H^k$ wave
operators for that equation, namely construct solutions that are $H^k$
bounded and that are asymptotic in the $H^k$ sense to prescribed model
asymptotics of the type $u_a = U(t) u_+$ with asymptotic state $u_+$ in a dense subspace of
$H^k$. The method of construction is that sketched in the introduction,
adapted to the fact that we are now working directly with $t \to
\infty$. One looks for $u$ in the form $u = u_a + v$. The evolution
equation for $v$ becomes obviously
\beq \label{5.18e} i\partial_t v = - (1/2) \Delta_A v -
\widetilde{R}(u_a) \eeq

\noi where
\bea
\label{5.19e}
\widetilde{R}(u_a) &=& \left ( i \partial_t + (1/2) \Delta_A \right ) u_a \nn \\
&=& - i A \cdot \nabla u_a - (1/2) A^2 u_a \ .
\eea

\noi One takes a large finite $t_0$, one constructs a solution
$v_{t_0}$ of (\ref{5.18e}) with suitably small initial data $v_0 =
v_{t_0}(t_0)$ at $t_0$, and one takes the limit of $v_{t_0}$ as $t_0
\to \infty$. The key of the proof consists in estimating $v_{t_0}$ for
$t \leq t_0$ uniformly in $t_0$ at the required level of regularity.
The estimates thereby obtained remain true in the limit and provide
asymptotic estimates for the solution $u$ at that level. The method
will be presented in detail in Section 6 below in the more interesting
case of the equation (\ref{1.12e}) for $w$ at the level of regularity
of $H^3$ and $H^4$. Here we only provide the basic step thereof,
namely we derive the estimates on $v$. For simplicity we derive them
directly on the limiting $v$, assuming that it tends to zero at
infinity in the relevant norms. As a preliminary step, we need decay
estimates on $\widetilde{R}(u_a)$. \\

\noi {\bf Lemma 5.1.} {\it Let $j \geq 0$ be an integer, let $\alpha$
be a multiindex and let $A$ satisfy the estimates
\beq \label{5.20e}   \parallel \partial_t^l\ 
\partial_x^{\beta} \ A \parallel_r \ \vee \ \parallel \partial_t^l\ 
\partial_x^{\beta} (x \cdot A)\parallel_r  \ \leq C \ t^{-1+2/r}\eeq

\noi for all $r$, $2 \leq r \leq \infty$, for $0 \leq l \leq j$, $0
\leq \beta \leq \alpha$ and for all $t\geq 1$. Let $u_+$, $xu_+ \in H_{\bar{r}}^{2j+|\alpha |}$ for some $\bar{r}$, $1
\leq \bar{r} \leq 2$. \par

Then $\widetilde{R} \equiv \widetilde{R}(U(t)u_+)$ satisfies the
estimate
\beq
\label{5.21e}
\parallel \partial_t^j \ \partial_x^{\alpha} \ \widetilde{R}\parallel_2 \ \leq C\ t^{-5/2+1/r}
\eeq

\noi where $1/r + 1/\bar{r} = 1$, for all $t \geq 1$.}\\

\noi {\bf Proof.} Using the generator of Galilei transformations
\beq
\label{5.22e}
J \equiv J(t) = x + it \nabla = U(t) x U(-t)
\eeq

\noi we rewrite $\widetilde{R}$ as 
\bea
\label{5.23e}
\widetilde{R} &=& - t^{-1} A \cdot (J-x) u_a - (1/2) A^2 u_a \nn \\
&=& - t^{-1} A \cdot U(t) x u_+ + t^{-1} (x \cdot A) U(t) u_+ - (1/2) A^2 U(t) u_+ \ .
\eea

\noi Using the basic estimate of the Schr\"odinger group
\beq
\label{5.24e}
\parallel U(t) f\parallel_r \ \leq (2 \pi |t|)^{- \delta (r)} \ \parallel f\parallel_{\bar{r}}
\eeq

\noi where $2 \leq r \leq \infty$, $1/r + 1/ \bar{r} = 1$ and $\delta (r) = 3/2 - 3/r$, we estimate
$$\parallel \partial_t^j \ \partial_x^{\alpha} \ \widetilde{R}\parallel_2 \
\leq C\ t^{-\delta(r)} \sum_{0 \leq l \leq j} \ \sum_{\beta \leq
\alpha} \Big \{ t^{-1} \parallel \partial_t^l \ \partial_x^{\beta}\ 
A\parallel_s \ \parallel \Delta^{j-l} \ \partial_x^{\alpha - \beta }
(xu_+)\parallel_{\bar{r}}$$ $$+ \left ( t^{-1} \parallel \partial_t^l\ 
\partial_x^{\beta} (x \cdot A)\parallel_s \ + \ \parallel \partial_t^l\ 
\partial_x^{\beta} \ A^2\parallel_s\right ) \parallel
\Delta^{j-l}\ \partial_x^{\alpha - \beta} u_+\parallel_{\bar{r}} \Big \}$$

\noi where $1/s = 1/2 - 1/r$ , 
$$\cdots \leq C \ t^{-5/2 + 1/r}$$

\noi by (\ref{5.20e}). \par\nobreak \hfill $\sq$ \par

We can now derive asymptotic estimates of higher norms for the
solutions of (\ref{1.1e}) with prescribed asymptotic behaviour
$U(t)u_+$ for sufficiently regular $u_+$.\\

\noi {\bf Proposition 5.2.} {\it Let $j \geq 0$ be an integer and let $A$ satisfy
$$\parallel \partial_t^l \
\partial_x^{\alpha} \ A \parallel_r \ \vee \ \parallel \partial_t^l \
\partial_x^{\alpha} (x\cdot A) \parallel_r \ \leq C\ t^{-1+2/r} \eqno (1.29)\equiv(5.25)$$
\noi for all $r$, $2 \leq r \leq \infty$, for $0 \leq l + |\alpha |/2
\leq j$ and for all $t \geq 1$. Let $u_+$, $x u_+ \in H^{2j} \cap
H_{\bar{r}}^{2j}$ for some $\bar{r}$, $1 \leq \bar{r} \leq 2$. Then
there exists a unique solution $u \in \displaystyle{\mathrel{\mathop {\rm \cap}_{0
\leq l \leq j}}}{\cal C}^{j-l} ([1, \infty ) , H^{2l})$ of (\ref{1.1e})
satisfying the estimates
$$\parallel \partial_t^{j-l} \Delta^l(u(t) - U(t)
u_+)\parallel_2 \ \leq C\ t^{-3/2+1/r} \eqno(5.26)$$
\noi where $1/r + 1/ \bar{r} = 1$, for $0 \leq l \leq j$ and for all $t
\geq 1$. The solution is actually unique in ${\cal C}([1, \infty ),
L^2)$ under the condition (5.26) for $j = 0$.}\\

\noi {\bf Proof.} As mentioned above, we concentrate on the derivation
of the estimates. We first prove (5.26) for $l = 0, 1$. For that
purpose we estimate $v = u - U(t) u_+$ by the same method as in the
proof of Proposition 5.1, starting from (\ref{5.18e}) instead of
(\ref{5.5e}), namely with $(u, f)$ replaced by $(v, - \widetilde{R})$.
Defining
$$y_l = \parallel \partial_t^l v \parallel_2\ \hbox{for $0 \leq l \leq
j$}\ , z_l = \parallel \Delta_A  \partial_t^{l-1} v \parallel_2 \
\hbox{for $1 \leq l \leq j$} \ ,\ z_0 = 0 \ ,$$

\noi and using the fact that  
$$\parallel \partial_t^l \widetilde{R} \parallel_2 \
\leq C\ t^{-1-\lambda} \quad \hbox{for $0 \leq l \leq j$} \eqno(5.27)$$

\noi with $\lambda = 3/2 - 1/r$, we obtain (see (\ref{5.10e}) (\ref{5.12e}))
$$|\partial_t y_j | \leq C\ t^{-1} \left \{ (y_{j-1} z_j)^{1/2} +
\sum_{0\leq l < j} (y_l + z_l) \right \} + C \ t^{-1-\lambda}
\eqno(5.28)$$

$$z_j \leq 2y_j + C\ t^{-1} \sum_{0\leq l < j} (y_l + z_l) + C\ t^{-1-\lambda} \ .
\eqno(5.29)$$

\noi We now prove the estimate (5.26) for $l = 0$ by induction
on $j$. The starting point $y_0 \leq C\ t^{-\lambda}$ is obtained by
integrating (5.28) with $j = 0$. We next assume that

$$y_l \vee z_l \leq C\ t^{-\lambda} \qquad \hbox{for $0 \leq l <j$}\ .
\eqno(5.30)$$

\noi Substituting (5.30) into (5.28) (5.29) yields 
$$|\partial_t y_j | \leq C\ t^{-1-\lambda /2}\ z_j^{1/2} + C\ t^{-1-\lambda}
\eqno(5.31)$$

$$z_j \leq 2y_j + C\ t^{-1-\lambda}
\eqno(5.32)$$

\noi from which the result follows by integration over time by the use
of Lemma 2.3. We next prove (5.26) for general $l$, $0 \leq l
\leq j$. In the same way as in the proof of Proposition 5.1, we define
$$z_j^l = \ \parallel \partial_t^{j-l} \Delta^l v \parallel_2 \qquad \hbox{for $0 \leq l \leq j$}$$

\noi and by the same computation, we estimate (see (\ref{5.16e}))
$$z_j^l \leq 2 z_j^{l-1} + C\ t^{-1} \left \{ \left ( z_j^l
z_{j-1}^{l-1}\right )^{1/2} + \sum_{m \leq k \leq j-1} z_k^m \right \}
+ C\ t^{-1-\lambda}\eqno(5.33)$$

\noi from which the result follows as before by induction on $l$ and
$j$. \par\nobreak \hfill $\sq$

\mysection{Scattering theory at the level of H$^{\bf k}$ with k $\geq$ 3 for w}
\hspace*{\parindent}
In this section, we begin the construction of the wave operators for
$u$ at the level of regularity corresponding to $w_*$ or $w$ in $H^k$
for $k \geq 3$ by studying the Cauchy problem for $w$, namely for
(\ref{1.12e}), with initial time zero, at that level of regularity. We
apply the indirect method sketched in the introduction in order to
circumvent the fact that Propositions 3.3 and 3.4 do not apply with
initial time zero under the assumptions made on $B$. More precisely, we
give ourselves an asymptotic behaviour for $w$ at $t = 0$ in the form
of a model $W$ defined in $(0,1]$ and we look for $w$ in the form $w =
W + q$ with $q$ tending to zero as $t \to 0$. The evolution equation
for $q$ becomes (\ref{1.31e}) with $R(\cdot )$ defined by
(\ref{1.32e}). That equation is of the form (\ref{3.1e}) with $(v, A,
V, f)$ replaced by $(q, B, - \hbox{\it \v B}, - R(W))$. In this
section, we prove the existence of $w$ at the relevant level of
regularity under general assumptions on $W$, the most important of
which are decay properties of $R(W)$. In the next section, we shall
construct $W$ satisfying those properties. \par

We first state the relevant result at the level of regularity of $H^3$.
\\

\noi {\bf Proposition 6.1.} {\it Let $I = (0, 1]$ and let $B$ satisfy 
\beq
\label{6.1e}
\parallel \partial_t^j \ \partial_x^{\alpha} \ B \parallel_r \ \vee \ \parallel \partial_t^j
\ \partial_x^{\alpha} \ \hbox{\it \v B} \parallel_r \ \leq b \ t^{-j-|\alpha | + 1/r}
\eeq

\noi for some constant $b$, for $2 \leq r \leq \infty$, $0 \leq j \leq
1$, $0 \leq |\alpha | \leq 1$ and for all $t \in I$. Let $\lambda_0$,
$\lambda_1$ and $\lambda$ satisfy
\beq \label{6.2e} \left \{ \begin{array}{l} \lambda_1 \leq \lambda_0 -
1/2 \\ \\ 0 < \lambda \leq \lambda_1 \wedge (\lambda_0 - 3/2) \wedge
((\lambda_0 + \lambda_1 )/2 - 1)\end{array} \right . \eeq

\noi and let $\lambda '_1 = \lambda_1 \wedge (\lambda_0 - 1)$.\par

Let $W \in {\cal C} (I, H^3) \cap {\cal C}^1 (I, H^1)$ be such that $R
\equiv R(W) \in ({\cal C} \cap H_{\infty}^1)(I, H^1)$ and that $R$
satisfies the estimates
\bea \label{6.3e} & \parallel \partial_t^j R \parallel_2 \ \leq C \
t^{\lambda_j - 1} &\hbox{for $j = 0, 1$} \ , \\ \label{6.4e} &
\parallel \nabla \partial_t^j R \parallel_2 \ \leq C \ t^{\lambda - j}
&\hbox{for $j = 0, 1$} \ , \eea

\noi for all $t \in I$. \par

Then there exists a unique solution $w \in {\cal C}(I,H^3) \cap {\cal
C}^1(I, H^1)$ of (\ref{1.12e}) in $I$ satisfying the estimates
\beq \label{6.5e}  \parallel \partial_t^j (w - W) \parallel_2 \ \leq
C \ t^{\lambda _j} \qquad \hbox{for $j = 0, 1$} \ , \eeq
\beq \label{6.6e}
\parallel \Delta (w - W) \parallel_2 \ \vee \  \parallel \Delta_B (w
- W) \parallel_2 \ \leq C \ t^{\lambda '_1} \ , \eeq
\beq  \label{6.7e} \parallel \nabla
\partial_t (w - W) \parallel_2 \ \vee \  \parallel \nabla \Delta_B (w -
W) \parallel_2 \   \vee \  \parallel \nabla \Delta (w - W) \parallel_2
\ \leq C\ t^{\lambda} \eeq

\noi for all $t \in I$. The solution is actually unique in $(L^{\infty}
\cap {\cal C}_w) (I, L^2)$ under the condition (\ref{6.5e}) for $j =
0$.} \\

\noi {\bf Remark 6.1.} The condition (\ref{6.2e}) is satisfied in
particular by the linear scale $0 < \lambda = \lambda_1 - 1 = \lambda_0
- 2$ which will occur in a natural way for the available $W$. (See
Section 7 below). In that case $\lambda '_1 = \lambda_1$.\\

\noi {\bf Proof.} By (\ref{6.1e}) (\ref{6.3e}) (\ref{6.4e}), the
assumptions of Proposition 3.3 are satisfied in $I = (0, 1]$ (but not
in $[0, 1]$) for the equation (\ref{1.31e}), namely for (\ref{3.1e})
with $(v, A, V, f)$ replaced by $(q, B, - \hbox{\it \v B}, - R)$. Let
$0 < t_0 \leq 1$ and let $q_{t_0} \in {\cal C}(I, H^3) \cap {\cal
C}^1(I,H^1)$ be the solution of (\ref{1.31e}) with $q_{t_0}(t_0) = 0$
obtained from Proposition 3.3. The proof will consist in taking the
limit $t_0 \to 0$ of $q_{t_0}$. For that purpose we shall estimate
$q_{t_0}$ in $H^3$ and $\partial_t q_{t_0}$ in $H^1$ uniformly in $t_0$
for $t_0 \leq t \leq 1$. Those estimates will rely on the identities
(\ref{3.6e}) (\ref{3.19e}) and (\ref{3.42e}) satisfied by $q_{t_0}$ in
$I$. We now estimate $q_{t_0}$, omitting the subscript $t_0$ for
brevity. We define 
\beq \label{6.8e} y_j = \ \parallel \partial_t^j q \parallel_2 \ ,  \quad j =
0, 1\ , \quad y = \parallel \nabla_B \partial_t q \parallel_2 \ . \eeq

\noi We first estimate $y_0$. From (\ref{3.6e}), we obtain
\beq
\label{6.9e}
\left | \partial_t y_0 \right |  \leq \ \parallel R \parallel_2
\eeq

\noi and therefore by integration and by (\ref{6.3e}) with $j = 0$
\beq
\label{6.10e}
y_0 \leq C \ t^{\lambda_0} \qquad \hbox{for $t_0 \leq t \leq 1$} \ .
\eeq

\noi We next estimate $y_1$. From (\ref{3.19e}), in the same way as in
Proposition 3.2, we obtain (see (\ref{3.28e}))
\beq
\label{6.11e}
\left | \partial_t y_1 \right |  \leq \ \parallel f_1 \parallel_2
\eeq

\noi where now (see (\ref{3.18e}))
\beq \label{6.12e} f_1 = i \left ( \partial_t B \right ) \cdot \nabla_B q -
\left ( \partial_t \hbox{\it \v B}\right ) q - \partial_t R \eeq

\noi so that by (\ref{6.1e})
\beq \label{6.13e} \left | \partial_t y_1 \right |  \leq t^{-1}b \left
(  \parallel \nabla_B q \parallel_2\ + \ y_0 \right ) + \ \parallel
\partial_t R \parallel_2 \ . \eeq

\noi Now by direct estimate of (\ref{1.31e})
\beq \label{6.14e} \parallel \nabla_B q \parallel_2^2 \ \leq \ \parallel
q \parallel_2\ \parallel \Delta_B q \parallel_2  \ \leq 2y_0 \left (
y_1 + b\ y_0 \ + \ \parallel R \parallel_2 \right )\eeq

\noi so that for $t \geq t_0$
\bea \label{6.15e} \left | \partial_t y_1 \right |  &\leq & C \ t^{-1}
\left \{ (y_0 \ y_1)^{1/2} + y_0 + \left ( y_0 \parallel R
\parallel_2\right )^{1/2} \right \} + \ \parallel \partial_t R
\parallel_2\nn \\ &\leq & C \left ( t^{-1+ \lambda_0/2} \ y_1^{1/2} +
t^{\lambda_0 - 3/2} + t^{\lambda_1 - 1} \right )\eea

\noi by (\ref{6.3e}) (\ref{6.10e}). Integrating (\ref{6.15e}) by Lemma
2.3 with $y_1(t_0) = \ \parallel R(t_0)\parallel_2$ satisfying
(\ref{6.3e}) yields
\beq \label{6.16e} y_1 \leq C \left ( t^{\lambda_1} +
t_0^{\lambda_0-1}\right ) \leq C\ t^{\lambda '_1} \quad \hbox{for $t_0
\leq t \leq 1$}  \eeq

\noi and in particular
\beq \label{6.17e} y_1 \leq C\ t^{\lambda_1} \quad \hbox{for $t_0 \vee
t_0^{(\lambda_0 - 1)/\lambda_1} \leq t \leq 1$}\ . \eeq

\noi Furthermore it follows from (\ref{6.14e}) (\ref{6.10e}) (\ref{6.16e}) (\ref{6.3e}) that 
\beq
\label{6.18e}
\parallel \Delta_B q \parallel_2 \ \leq C \ t^{\lambda '_1} \ ,
\eeq
\beq
\label{6.19e}
\parallel \nabla_B q \parallel_2 \ \leq C \ t^{(\lambda_0 + \lambda '_1)/2} \ ,
\eeq
\beq \label{6.20e} \parallel \Delta q \parallel_2 \ \leq \ \parallel \Delta_B
q \parallel_2 \ + \ 2b \parallel \nabla_B q \parallel_2 \ + \ b^2 \parallel
q \parallel_2 \ \leq C\ t^{\lambda '_1} \ . \eeq

\noi We shall also need an estimate of $\parallel \nabla_B \otimes
\nabla_B q \parallel_2$. Now 
$$\parallel \nabla_B \otimes \nabla_B q \parallel_2^2 \ = - \sum_j <q,
\nabla_{Bj} \Delta_B \nabla_{Bj} q> \ ,$$
$$\nabla_{Bj} \Delta_B \nabla_{Bj} = \Delta_B^2 - i G_{jl} \nabla_{Bl}
\nabla_{Bj} - i \nabla_{Bl} G_{jl} \nabla_{Bj}$$

\noi where
\beq
\label{6.21e}
G_{jl} = i \left [ \nabla_{Bj} , \nabla_{Bl}\right ] = \nabla_j B_l - \nabla_l B_j
\eeq

\noi so that
$$\parallel \nabla_B \otimes \nabla_B q \parallel_2^2 \ \leq \ \parallel
\Delta_B q \parallel_2^2 \ + \ \parallel G \parallel_{\infty} \
\parallel q \parallel_2 \left ( \parallel\Delta_B q\parallel_2 \ + \
\parallel \nabla_B \otimes \nabla_B q \parallel_2 \right )$$

\noi and therefore
\beq
\label{6.22e}
\parallel \nabla_B \otimes \nabla_B q \parallel_2 \ \leq \ \parallel
\Delta_B q \parallel_2 \ + \ \parallel G \parallel_{\infty} \
\parallel q \parallel_2 \ \leq C \ t^{\lambda '_1}
\eeq

\noi by (\ref{6.10e}) (\ref{6.18e}).\par

We now turn to the estimate of $y$ defined by (\ref{6.8e}). From
(\ref{3.42e}), in the same way as in the proof of Proposition 3.3, we
obtain (see (\ref{3.47e}))
\beq
\label{6.23e}
\left [ \partial_t y \right | \leq \ \parallel g \parallel_2
\eeq

\noi where now (see (\ref{3.41e}))
\bea \label{6.24e} g &=& \left ( \partial_t B - \nabla \hbox{\it \v
B}\right ) \partial_t q + i \left ( \partial_t B \right ) \cdot \nabla_B
\otimes \nabla_B q + i \left ( \nabla \partial_t B \right ) \cdot
\nabla_B q\nn \\ &&- \left ( \partial_t \hbox{\it \v B}\right )
\nabla_B q - \left ( \nabla \partial_t \hbox{\it \v B}\right )q -
\nabla_B \partial_t R \eea

\noi so that
\bea \label{6.25e} &&\left | \partial_t y \right | \leq \left ( \parallel
\partial_t B \parallel_{\infty} \ + \ \parallel \nabla \hbox{\it \v B}
\parallel_{\infty} \right ) y_1 \ + \ \parallel \partial_t B
\parallel_{\infty} \ \parallel \nabla_B \otimes \nabla_B q \parallel_2
\nn \\ &&+ \left ( \parallel \nabla \partial_t B \parallel_{\infty} \ + \
\parallel \partial_t \hbox{\it \v B}\parallel_{\infty} \right )
\parallel \nabla_B q \parallel_2 \ + \ \parallel \nabla \partial_t
\hbox{\it \v B} \parallel_{\infty} \ y_0 \ + \ \parallel \nabla_B
\partial_t R \parallel_2\nn \\ &&\leq  bt^{-1} \  \left \{ 2 y_1 \ + \
\parallel \nabla_B \otimes \nabla_B q \parallel_2 \ + \left ( t^{-1} +
1 \right ) \parallel \nabla_B q \parallel_2 \ + \ t^{-1} y_0 \right \} \nn
\\ &&+ \ \parallel \nabla \partial_t R \parallel_2 \ + \ b \parallel
\partial_t R \parallel_2 \nn \\ &&\leq C \left ( t^{-1 + \lambda '_1} +
t^{-2+(\lambda_0 + \lambda '_1)/2} + t^{-2 + \lambda_0} + t^{-1 +
\lambda} + t^{-1 + \lambda_1} \right ) \leq C\ t^{-1 + \lambda} \eea

\noi by (\ref{6.3e}) (\ref{6.4e}) (\ref{6.10e}) (\ref{6.16e})
(\ref{6.19e}) (\ref{6.22e}) and (\ref{6.2e}). Integrating
(\ref{6.25e}) with initial condition
$$y(t_0) = \ \parallel \nabla_B R(t_0) \parallel_2 \ \leq C \ t^{\lambda}_0$$

\noi yields
\beq
\label{6.26e}
y \leq C\ t^{\lambda} \quad \hbox{for $t_0 \leq t \leq 1$} \ .
\eeq

\noi It follows therefrom and from (\ref{1.31e}) (\ref{6.4e}) (\ref{6.19e}) that
\beq \label{6.27e} \parallel \nabla_B \Delta_B q \parallel_2 \ \leq 2
\left ( y \ + \ \parallel \nabla_B \hbox{\it \v B} q \parallel_2 \ + \
\parallel \nabla_B R \parallel_2 \right ) \leq C\ t^{\lambda} \ . \eeq

\noi Furthermore by similar elementary estimates and by (\ref{6.2e})
\beq
\label{6.28e}
\parallel \nabla \Delta_B q \parallel_2 \ \vee \ \parallel \Delta_B\nabla  q \parallel_2 \ \vee \
\parallel \nabla \Delta q \parallel_2  \ \leq C\ t^{\lambda} \ .
\eeq

\noi We can now take the limit $t_0 \to 0$. We come back to the
original notation $q_{t_0}$ for that part of the argument. The solution
$q_{t_0}$ and its time derivative are estimated in $H^3$ and $H^1$
respectively by (\ref{6.10e}) (\ref{6.20e}) (\ref{6.28e}) and by
(\ref{6.16e}) (\ref{6.26e}) uniformly in $t_0$ for $t_0 \leq t \leq 1$.
Let $0 < t_1 \leq t_0 \leq 1$. Then the $L^2$ norm of the difference
$q_{t_1} - q_{t_0}$ is conserved so that for all $t \in [t_0, 1]$
\beq \label{6.29e} \parallel q_{t_1}(t) - q_{t_0} (t) \parallel_2 \ = \
\parallel q_{t_1}(t_0) \parallel_2 \  \leq C\ t_0^{\lambda_0} \eeq

\noi by (\ref{6.10e}). Therefore $q_{t_0}$ converges in norm in
$L^{\infty}([T, 1],L^2)$ for all $T > 0$ to some $q \in {\cal C} (I,
L^2)$. By the previous uniform estimates and a standard compactness
argument, $q \in L^{\infty} (I, H^3) \cap H_{\infty}^1(I, H^1)$ and $q$
satisfies the same estimates. Furthermore $q$ also satisfies
(\ref{1.31e}) and therefore $q$ can be chosen in ${\cal C}_w(I,H^3)
\cap {\cal C}_w^1(I, H^1)$. By Proposition 3.3, actually $q \in {\cal
C}(I, H^3) \cap {\cal C}^1(I,H^1)$. Together with the estimates, this
proves that $q \in {\cal C}([0, 1] ,H^3) \cap {\cal C}^1([0,1],H^1)$
with $q(0) = \partial_t q(0) = 0$. Returning to the variables $w$
proves the existence part of the proposition and the estimates
(\ref{6.5e})-(\ref{6.7e}), except for the fact that the estimate
(\ref{6.15e}) for $y_1$ used so far for $t_0 \leq t \leq 1$ has
$t^{\lambda ' _1}$ instead of $t^{\lambda_1}$. However the final
estimate (\ref{6.5e}) for $j = 1$ follows from (\ref{6.17e}) in the
limit $t_0 \to 0$. \par

The uniqueness statement in the proposition follows from that of
Proposition 3.3 in $I$ and from (\ref{6.5e}) for $j = 0$.\par\nobreak \hfill $\sq$\par

We now state the corresponding result at the level of regularity of $H^4$. \\

\noi {\bf Proposition 6.2.} {\it Let $I = (0, 1]$ and let $B$ satisfy (\ref{6.1e}) for $2 \leq r \leq \infty$, $0 \leq j + |\alpha |  \leq
2$ and for all $t \in I$. Let $\lambda_j$, $j = 0,1,2$ satisfy
\beq \label{6.30e} \left \{ \begin{array}{l} \lambda_0 > 2 \ , \ 1 < 
\lambda_1 \leq \lambda_0  - 1/2 \\ \\ 0 < \lambda_2 \leq \lambda '_1 -
1/2 \equiv \left ( \lambda_1 - 1/2 \right ) \wedge \left ( \lambda_0 -
3/2\right ) \end{array} \right . \eeq

\noi and let $\lambda '_2 = \lambda_2 \wedge (\lambda_1 - 1)$.\par

Let $W \in {\cal C}^1 (I, H^2) \cap {\cal C}^2 (I, L^2)$ be such that $R
\equiv R(W) \in ({\cal C}^1  \cap H_{\infty}^2)(I, L^2)$ and that $R$
satisfies the estimates
\beq \label{6.31e}  \parallel \partial_t^j R \parallel_2 \ \leq C \
t^{\lambda_j - 1} \qquad \hbox{for $j = 0, 1, 2$}\eeq

\noi for all $t \in I$. Then\par

(1) There exists a unique solution $w \in {\cal C}^1(I,H^2) \cap {\cal
C}^2(I, L^2)$ of (\ref{1.12e}) in $I$ satisfying the estimates (\ref{6.6e}) and
\beq \label{6.32e}  \parallel \partial_t^j (w - W) \parallel_2 \ \leq
C \ t^{\lambda_j} \qquad \hbox{for $j = 0, 1, 2$} \ , \eeq
\beq \label{6.33e}
\parallel \partial_t \Delta_B (w - W) \parallel_2 \ \vee \  \parallel \Delta_B \partial_t (w
- W) \parallel_2 \  \vee \  \parallel \Delta \partial_t (w -
W) \parallel_2 \ \leq C\ t^{\lambda '_2} \eeq

\noi for all $t \in I$. The solution is actually unique in $(L^{\infty}
\cap {\cal C}_w) (I, L^2)$ under the condition (\ref{6.32e}) for $j =
0$.\par

(2) Assume in addition that $B, \hbox{\it \v B}, R \in {\cal C}(I,H^2)$
and that $R$ satisfies the estimate
\beq
\label{6.34e}
\parallel \Delta R \parallel_2 \ \leq C\ t^{\widetilde{\lambda}_2 - 1}
\eeq

\noi for some $\widetilde{\lambda}_2$ satisfying
\beq \label{6.35e} \widetilde{\lambda}_2 \leq \left ( \lambda_0 - 3/2\right )
\wedge \left ( (\lambda_0 + \lambda_1)/2 - 1 \right ) \wedge \lambda_1
\wedge \left ( \lambda_2 + 1\right )
\eeq

\noi and for all $t \in I$. Then $w - W \in {\cal C}(I,H^4)$ and $w$ satisfies
\beq
\label{6.36e}
\parallel \Delta^2 (w - W) \parallel_2 \ \leq C\ t^{\widetilde{\lambda}_2 - 1}
\eeq

\noi for all $t \in I$. In particular $w - W \in {\cal C}([0, 1],
H^{2(1 + \theta )})$ for $0 \leq \theta \leq 1$, $\theta (1 + \lambda
'_1 - \widetilde{\lambda}_2 ) < \lambda '_1$.} \\

\noi {\bf Remark 6.2.} Under the condition (\ref{6.30e}), the
condition (\ref{6.35e}) is satisfied for $\widetilde{\lambda}_2 =
\lambda_2$. Furthermore, the conditions (\ref{6.30e}) (\ref{6.35e}) are
satisfied in particular by the linear scale $0 < \widetilde{\lambda}_2
= \lambda_2 = \lambda_1 - 1 = \lambda_0 - 2$ (see Remark 6.1). In that
case $\lambda '_1 = \lambda_1$ and $\lambda '_2 = \lambda_2$.\\

\noi {\bf Proof.} \underbar{Part (1)}. It follows from (\ref{6.1e})
(\ref{6.31e}) that the assumptions of Proposition 3.4 are satisfied in
$I = (0, 1]$ (but not in $[0, 1])$ for the equation (\ref{1.31e}),
namely for (\ref{3.1e}) with $(v, A, V, f)$ replaced by $(q, B, -
\hbox{\it \v B}, - R)$. The proof proceeds again by taking the limit
$t_0 \to 0$ of a solution $q_{t_0}$ of
(\ref{1.31e}) with suitable data at $t_0$ for some $t_0 \in I$. We
choose the initial data at $t_0$ in the same way as in Proposition 3.4,
namely
\beq \label{6.37e} q_{t_0} (t_0) = q_0 = \left ( t_0^{-1} + b +
K(t_0)\right )^{-1} R(t_0) \eeq

\noi where now $K = - (1/2) \Delta_B - \hbox{\it \v B}$ and $b$ is the
constant occurring in (\ref{6.1e}), so that $K + b$ is a positive
operator. The choice (\ref{6.37e}) is the special case of (\ref{3.66e})
where $(v_0, \rho , f, z) = (q_0, t_0^{-1} + b , - R ,0)$. With that
choice
\beq \label{6.38e} \parallel q_{t_0} (t_0)\parallel_2 \ \leq t_0
\parallel R(t_0) \parallel_2 \ \leq C\ t_0^{\lambda_0} \ , \eeq
\beq \label{6.39e} i \left ( \partial_t q_{t_0}\right )  (t_0) = -
\left ( t_0^{-1} + b \right ) q_0 = - \left ( t_0^{-1} + b \right ) \left (
t_0^{-1} + b + K(t_0)\right )^{-1} R(t_0) \eeq

\noi so that \beq \label{6.40e} \parallel \partial_t q_{t_0}
(t_0)\parallel_2 \ \leq \left ( 1 + b \ t_0\right ) \parallel R(t_0)
\parallel_2 \ \leq C\ t_0^{\lambda_0-1} \eeq

\noi and
\bea \label{6.41e} &&- \left  ( \partial_t^2 q_{t_0}\right )  (t_0) =
\left ( K i \partial_t q_{t_0} + i \left ( \partial_t K\right ) q_{t_0}
- i \partial_t R \right )(t_0) \nn \\ &&=\left ( \left ( - \left ( t_0^{-1} + b
\right ) K + i \partial_t K \right ) \left ( t_0^{-1} + b + K \right
)^{-1} R - i \partial_t R \right ) (t_0) \ . \eea

\noi We need to estimate the last quantity in $L^2$. Since $K + b$ is a
positive operator, we have
$$- b \ t_0 \leq K \left ( t_0^{-1} + b + K\right )^{-1} \leq 1$$

\noi so that 
\beq \label{6.42e} \parallel \left ( t_0^{-1} + b \right )
K \left ( t_0^{-1} + b + K \right )^{-1} R(t_0) \parallel_2 \ \leq
t_0^{-1} \left ( 1 + b \ t_0\right )^2 \parallel R(t_0) \parallel_2 \ .
\eeq

\noi On the other hand, from (\ref{6.1e}) and from 
\beq \label{6.43e}
\partial_t K = i \left ( \partial_t B \right ) \cdot \nabla_B - \partial_t
\hbox{\it \v B} \eeq

\noi we obtain \beq \label{6.44e} \parallel \left ( \partial_t K \right
) q_0 \parallel _2 \ \leq \ t_0^{-1} b \left ( \parallel \nabla_B q_0
\parallel_2 \ + \ \parallel  q_0 \parallel _2 \right ) \ . \eeq

\noi Now
\bea \label{6.45e} &&\parallel \nabla_B q_0 \parallel _2^2 \ = 2 <q_0 ,
(K + \hbox{\it \v B})q_0 >  \ \leq 2 <q_0, (K + b) q_0 > \nn \\ &&= 2
<R(t_0), \left ( t_0^{-1} + b + K \right )^{-2} (K + b) R(t_0) > \ \leq
t_0 \parallel  R(t_0) \parallel _2^2 \ . \eea

\noi Collecting (\ref{6.41e})-(\ref{6.45e}) and using (\ref{6.31e})
(\ref{6.38e}) yields 
\bea \label{6.46e} \parallel \partial_t^2
q_{t_0}(t_0) \parallel_2  &\leq &\left ( t_0^{-1} (1 + b \ t_0)^2 + b
\ t_0^{-1/2} + b \right ) \parallel R(t_0)\parallel _2 \ + \ \parallel 
\partial_t R(t_0)\parallel _2 \nn \\ &\leq & C \left ( t_0^{\lambda_0 -
2} + t_0^{\lambda_1 - 1} \right ) \leq C \ t_0^{\lambda '_1 - 1} \ .
\eea

\noi Let $q_{t_0} \in {\cal C}^1 (I, H^2) \cap {\cal C}^2 (I, L^2)$ be
the solution of (\ref{1.31e}) with initial data (\ref{6.37e}) at $t_0$
obtained from Proposition 3.4. In order to take the limit $t_0 \to 0$,
we need to estimate $q_{t_0}$ and $\partial_t q_{t_0}$ in $H^2$ and
$\partial_t^2 q_{t_0}$ in $L^2$ uniformly in $t_0$ for $t_0 \leq t \leq
1$. We again omit the subscript $t_0$ on $q$ and we define
\beq
\label{6.47e}
y_j = \ \parallel \partial_t^j q \parallel _2 \ .
\eeq

\noi The estimates (\ref{6.38e}) (\ref{6.40e}) (\ref{6.46e}) yield
\beq
\label{6.48e}
y_0(t_0 ) \leq C\ t_0^{\lambda_0} \quad , \quad y_1 (t_0) \leq C\ t_0^{\lambda_0 - 1} \ ,
\eeq
\beq
\label{6.49e}
y_2(t_0) \leq C \ t_0^{\lambda '_1 - 1} \ .
\eeq

\noi We have already estimated $y_0$, $y_1$ and $\parallel \Delta q
\parallel_2$ , $\parallel \Delta_B q \parallel_2$ , in the proof of
Proposition 3.3 by using only (\ref{6.1e}) and (\ref{6.3e}) for $j = 0,
1$. Here the initial condition at $t_0$ is different, but because of
(\ref{6.48e}) it makes no difference in the basic estimates
(\ref{6.10e}) (\ref{6.16e}). It remains only to estimate $y_2$ and
$\parallel  \Delta \partial_t q\parallel _2$. From (\ref{3.74e}) with $j
= 2$, we obtain as in the proof of Proposition 3.4 (see (\ref{3.86e})) 
\beq
\label{6.50e}
\left | \partial_t y_2 \right | \leq \ \parallel f_2\parallel_2
\eeq

\noi where now (see (\ref{3.62e}) (\ref{3.63e}) and (\ref{6.43e}))
\beq \label{6.51e} f_2 = 2i \left ( \partial_t B \right ) \cdot \nabla_B
\partial_t q - 2 \left ( \partial_t \hbox{\it \v B}\right ) \partial_t
q + i \left ( \partial_t^2 B \right ) \cdot \nabla_B q + \left ( \left (
\partial_t B\right )^2 - \left ( \partial_t^2 \hbox{\it \v B}\right )
\right ) q - \partial_t^2 R \eeq

\noi so that by (\ref{6.1e})
\beq \label{6.52e} \left | \partial_t y_2 \right | \leq 2 b t^{-1}
\left ( \parallel \nabla_B \partial_t q \parallel_2 \ + \ y_1 \right ) +
\ bt^{-2} \left ( \parallel  \nabla_B q \parallel_2 \ + \ (1 + b)y_0\right )  + \
\parallel  \partial_t^2 R \parallel_2 \ . \eeq

\noi On the other hand
\beq \label{6.53e} \parallel \nabla_B \partial_t q \parallel_2^2 \ \leq
y_1 \parallel \Delta_B \partial_t q \parallel_2 \ \leq y_1 \left (
\parallel  \partial_t \Delta_B q \parallel_2 \ + \ 2bt^{-1} \parallel
\nabla_B q \parallel _2 \right ) \ . \eeq

\noi By a direct estimate of the time derivative of (\ref{1.31e}), we obtain
\beq \label{6.54e} \parallel \partial_t \Delta_B  q \parallel_2 \ \leq
2 \left ( y_2 + b \left (y_1 + t^{-1} y_0 \right ) + \parallel
\partial_t R \parallel_2 \right )  \eeq

\noi and therefore
\beq
\label{6.55e}
\parallel \nabla_B \partial_t q \parallel_2^2 \ \leq
2y_1 \left ( y_2 + b \left ( t^{-1} \parallel \nabla_B q \parallel_2 \ + \ y_1 + t^{-1} y_0 \right ) +  
\parallel  \partial_t R \parallel_2 \right ) \ .
\eeq

\noi Substituting (\ref{6.55e}) into (\ref{6.52e}) yields
\bea \label{6.56e} \left | \partial_t y_2 \right | &\leq& C\  t^{-1} \Big \{
\left ( y_1 y_2 \right )^{1/2} + \left ( y_1 \ t^{-1} \parallel \nabla_B q \parallel_2 \right )^{1/2} + y_1 + \left ( t^{-1} y_1 \ y_0\right )^{1/2} \nn \\
&&+ \left ( y_1 \parallel  \partial_t R  \parallel_2 \right )^{1/2} + t^{-1}\parallel \nabla_B q \parallel_2\ + \ t^{-1} y_0 \Big \} + \
\parallel  \partial_t^2 R \parallel_2 \eea

\noi so that by (\ref{6.31e}) (\ref{6.10e}) (\ref{6.16e}) (\ref{6.19e})
\bea \label{6.57e} \left | \partial_t y_2 \right | &\leq& C \Big \{ t^{-1+ \lambda '_1/2}
  \ y_2^{1/2} + t^{-1 + \lambda '_1} + t^{-3/2 + (\lambda_1 + \lambda '_1)/2}\nn \\
&&+ \ t^{-2+(\lambda_0 + \lambda '_1)/2} + t^{-2 + \lambda_0} + t^{-1 + \lambda_2} \Big \} \nn \\
&\leq& C\left ( t^{-1 + \lambda '_1/2} \ y_2^{1/2} + t^{-1 + \lambda_2} \right ) \eea

\noi provided
\beq \label{6.58e} \lambda_2 \leq \lambda '_1 \wedge \left ( \left (
\lambda_1 + \lambda '_1 \right )/2 - 1/2 \right ) \wedge \left ( \left
( \lambda_0 + \lambda '_1 \right )/2-1 \right ) \eeq

\noi which reduces to the last inequality in (\ref{6.30e}) by an
elementary computation. Integrating (\ref{6.57e}) by Lemma 2.3 with
initial condition satisfying (\ref{6.49e}) yields
\beq \label{6.59e} y_2 \leq C \left ( t^{\lambda_2} + t_0^{\lambda '_1
- 1} \right ) \qquad \hbox{for $t_0 \leq t \leq 1$} \eeq

\noi and therefore
\beq \label{6.60e} y_2 \leq C\ t^{\lambda_2} \quad \hbox{for $t_0 \vee
t_0^{(\lambda ' _1 - 1)/\lambda_2} \equiv \widetilde{t}_0 \leq t \leq
1$} \ . \eeq

\noi It then follows from (\ref{6.54e}) (\ref{6.10e}) (\ref{6.16e})
(\ref{6.31e}) (\ref{6.60e}) that
\beq \label{6.61e} \parallel \partial_t \Delta_B q \parallel_2 \ \leq C
\ t^{\lambda '_2} \quad \hbox{for $\widetilde{t}_0 \leq t \leq 1$} \eeq

\noi and from (\ref{6.53e}) (\ref{6.19e}) (\ref{6.61e}) that 
\beq
\label{6.62e}
\parallel \Delta_B \partial_t  q \parallel_2 \ \leq C
\ t^{\lambda '_2} \quad \hbox{for $\widetilde{t}_0 \leq t \leq 1$}\ .
\eeq

\noi Finally by (\ref{6.16e}) (\ref{6.62e}) 
\beq \label{6.63e}
\parallel \Delta \partial_t  q \parallel_2 \ \leq \ \parallel \Delta_B
\partial_t  q \parallel_2 \ + \ 2b \left ( y_1 \parallel \Delta_B
\partial_t  q \parallel_2 \right )^{1/2} + \ b^2 y_1 \leq C\ t^{\lambda
'_2} \eeq

\noi for $\widetilde{t}_0 \leq t \leq 1$. The estimates (\ref{6.10e})
(\ref{6.16e}) (\ref{6.20e}) (\ref{6.60e}) (\ref{6.63e}) provide uniform
estimates in $t_0$ of $q$, $\partial_tq$ in $H^1$ and of
$\partial_t^2q$ in $L^2$ for $\widetilde{t}_0 \leq t \leq 1$. We can
now take the limit $t_0 \to 0$. The argument is the same as in the
proof of Proposition 6.1 and will be omitted.\\

\noi \underbar{Part (2)}. The continuity of $q = w - W$ in $H^4$
follows from (\ref{1.31e}), from the continuity of $B$ $\hbox{\it \v
B}$, $R$ in $H^2$ and from Part (1). In order to prove the estimate
(\ref{6.36e}), it is sufficient to estimate $\Delta^2 q_{t_0}$ in $L^2$
uniformly in $t_0$ for $\widetilde{t}_0 \leq t \leq 1$. We omit again
the subscript $t_0$. We know already from (\ref{6.63e}) that 
\beq
\label{6.64e}
\parallel \Delta (Kq - R) \parallel_2 \ \leq C\ t^{\lambda '_2} \ .
\eeq

\noi We next estimate 
\bea \label{6.65e} &&\parallel \Delta^2 q \parallel_2 \ \leq 2
\parallel \Delta Kq \parallel_2 \ + \ 2 \parallel \Delta (B \cdot
\nabla q) \parallel_2 \  + \  \parallel \Delta\left ( ( B^2 - 2
\hbox{\it \v B})q\right ) \parallel_2\nn \\ &&\leq 2 \parallel \Delta
Kq \parallel_2 \ + \ 2b \left ( \parallel \Delta q \parallel_2\ \parallel
\Delta^2 q \parallel_2\right )^{1/2} + \nn \\ &&+ \ C \left ( t^{-2}
\left ( \parallel \nabla_B q \parallel_2\ + \ \parallel q
\parallel_2\right ) + t^{-1} \parallel \Delta q \parallel_2\right ) \ .
\eea

\noi By Lemma 2.2 and (\ref{6.19e}) (\ref{6.20e}) (\ref{6.64e}), this implies 
\beq \label{6.66e} \parallel \Delta^2 q \parallel_2 \ \leq\  C \left ( 
\parallel \Delta R \parallel_2 \ + t^{\lambda '_2} + t^{-2 + (\lambda_0
+ \lambda '_1)/2} + t^{-1 + \lambda '_1}\right ) \leq \ C \
t^{\widetilde{\lambda}_2-1} \eeq

\noi under the assumption (\ref{6.34e}), provided
\beq \label{6.67e} \widetilde{\lambda}_2 \leq \left ( \lambda '_2 + 1
\right ) \wedge \left ( \left ( \lambda_0 + \lambda '_1 \right ) /2 - 1
\right ) \wedge \lambda '_1 \ , \eeq

\noi a condition which reduces to (\ref{6.35e}). The estimate
(\ref{6.66e}) holds for $\widetilde{t}_0 \leq t \leq 1$ uniformly in
$t_0$. \par

The last statement of Part (2) follows from (\ref{6.6e}) and
(\ref{6.66e}) by interpolation.\par\nobreak \hfill $\sq$

\mysection{Choice of W and remainder estimates}
\hspace*{\parindent}
In this section we continue the program started in
Section 6 by constructing model functions $W$ satisfying the
assumptions of Propositions 6.1 and 6.2. In all this section we assume
$A$ to satisfy the free wave equation and therefore to be given by
(\ref{2.12e}) for suitable $(A_+, \dot{A}_+)$. We recall that $A$ and
$B$ are related by (\ref{1.11e}). We first choose $W$ in the form $W =
U(t) \overline{w_+}$ and we obtain sufficient conditions on $(w_+, A_+,
\dot{A}_+)$ to ensure the required assumptions. Those conditions will
require a support condition on $w_+$. We begin by deriving sufficient
conditions in terms of $\chi$, $B$ and $\hbox{\it \v B}$, where $\chi$
is the characteristic function of the support of $w_+$. \\

\noi {\bf Proposition 7.1.} {\it Let $0 < \lambda \leq 1$ and let $w_+
\in H^5$. Let $B$ satisfy (\ref{2.27e}) (equivalently (\ref{6.1e})) for
$2 \leq r \leq \infty$, $0 \leq j + |\alpha | \leq 2$ and in addition
\beq \label{7.1e} \parallel \chi \partial_t^j \ \partial_x^{\alpha}
B(t)\parallel_2 \ \vee \ \parallel \chi \partial_t^j \
\partial_x^{\alpha} \hbox{\it \v B} (t) \parallel_2 \ \leq C\ t^{1 +
\lambda - j - |\alpha |} \eeq

\noi for $0 \leq j + |\alpha| \leq 2$ and for all $t \in (0, 1]$. Then
the following inequality holds
\beq \label{7.2e} \parallel \partial_t^j \ \partial_x^{\alpha}
R(U(t)\overline{w_+}) \parallel_2 \ \leq C\ t^{1+\lambda - j - |\alpha
|} \eeq

\noi for $0 \leq j + |\alpha | \leq 2$ and for all $t \in (0, 1]$.}\\

\noi {\bf Proof.} From (\ref{1.32e}) with $W = U(t)\overline{w_+}$, we
obtain
\beq \label{7.3e} R(W) = - i B \cdot \nabla W + \left ( \hbox{\it \v B}
- B^2/2\right ) W \ . \eeq

\noi It will be sufficient to prove an estimate of the type
\beq \label{7.4e} \parallel \partial_t^j \ \partial_x^{\alpha}
BW\parallel_2 \ \leq C\ t^{1+\lambda - j - |\alpha |} \eeq

\noi for the relevant $j$, $\alpha$, for $w_+$ in $H^4$ and under the
assumptions made on $B$. The final estimate (\ref{7.2e}) will then be
obtained by applying that special case with $(B, w_+)$ replaced by $(B,
\nabla w_+)$, by $(\hbox{\it \v B}, w_+)$ and by $(B^2, w_+)$, given
the fact that the estimates available for $B$ imply the same estimates
for $B^2$. By a Taylor expansion of $U(t)$ to second order

\beq \label{7.5e} BW = B \overline{w_+} + i (t/2) B \Delta
\overline{w_+} - (1/4) B \int_0^t dt'(t-t') U(t') \Delta^2
\overline{w_+} \ , \eeq

\noi we estimate

\beq \label{7.6e} \parallel BW \parallel_2 \ \leq \ \parallel\chi B
\parallel_2\ \parallel w_+ \parallel_{\infty} \ + \ t\parallel \chi B
\parallel_2\ \parallel\Delta w_+\parallel_{\infty} \ + \ t^2 \parallel
B \parallel_{\infty} \ \parallel \Delta^2 w_+ \parallel_2 \eeq

\noi which yields (\ref{7.4e}) for $j = |\alpha | = 0$.\par

In the case $j + |\alpha | \not= 0$, we obtain various terms which we
estimate differently. The terms with all derivatives on $B$ are
estimated exactly as in (\ref{7.5e}) (\ref{7.6e}), with an extra power
$t^{-j-|\alpha |}$ coming from the assumptions on $B$. The terms with
one derivative on $W$ are estimated more simply by a Taylor expansion
of $U(t)$ to first order so that
$$\parallel B\nabla W \parallel_2 \ \leq \ \parallel\chi B
\parallel_2\ \parallel \nabla w_+ \parallel_{\infty} \ + \ t\parallel  B
\parallel_{\infty}\ \parallel\nabla \Delta w_+\parallel_{2} \ \leq \ C\ t$$
$$\parallel B\partial_t W\parallel_2 \ \leq \ \parallel\chi B
\parallel_2\ \parallel \Delta w_+ \parallel_{\infty} \ + \ t\parallel  B
\parallel_{\infty}\ \parallel\Delta^2 w_+\parallel_{2} \ \leq \ C\ t$$

\noi in the cases $j = 0$, $|\alpha | = 1$ and $j = 1$, $\alpha = 0$
respectively, which completes the proof of (\ref{7.4e}) for $j +
|\alpha | = 1$. Similar estimates hold and take care of the terms with
one derivative on $W$ if $j + |\alpha | = 2$. Finally, the terms with
two derivatives on $W$ in the case $j + |\alpha | = 2$ are estimated
simply by
$$\parallel B\partial_t^j \ \partial_x^{\alpha}  W\parallel_2 \ \leq \
\parallel B \parallel_{\infty}\ \parallel \Delta^j \partial_x^{\alpha}
w_+ \parallel_{2}  \ \leq \ C$$

\noi thereby completing the proof of (\ref{7.4e}) in that
case.\par\nobreak \hfill $\sq$ \par

We next complete the argument by giving sufficient conditions on $w_+$
and on $(A_+, \dot{A}_+)$ that ensure (\ref{7.1e}). The following
proposition is a slight extension of Lemma 5.2, part (2) in
\cite{5r}.\\

\noi {\bf Proposition 7.2.} {\it Let $\lambda \geq 0$ and let $w_+$
satisfy the support condition
$${\rm Supp}\ w_+ \subset \{x:||x|-1| \geq \eta \} \eqno(1.33)\equiv (7.7)$$

\noi for some $\eta$, $0 < \eta < 1$. Let $j \geq 0$ be an integer, let
$\alpha$ be a multiindex and let $\chi_R$ be the characteristic
function of the set $\{x:|x|\geq R\}$. Let $(A_+, \dot{A}_+)$ satisfy
$$\parallel \chi_R \partial_x^{\alpha} (x \cdot
\nabla)^{j'}A_+\parallel_2\ \vee \ \parallel \chi_R \partial_x^{\alpha}
(x \cdot \nabla)^{j'}(x \cdot A_+)\parallel_2\ \leq C\ R^{-\lambda -
1/2} \eqno(7.8)$$
$$\parallel \chi_R (x \cdot
\nabla)^{j'}\dot{A}_+\parallel_{6/5}\ \vee \ \parallel \chi_R 
(x \cdot \nabla)^{j'}(x \cdot \dot{A}_+)\parallel_{6/5}\ \leq C\ R^{-\lambda -
1/2} \eqno(7.9)$$

$$\parallel \chi_R \partial_x^{\alpha '} (x \cdot
\nabla)^{j'}\dot{A}_+\parallel_2\ \vee \ \parallel \chi_R \partial_x^{\alpha '}
(x \cdot \nabla)^{j'}(x \cdot \dot{A}_+)\parallel_2\ \leq C\ R^{-\lambda -
1/2} \eqno(7.10)$$

\noi for some $\alpha ' \leq \alpha$ with $|\alpha '| = |\alpha | - 1$
if $\alpha \not= 0$, for all $j'$, $0 \leq j' \leq j$ and for all $R
\geq R_0$ for some $R_0 > 0$. Then (\ref{7.1e}) holds for all $t \in
(0, 1]$.}\\

\noi {\bf Proof.} In the special case $j = 0$ and as regards $B$, the
result is that of Lemma 5.2, part (2) of \cite{5r} to which we refer
for the proof. That proof is a simple application of the finite
propagation speed for the wave equation. The case of general $j \geq 0$
follows therefrom and from (\ref{2.21e}) (\ref{2.19e}). Finally the
result for $\hbox{\it \v B}$ follows from that for $B$ and from
(\ref{2.20e}) (\ref{2.13e}).\par\nobreak \hfill $\sq$ \par

We next try to eliminate the support condition (7.7) on $w_+$.
For that purpose we choose a more complicated $W$. We take $W$ in the
following form
$$W = \left ( 1 - i h \cdot \nabla + \hbox{\it \v h}\right )
W_0 \eqno(7.11)$$

\noi where $h$ and $\hbox{\it \v h}$ are defined by
$$\Delta h = - 2 B \qquad , \quad \Delta \hbox{\it \v h} = - 2 \hbox{\it \v B}
\eqno(7.12)$$

\noi and where $W_0$ is a solution of the equation
$$\left ( i \partial_t + (1/2) \Delta - B^2/2\right ) W_0 = 0 \ .
\eqno(7.13)$$

From (7.12) it follows that $h$ and $\hbox{\it \v h}$ can be
made to satisfy estimates similar to the estimates (\ref{2.27e}) (or
(\ref{6.1e})) for $B$ and $\hbox{\it \v B}$ improved by a factor $t^2$.
This will be proved in Proposition 7.4 below. Anticipating on that fact
and on suitable estimates of $W_0$ which will be proved in Proposition
7.5 below, we now show that $R(W)$ satisfies the estimates required in
Propositions 6.1 and 6.2. \\

\noi {\bf Proposition 7.3.} {\it Let $I = (0, 1]$. Let $B$ and $h$,
$\hbox{\it \v h}$ defined by (7.12) satisfy the estimates
(\ref{2.27e}) (or (\ref{6.1e})) and
$$\parallel \partial_t^j \ \partial_x^{\alpha }h
\parallel_r \ \vee\ \parallel \partial_t^j \ \partial_x^{\alpha
}\hbox{\it \v h} \parallel_r\ \leq C \ t^{2-j-|\alpha|+1/r} \eqno(7.14)$$

\noi for $2 \leq r \leq \infty$, $0 \leq j + |\alpha | \leq 3$ and for
all $t \in I$. Let $r_0 > 3$ and let $W_0 \in {\cal C}(I,H^4) \cap
{\cal C}^1(I,H^3) \cap {\cal C}^2(I,H^2)$ be a solution of
(7.13) in $I$ satisfying the estimates

$$\parallel \partial_t W_0 \parallel_2 \ \vee \
\parallel W_0;H^2\parallel\ \vee\ \parallel\nabla^2 W_0\parallel_{r_0} \
\leq C\ , \eqno(7.15)$$
$$\parallel \partial_t^2 W_0 \parallel_2 \ \vee \
\parallel \partial_t \Delta W_0\parallel_2 \ \vee\ \parallel\nabla\Delta  W_0\parallel_{2} \
\leq C\ t^{-1/2}\ , \eqno(7.16)$$
$$\parallel \partial_t^2 \Delta W_0 \parallel_2 \ \vee \
\parallel \partial_t \nabla \Delta W_0\parallel_2 \ \vee\ \parallel\Delta^2  W_0\parallel_{2} \
\leq C\ t^{-3/2} \eqno(7.17)$$

\noi for all $t \in I$. Then $R(W)$ satisfies the estimates
$$\parallel \partial_t^j  R(W)\parallel_2\ \leq C\
t^{1+ \lambda - j} \qquad \hbox{\it for $0 \leq j \leq 2$} \eqno(7.18)$$
$$\parallel \partial_t^j  \nabla R(W)\parallel_2\ \leq C\
t^{\lambda - j} \qquad \hbox{\it for $j =0, 1$} \eqno(7.19)$$
$$\parallel \Delta  R(W)\parallel_2\ \leq C\
t^{-1+ \lambda} \eqno(7.20)$$

\noi with $\lambda = 1/2 - 1/r_0$, for all $t\in I$. In particular
$R(W)$ satisfies the estimates (\ref{6.3e}) (\ref{6.4e}) and
(\ref{6.31e}) (\ref{6.34e}) of Propositions 6.1 and 6.2 with
$\lambda_j = \lambda + 2 - j$ and $\widetilde{\lambda}_2 = \lambda$.}\\

\noi {\bf Proof.} Substituting (7.11) into the definition (\ref{1.32e}) of $R(W)$ we obtain 
$$R(W) = \left ( 1 - i h \cdot \nabla + \hbox{\it \v
h}\right ) \left ( i \partial_t + (1/2)\Delta - B^2/2\right ) W_0 - ih
B\cdot (\nabla B)W_0$$ $$- i h \left ( - i B \cdot \nabla +
\hbox{\it \v B}\right ) \nabla W_0 + (1 + \hbox{\it \v h}) \left ( - i
B \cdot \nabla + \hbox{\it \v B}\right ) W_0$$
$$- i \left ( i
\partial_t h + (1/2)\Delta h - i B \cdot \nabla h + \nabla h \nabla
\right ) \cdot \nabla W_0$$
$$+ \left ( i \partial_t \hbox{\it \v
h} + (1/2) \Delta \hbox{\it \v h} - i B \cdot \nabla \hbox{\it \v h} +
\nabla \hbox{\it \v h} \cdot \nabla \right ) W_0 \ . \eqno(7.21)$$

\noi By (7.12) (7.13), $R(W)$ reduces to 
$$R(W) = \left ( i \partial_t \hbox{\it \v h} - i B
\cdot \nabla \hbox{\it \v h} - i h B \cdot \nabla B + \hbox{\it \v h}
\hbox{\it \v B}\right ) W_0$$ $$+ \left ( \partial_t h - B \cdot
\nabla h + \nabla \hbox{\it \v h} - i (h \hbox{\it \v B} + \hbox{\it \v
h} B)\right ) \cdot \nabla W_0 - (i \nabla h + hB)\cdot \nabla^2 W_0$$
$$\equiv N_0 W_0 + N_1 \nabla W_0 + N_2 \nabla^2 W_0 \ . \eqno(7.22)$$

\noi The contractions in (7.21) (7.22) have been left
partly unspecified since they will disappear in the estimates.\par

We first estimate
$$\parallel R(W) \parallel_2 \leq  \parallel N_0
\parallel_2 \ \parallel W_0 \parallel_{\infty} \ + \ \parallel N_1
\parallel_2\ \parallel \nabla W_0 \parallel_{\infty} \ + \ \parallel
N_2 \parallel_s \ \parallel \nabla^2 W_0 \parallel_{r_0}$$ $$\leq C
\left ( t^{3/2} + t^{1+\lambda}\right ) \leq C \ t^{1+\lambda} = M_0\eqno(7.23)$$

\noi where $\lambda = 1/s = 1/2 -1/r_0$, by (\ref{2.27e})
(7.14) (7.15). \par

We next estimate $\partial_t R(W)$. The terms with $\partial_t$ applied
to $B$, $\hbox{\it \v B}$ or to $h$, $\hbox{\it \v h}$ are estimated in
a way similar to (7.23), thereby yielding a contribution
$Ct^{-1}M_0$, so that
$$\parallel \partial_tR(W) \parallel_2 \ \leq  C\ t^{-1}M_0 \ + \ \parallel N_0
\parallel_{\infty} \ \parallel \partial_t W_0 \parallel_{2} \ + \ \parallel N_1
\parallel_{\infty} \ \parallel \partial_t \nabla W_0 \parallel_2^{1/2} \  \parallel$$
$$+\ \parallel N_2 \parallel_{\infty} \ \parallel \partial_t \nabla^2 W_0 \parallel_2 \ \leq C
\left ( t^{-1} M_0 + t^{1/2}\right ) \leq C \ t^{\lambda} = M_1\eqno(7.24)$$

\noi by (\ref{2.27e}) (7.14) (7.15) (7.16) and by
using in particular the fact that $\parallel N_j \parallel_{\infty}
\leq Ct$ for $j = 0, 1,2$. \par

We estimate similarly
$$\parallel \nabla R(W) \parallel_2 \ \leq  C\ t^{-1}M_0 \ + \ \parallel N_0
\parallel_{\infty} \ \parallel \nabla  W_0 \parallel_{2} \ + \ \parallel N_1
\parallel_{\infty} \ \parallel \nabla^2 W_0 \parallel_2$$ 
$$+\ \parallel N_2 \parallel_{\infty} \ \parallel  \nabla\Delta  W_0 \parallel_2 \ \leq C
 \ t^{\lambda} = M_1\ .\eqno(7.25)$$

We next estimate the second order derivatives of $R(W)$. The terms with
one or two derivatives applied to $B$, $\hbox{\it \v B}$ or to $h$,
$\hbox{\it \v h}$ are estimated in a way similar to (7.24)
(7.25), thereby yielding a contribution $Ct^{-1}M_1$, and only
the terms with two derivatives on $W_0$ need separate consideration.
Estimating again the $W_0$ factors in $L^2$ and using again the fact
that $\parallel N_j\parallel_{\infty} \leq Ct$, we obtain 
$$\parallel \partial_t^2 R(W) \parallel_2 \leq  C\
t^{-1+\lambda }+ C\ t \left (  \parallel \partial_t^2 W_0\parallel_{2} \
+ \ \parallel \partial_t^2 \nabla W_0 \parallel_{2} \ + \ \parallel
\partial_t^2 \Delta  W_0 \parallel_2\right )$$ $$\leq  C \ t^{-1 +
\lambda} + C\ t^{-1/2} \leq C\ t^{-1+\lambda} \eqno(7.26)$$
$$\parallel \partial_t\nabla  R(W) \parallel_2 \leq  C\
t^{-1+\lambda }+ C\ t \left (  \parallel \partial_t \nabla W_0\parallel_{2} \
+ \ \parallel \partial_t \Delta W_0 \parallel_{2} \ + \ \parallel
\partial_t \nabla \Delta  W_0 \parallel_2\right )$$ $$\leq C \ t^{-1 +
\lambda} \eqno(7.27)$$
$$\parallel \Delta   R(W) \parallel_2 \leq  C\
t^{-1+\lambda }+ C\ t \left (  \parallel \Delta  W_0\parallel_{2} \ + \
\parallel \nabla  \Delta W_0 \parallel_{2} \ + \ \parallel \Delta^2 
W_0 \parallel_2\right )$$ $$\leq  C \ t^{-1 + \lambda} \eqno(7.28)$$

\noi which yields the required second order estimates.\par\nobreak \hfill $\sq$\par

We next derive the estimates (7.14). This will be done
conveniently by using homogeneous Besov spaces \cite{1r} \cite{15r}.
For that purpose we introduce a Paley-Littlewood dyadic decomposition
in the following standard way. Let $F\psi \equiv \widehat{\psi} \in
{\cal C}_0^{\infty} ({I\hskip-1truemm R}^3)$ with $0 \leq
\widehat{\psi} \leq 1$, $\widehat{\psi}(\xi ) = 1$ for $|\xi | \leq 1$
and $\widehat{\psi}(\xi ) = 0$ for $|\xi | \geq 2$. Let
$\widehat{\varphi}_0(\xi ) = \widehat{\psi}(\xi ) - \widehat{\psi}(2\xi
)$ and for any $j \in \Z$, $\widehat{\varphi}_j(\xi ) =
\widehat{\varphi}_0 (2^{-j}\xi)$ so that $\widehat{\varphi}_j$ is
supported in $\{\xi : 2^{j-1} \leq |\xi | \leq 2^{j+1}\}$ and for any
$\xi \in {I\hskip-1truemm R}^3\setminus \{0\}$, 
$$\sum_j \widehat{\varphi}_j(\xi ) = 1$$

\noi with at most two nonvanishing terms in the sum for each $\xi$. The homogeneous
Besov space $\dot{B}_{r,s}^{\rho}$ is then defined for any $\rho \in
{I\hskip-1truemm R}$ and $1 \leq r$, $s \leq \infty$ by
$$\dot{B}_{r,s}^{\rho} = \{v: \parallel v;
\dot{B}_{r,s}^{\rho}\parallel\ = \ \parallel2^{\rho j} \varphi_j *
v;l_j^s(L_x^r)\parallel\ < \infty \}  \eqno(7.29)$$

\noi where $F\varphi_j = \widehat{\varphi}_j$ and $*$ denotes the
convolution in ${I\hskip-1truemm R}^3$.\par

We can now state the result as follows.\\

\noi {\bf Proposition 7.4.} {\it Let $I = (0, 1]$, let $j \geq 0$ be an
integer and let $\alpha$ be a multiindex. \par

(1) Assume that $(A_+, \dot{A}_+)$ satisfies the conditions
$${\cal A} \in \dot{B}_{2,2}^{-2} \cap
\dot{B}_{1,1}^{0}\qquad , \qquad \dot{\cal A} \in \dot{B}_{2,2}^{-3} \cap
\dot{B}_{1,1}^{-1} \eqno(7.30)$$
\noi for 
$${\cal A} = \partial_x^{\alpha}(x \cdot \nabla)^{j'} A_+
\qquad , \qquad {\cal A} = \partial_x^{\alpha}(x \cdot \nabla)^{j'} (x\cdot
A_+) \eqno(7.31)$$
$$\dot{\cal A} = \partial_x^{\alpha}(x \cdot \nabla)^{j'} \dot{A}_+
\qquad , \qquad \dot{\cal A} = \partial_x^{\alpha}(x \cdot \nabla)^{j'} (x\cdot
\dot{A}_+) \eqno(7.32)$$
\noi for $0 \leq j' \leq j$. Then $h$, $\hbox{\it \v h}$ satisfy the
estimates (7.14) for $2 \leq r \leq \infty$ and for all $t \in
I$.\par

(2) Let ${\cal A}$ and $\dot{\cal A}$ satisfy
$$\omega^{\varepsilon} {\cal A} \in L^1 \quad , \quad
<x>^{\theta} {\cal A} \in L^1 \quad , \quad \int dx\ {\cal A}(x) = 0 \eqno(7.33)$$
$$<x>^{1+ \theta} \dot{\cal A} \in L^1 \quad , \quad \int
dx \ \dot{\cal A}(x) = \int dx\ x\dot{\cal A}(x) = 0 \eqno(7.34)$$
\noi for some $\theta > 1/2$. Then ${\cal A}$, $\dot{\cal A}$ satisfy (7.30).}\\

\noi {\bf Proof.} From (\ref{1.11e}) and (7.12), it follows that
$$2\omega^{-2} A = - 2 t^{-1} \ \omega^{-2} D_0B(1/t)$$
$$= - t^{-1} D_0 \left (t^{-2}h\right ) (1/t)
\eqno(7.35)$$

\noi and similarly
$$2\omega^{-2} (x\cdot A) =  -  t^{-1} D_0 \left (t^{-2}\hbox{\it \v h}\right ) (1/t) \ .
\eqno(7.36)$$

We are now in the same situation as in Lemma 2.4, where the estimates
(\ref{2.27e}) for $B$, $\hbox{\it \v B}$ were obtained from the estimates
(\ref{2.26e}) for $A$. In order to prove the estimates (7.14)
for $h$, $\hbox{\it \v h}$, it is sufficient to show that
$\omega^{-2}A$ and $\omega^{-2}(x \cdot A)$ satisfy the same estimates
(\ref{2.26e}) as $A$ and $(x\cdot A)$. Since $x \cdot A$ satisfies the
wave equation as well as $A$, it is sufficient to prove those estimates
for $A$. Furthermore by (\ref{2.15e})-(\ref{2.18e}), it is sufficient
to consider the case $j = 0$, $\alpha = 0$.\\

\noi \underbar{Part (1)}. From the basic estimate (see (3.13) in [4])
$$\parallel \exp (i\omega t) \varphi_j * f \parallel_r \ \leq C
|t|^{-1+2/r} \ 2^{j(2-4/r)} \parallel\varphi_j * f
\parallel_{\bar{r}}$$

\noi with $2 \leq r \leq \infty$, $1/r + 1/\bar{r} = 1$, we obtain
$$\parallel \varphi_j * \omega^{-2} A\parallel_r
\leq C\ t^{-1+2/r} \Big \{ 2^{-4j/r} \parallel \varphi_j * A_+
\parallel_{\bar{r}}$$ $$+ 2^{-j (1+4/r)} \parallel\varphi_j *
\dot{A}_+ \parallel_{\bar{r}} \Big \} . \eqno(7.37)$$

\noi Taking the $l^2$ norm for $r = 2$ and the $l^1$ norm for $r = \infty$ yields 
$$\parallel \omega^{-2}A; \dot{B}_{2,2}^0 \parallel\ \leq C \left (
\parallel A_+; \dot{B}_{2,2}^{-2} \parallel\ + \ \parallel \dot{A}_+;
\dot{B}_{2,2}^{-3} \parallel\right ) \ ,$$
$$\parallel \omega^{-2}A; \dot{B}_{\infty ,1}^0 \parallel\ \leq C \ t^{-1}\left (
\parallel A_+; \dot{B}_{1,1}^{0} \parallel\ + \ \parallel \dot{A}_+;
\dot{B}_{1,1}^{-1} \parallel\right ) \ .$$

\noi By interpolation and by the standard embedding properties of Besov spaces, this implies
$$\parallel \omega^{-2}A\parallel_r \ \leq C\ t^{-1+2/r}$$

\noi where the last constant depends on the relevant norms of $(A_+,
\dot{A}_+)$.\\

\noi \underbar{Part (2)}. As in Part (1) it is sufficient to consider
the case where $({\cal A}, \dot{\cal A}) = (A_+, \dot{A}_+)$. Using the Young
inequality and the homogeneity relation
$$\parallel \varphi_j\parallel_{\bar{r}}\ = 2^{3j/r} \parallel
\varphi_0\parallel_{\bar{r}}$$

\noi we estimate the bracket in (7.37) by 
$$ 2^{-4j/r} \parallel \varphi_j *
A_+\parallel_{\bar{r}} \ + \ 2^{-j(1+4/r)} \parallel \varphi_j *
\dot{A}_+ \parallel_{\bar{r}}$$ $$\leq C \left \{ 2^{-\varepsilon j
- j/r} \parallel\omega^{\varepsilon}A_+ \parallel_1 \ + \ 2^{-j(1+1/r)}
\parallel\dot{A}_+ \parallel_1 \right \}\ . \eqno(7.38)$$

\noi Therefore the high frequency part of the Besov norms, more
precisely the summation over $j \geq 0$, is controlled by the
conditions $\omega^{\varepsilon} A_+ \in L^1$, $\dot{A}_+ \in L^1$.\par

We now consider the low frequency part of the Besov norms. Using the
vanishing integral condition on $A_+$, we rewrite
$$\left ( \varphi_j * A_+\right ) (x) = \int dy \left ( \varphi_j (x -
y) - \varphi_j (x)\right ) A_+(y)$$

\noi so that 
$$\parallel \varphi_j * A_+\parallel_{\bar{r}}\ \leq
\ \mathrel{\mathop {\rm Sup}_{y}}\ |y|^{-\theta} \ \parallel \varphi_j
(. - y) - \varphi_j \parallel_{\bar{r}} \ \parallel |x|^{\theta}
A_+\parallel_{1}$$ $$= 2^{j(\theta + 3/r)}  \mathrel{\mathop {\rm
Sup}_{y}}\ |y|^{-\theta} \ \parallel \varphi_0 (. - y) - \varphi_0
\parallel_{\bar{r}} \ \parallel |x|^{\theta} A_+\parallel_{1}$$ $$=
C\ 2^{j(\theta + 3/r)} \parallel |x|^{\theta}A_+ \parallel_1 \eqno(7.39)$$

\noi by homogeneity and the fact that the last Sup is finite for $0
\leq \theta \leq 1$. This implies the summability over $j \leq 0$ of
the contribution of $A_+$ to the bracket in (7.37) for $2 \leq r
\leq \infty$ and $\theta > 1/2$. Similarly using the vanishing integral
conditions on $\dot{A}_+$, we rewrite
$$\left ( \varphi_j * \dot{A}_+\right ) (x) = \int dy \left ( \varphi_j (x -
y) - \varphi_j (x) + y \cdot \nabla \varphi_j (x) \right ) \dot{A}_+(y)$$

\noi so that
$$\parallel \varphi_j *
\dot{A}_+\parallel_{\bar{r}}\ \leq \ \mathrel{\mathop {\rm Sup}_{y}}\
|y|^{-(1 + \theta )} \ \parallel \varphi_j (. - y) - \varphi_j + y
\cdot \nabla \varphi_j \parallel_{\bar{r}} \ \parallel |x|^{1 + \theta}
\dot{A}_+\parallel_{1}$$ $$= 2^{j(1 + \theta + 3/r)}  \mathrel{\mathop
{\rm Sup}_{y}}\ |y|^{-(1 + \theta )} \ \parallel \varphi_0 (. - y) -
\varphi_0 + y \cdot \nabla \varphi_0 \parallel_{\bar{r}} \ \parallel
|x|^{1 + \theta} \dot{A}_+\parallel_{1}$$ $$= C\ 2^{j(1 + \theta +
3/r)} \parallel |x|^{1 + \theta}\dot{A}_+ \parallel_1 \eqno(7.40)$$

\noi by homogeneity and by the finiteness of the last Sup for $0 \leq
\theta \leq 1$. This implies the summability over $j \leq 0$ of the
contribution of $\dot{A}_+$ to the bracket in (7.37) for $2 \leq
r \leq \infty$ and $\theta > 1/2$.\par\nobreak \hfill $\sq$\par

We now turn to the study of the equation (7.13) and we prove
that it admits solutions $W_0$ satisfying the requirements of
Proposition 7.3. We rewrite that equation in a form similar to
(\ref{3.1e}), namely
$$i \partial_t v = - (1/2) \Delta v + Vv
\eqno(7.41)$$
\noi with $V = (1/2)B^2$. \\

\noi {\bf Proposition 7.5.} {\it Let $I = (0, 1]$. Let $V \in {\cal
C}(I,H^4)\cap {\cal C}^1(I,H^2 \cap L^{6/5}) \cap {\cal C}^3(I, L^2)$
satisfy the estimates
$$\parallel \partial_t^j \ \partial_x^{\alpha} V \parallel_r \ \leq C \ t^{-j-|\alpha | + 1/r}
\eqno(7.42)$$

\noi for $0 \leq j + |\alpha | \leq 1$ and $r = \infty$, for $0 \leq j
+ |\alpha | \leq 3$ and $r = 2$, and for $\alpha = 0$, $j = 1$ and $6/5
\leq r \leq 2$. Let $v_1 \in H^6$. Then there exists a unique solution
$v \in {\cal C}(I, H^4) \cap {\cal C}^1(I, H^3) \cap {\cal C}^2(I, H^2)
\cap {\cal C}^3(I, L^2)$ of (7.41) in $I$ with $v(1) = v_1$,
satisfying the following estimates
$$\parallel v(t)\parallel_2 \ = \ \parallel v_1\parallel_2
$$ 
$$\parallel \partial_t v \parallel_2 \ \vee \  \parallel
\Delta v\parallel_2\ \leq C \eqno(7.43)$$
$$\parallel \partial_t^2 v
\parallel_2 \ \vee \  \parallel \partial_t \Delta v\parallel_2\ \ \vee
\ \parallel \nabla \Delta v\parallel_2\ \leq C\ t^{-1/2} 
\eqno(7.44)$$
$$\parallel \partial_t^3 v \parallel_2 \ \vee \ 
\parallel \partial_t^2 \Delta v\parallel_2\  \vee \ \parallel
\partial_t \nabla \Delta v\parallel_2\ \vee \ \parallel \Delta^2
v\parallel_2\ \leq C\ t^{-3/2} \eqno(7.45)$$ 
$$\parallel \Delta v \parallel_r \ \leq
\left \{ \begin{array}{ll} C &\hbox{\it for $2 \leq r < 4$} \\ \\ C \
\ell n \ t &\hbox{\it for $r = 4$} \\ \\ C\ t^{-1/2 + 2/r} &\hbox{\it
for $4 < r < 6$} \end{array}\right . \eqno(7.46)$$

\noi for all $t \in I$. The solution is actually unique in ${\cal
C}(I,L^2)$.}\\

\noi {\bf Proof.} The existence can be proved by an extension of the
method of Proposition 3.4 to the level of $H^6$ using the third
derivative $\partial_t^3v$, simplified by the fact that here $A = 0$
and $f = 0$. For that purpose, one has to ensure in particular that
$(\partial_t^3v)(1) \in L^2$. Now by an easy computation

$$ \left ( - i \partial_t^3 v \right ) (1) = \left (
K^3 + 3 i \left ( \partial_t V \right ) K - i \left ( \nabla
\partial_tV\right ) \cdot \nabla - (i/2) \left ( \Delta \partial_t
V\right ) - \partial_t^2 V\right ) v_1 \eqno(7.47)$$

\noi where now $K = - (1/2) \Delta + V$. Under the assumptions made on
$V$, this belongs to $L^2$ for $v_1 \in H^6$. In particular the
condition $V \in {\cal C}(I,H^4)$ implies that ${\cal D}(K^3) = H^6$.
The solution $v$ comes out with additional regularity properties that
are of no interest here and have not been stated. We skip the details
and we concentrate on the derivation of the estimates
(7.43)-(7.46). \par

In the same way as in the proof of Proposition 3.2, we first estimate
\begin{eqnarray*}
\parallel \Delta v  \parallel_2 & \leq & 2 \left ( \parallel
\partial_t v\parallel_2\ + \ \parallel V\parallel_{\infty}\ \parallel
v\parallel_2\right ) \ , \\
\partial_t\parallel \partial_t v \parallel_2 & \leq & \parallel
\partial_t V\parallel_2\  \parallel v\parallel_{\infty}\ \leq \
C\parallel \partial_t V \parallel_2 \ \parallel v\parallel_2^{1/4}\
\parallel \Delta v\parallel_2^{3/4}\\
&\leq & C \parallel \partial_t V\parallel_2 \left ( \parallel \partial_t v\parallel_2\ + \ C \right )^{3/4}\\
&\leq & C\ t^{-1/2} \left ( \parallel \partial_t v\parallel_2\ + \ C \right )^{3/4}
\end{eqnarray*}

\noi from which (7.43) follows by integration.\par

We next estimate as in the proof of Proposition 3.4
\begin{eqnarray*}
\partial_t\parallel \partial_t^2 v \parallel_2 & \leq & \parallel
\partial_t^2 V\parallel_2\  \parallel v\parallel_{\infty}\ + \
2\parallel \partial_t V \parallel_{\infty}\  \parallel \partial_t v\parallel_2 \\
&\leq & C\left (  t^{-3/2} + t^{-1}\right ) \leq C\ t^{-3/2} 
\end{eqnarray*}

\noi which implies the first inequality in (7.44) by integration. Furthermore
\begin{eqnarray*} \parallel \partial_t \Delta v \parallel_2 & \leq & 2
\left \{ \parallel \partial_t^2 v\parallel_2\  + \ \parallel \partial_t
V\parallel_{2}\ \parallel v \parallel_{\infty}\  + \ \parallel
V\parallel_{\infty} \ \parallel \partial_t v\parallel_{2}\right \} \\
&\leq & C\left (  t^{-1/2} + 1\right ) \leq C\ t^{-1/2} \end{eqnarray*}
\begin{eqnarray*} \parallel \nabla \Delta v \parallel_2 & \leq & 2
\left \{ \parallel \nabla \partial_t v\parallel_2\  + \ \parallel \nabla
V\parallel_{2}\ \parallel v \parallel_{\infty}\  + \ \parallel
V\parallel_{\infty} \ \parallel\nabla v\parallel_{2}\right \} \\
&\leq & C\left (  t^{-1/4} + t^{-1/2} + 1\right ) \leq C\ t^{-1/2} \end{eqnarray*}

\noi which completes the proof of (7.44).

We next estimate in a similar way
\begin{eqnarray*} \partial_t\parallel \partial_t^3 v \parallel_2 & \leq
& \parallel \partial_t^3 V\parallel_2\  \parallel v\parallel_{\infty}\
+ \ 3\parallel \partial_t^2 V \parallel_{2}\  \parallel \partial_t
v\parallel_{\infty}\ +\  3  \parallel \partial_t V\parallel_{\infty}\
\parallel \partial_t^2 v\parallel_{2}\\ &\leq & C\left (  t^{-5/2} +
t^{-15/8} + t^{-3/2}\right ) \leq C\ t^{-5/2} \end{eqnarray*}

\noi which implies the first inequality in (7.45) by integration. Furthermore
\begin{eqnarray*} &&\parallel \partial_t^2 \Delta v \parallel_2  \ \leq 
\ 2 \Big \{ \parallel \partial_t^3 v\parallel_2\  + \ \parallel
\partial_t^2 V\parallel_{2}\ \parallel v \parallel_{\infty}\  + \
2\parallel \partial_t V\parallel_{\infty} \ \parallel \partial_t
v\parallel_{2} \\ &&+\ \parallel V \parallel_{\infty}\ \parallel
\partial_t^2  v \parallel_2\Big\} \leq C\left (  t^{-3/2} + t^{-3/2} +
t^{-1}+ t^{-1/2}\right ) \leq C\ t^{-3/2} \ , \end{eqnarray*}
\begin{eqnarray*} \parallel \partial_t \nabla  \Delta v \parallel_2  & \leq& 
2 \Big \{ \parallel \partial_t^2 \nabla v\parallel_2\  + \ \parallel
\partial_t\nabla  V\parallel_{2}\ \parallel v \parallel_{\infty}\  + \
\parallel \partial_t V\parallel_{\infty} \ \parallel \nabla v \parallel_{2}\\
&& +\ \parallel \nabla V \parallel_{\infty}\ \parallel
\partial_t  v \parallel_2\ + \ \parallel V\parallel_{\infty}\ \parallel \partial_t \nabla v \parallel_2 \Big\} \\
&\leq& C\left (  t^{-1} + t^{-3/2} +
t^{-1}+ t^{-1/4}\right ) \leq C\ t^{-3/2} \ , \end{eqnarray*}
\begin{eqnarray*} &&\parallel \Delta^2  v \parallel_2  \ \leq 
\ 2 \Big \{ \parallel \partial_t \Delta v\parallel_2\  + \ \parallel
\Delta V\parallel_{2}\ \parallel v \parallel_{\infty}\  + \
2\parallel \nabla V\parallel_{\infty} \ \parallel \nabla
v\parallel_{2} \\ &&+\ \parallel V \parallel_{\infty}\ \parallel
\Delta  v \parallel_2\Big\} \leq C\left (  t^{-1/2} + t^{-3/2} +
t^{-1}+ 1\right ) \leq C\ t^{-3/2} \ , \end{eqnarray*}

\noi which completes the proof of (7.45).\par

We finally prove (7.46). Since
$$\parallel \Delta  v \parallel_r  \ \leq 
2 \left ( \parallel \partial_t  v\parallel_r\  + \ \parallel
 V\parallel_{\infty}\ \parallel v \parallel_{r}\right ) \  \leq \
\left (\parallel \partial_t v\parallel_{r} \ + \ C \right )
\eqno(7.48)$$

\noi it suffices to estimate $\partial_tv \in L^r$. We start from the
integral relation 
$$i\partial_t v = U(t-1) \left ( -(1/2)\Delta +
V(1)\right ) v_1 - \int_t^1 dt'\ U(t-t') \partial_t (Vv)(t')\ . \eqno(7.49)$$

\noi We estimate for $2 \leq r \leq 6$
$$\parallel U(t-1) \left ( -(1/2)\Delta + V(1)\right )
v_1 \parallel_r \ \leq C\parallel \Delta v_1;H^1 \parallel\ + \
\parallel  V(1)v_1;H^1 \parallel \ \leq C\ . \eqno(7.50)$$

\noi On the other hand from the basic estimate (\ref{5.24e}) of the
Schr\"odinger evolution group with $2 \leq r \leq \infty$, $1/r + 1/\bar{r} = 1$
and $\delta (r) = 3/2 - 3/r$, we obtain
$$\parallel \int_t^1 dt'\ U(t-t') \partial_t (Vv) (t'
)\parallel_r$$ $$\leq C \int_t^1 dt'(t-t')^{-\delta} \left (
\parallel\partial_tV(t') \parallel_{\bar{r}} \ \parallel v(t')
\parallel_{\infty}\ + \ \parallel V(t')\parallel_{3/\delta} \ \parallel
\partial_t v(t') \parallel_2 \right ) $$
$$\leq C \int_t^1
dt'(t-t')^{-\delta} \left ( t'^{-1/r} + t'^{1/2-1/r} \right ) \eqno(7.51)$$

$$\int_t^1 dt' (t-t')^{-\delta} \ t'^{-1/r} =
\int_t^{2t\wedge 1} + \int_{2t \wedge 1}^1 \leq C \left \{
t^{-1/2+2/r} + \int_{2t \wedge 1}^1 dt' \ t'^{-3/2+2/r} \right \} \ .
\eqno(7.52)$$

\noi Integrating (7.52) and collecting the result and
(7.48) (7.50) (7.51) yields (7.46).\\

\noi {\bf Remark 7.1.} The assumptions on $V$ in Proposition 7.5 are
unnecessarily restrictive. Under the condition $V \in {\cal C}^2(I,L^2
+ L^{\infty})$, $\partial_t^3 V \in L_{loc}^1(I,L^2+L^{\infty})$, one
can prove the existence of a unique solution $v \in {\cal C}^3(I,L^2)
\cap {\cal C}^2(I, H^2)$ of (7.41) in $I$. Furthermore the
additional estimates (7.43)-(7.46) on $v$ can be derived
by using only a subset of the assumptions made on $V$. On the other
hand, the assumptions made on $V$ are easily seen to follow from
sufficient regularity assumptions on $B$ and from estimates of the type
(\ref{2.27e}). \par

We now discuss briefly the situation that arises from Propositions
7.3-7.5 as regards the construction of an asymptotic $W$ satisfying the
assumptions of Propositions 6.1 and 6.2. We have actually constructed a
$W$ satisfying the required assumptions by (7.11) (7.12)
and by taking for $W_0$ the solution $v$ of (7.41) obtained in
Proposition 7.5. However that $W$ cannot be parametrized by the
asymptotic state $u_+$ of $u$, and is parametrized instead by $v_1 =
W_0(1)$. Now it follows from Proposition 3.2 that the Cauchy problem
for (7.41) is well posed in $[0,1]$ in $H^2$. In particular the
previous $W_0$ has an $H^2$ limit $\overline{w_+} \in H^2$ as $t \to 0$,
which can be identified with $\overline{Fu_+}$, and $W_0$ can be
reconstructed from that $\overline{w_+}$ by solving the $H^2$ Cauchy
problem with initial time zero. The weak point however is that we are 
unable to characterize those $w_+ \in H^2$ that arise from solutions of
(7.41) with the regularity at the level of $H^6$ which is
required for the needs of Propositions 6.1 and 6.2.

\mysection{Wave operators and asymptotics for u}
\hspace*{\parindent}
In this section we collect the implications of Sections 6 and 7 on the
theory of scattering at the level of regularity of $FH^3$ and $FH^4$
for (\ref{1.1e}). The main results consist in obtaining solutions of
(\ref{1.1e}) with given asymptotic behaviour at infinity in time, which
is essentially equivalent to the construction of the wave operators. We
consider separately the simple case of asymptotics provided by
Propositions 7.1 and 7.2 and the more complicated case provided by
Propositions 7.3-7.5. On the other hand since in both cases we derive
all the required estimates on $R(W)$ for Propositions 6.1 and 6.2
together, we also state the implications of those two
propositions together. In all this section $A$ is assumed to be a solution of
the free wave equation. \par

We begin with the case of the simple asymptotics provided by $W = U(t)
\overline{w_+}$.\\

\noi {\bf Proposition 8.1.} {\it Let $A$ be a solution of the free wave
equation satisfying the estimates (\ref{2.26e}) and the decay
assumptions (7.8)-(7.10) for $0 \leq j + |\alpha | \leq
2$ and for some $\lambda \in (0,1)$. Let $u_+ \in FH^5$ and let $w_+ =
Fu_+$ satisfy the support condition (7.7). Then there exists a
unique solution $u$ of (\ref{1.1e}) such that $\widetilde{u} \in {\cal
C}([1, \infty ), FH^4) \cap {\cal C}^1([1 , \infty ),FH^2) \cap {\cal
C}^2([1, \infty ), L^2)$ and satisfying the estimates
\beq
\label{8.1e}
\parallel \widetilde{u} - u_+ \parallel_2 \ \leq C\ t^{-2-\lambda}
\eeq
\beq \label{8.2e} \parallel |x|^{2+j} (\widetilde{u} - u_+) \parallel_2
\ \leq C\ t^{-1-\lambda + j}\qquad \hbox{\it for $0 \leq j \leq 2$} \eeq
\beq \label{8.3e} \parallel |x|^{j} \partial_t (\widetilde{u} - u_+) \parallel_2
\ \leq C\ t^{-3-\lambda + j}\qquad \hbox{\it for $0 \leq j \leq 1$} \eeq
\beq \label{8.4e} \parallel x^{2} \left ( i t^2 \partial_t + (1/2) x^2
\right ) (\widetilde{u} - u_+) \parallel_2 \ \leq C\ t^{-\lambda } \eeq
\beq \label{8.5e} \parallel \partial_t \left ( i t^2 \partial_t + (1/2)  x^2
\right ) (\widetilde{u} - u_+) \parallel_2 \ \leq C\ t^{-2-\lambda } \eeq

\noi for all $t \geq 1$. The solution is actually unique in ${\cal
C}([1, \infty ), L^2)$ under the condition (\ref{8.1e}).}\\

\noi {\bf Proof.} The results follow immediately from Propositions 6.1,
6.2 and 7.1, 7.2 through the change of variables (\ref{1.19e}). The
latter implies in particular that

\beq \label{8.6e} U(-t) i \partial_t w(t) = \overline{F((it^2\partial_t
+ (1/2) x^2) \widetilde{u})(1/t)}\ . \eeq

\noi The estimates (\ref{8.1e})-(\ref{8.5e}) are a rewriting of
(\ref{6.5e})-(\ref{6.7e}) and of (\ref{6.32e}) (\ref{6.33e})
(\ref{6.36e}).\par\nobreak \hfill $\sq$ \par

We now turn to the more complicated situation covered by Propositions 7.3-7.5.\\

\noi{\bf Proposition 8.2.} {\it Let $A$ be a solution of the free wave
equation with $A \in {\cal C} ({I\hskip-1truemm R}^+, H^4)$ satisfying the estimates
(\ref{2.26e}) and the conditions (7.33) (7.34) for $0
\leq j + |\alpha | \leq 3$. Let $0 < \lambda < 1/4$. Let $v_1 \in H^6$,
let $W_0$ be the solution of (7.13) in $(0, 1]$ with $W_0(1) =
v_1$ obtained from Proposition 7.5. Define $W$ by (7.11)
(7.12) and $\widetilde{u}_a$ by
\beq \label{8.7e} F\widetilde{u}_a = \overline{\widetilde{W}(1/t)}\ .
\eeq

\noi Then the same conclusions as in Proposition 8.1 hold with $u_+$
replaced by $\widetilde{u}_a$ in all the estimates.}\\

\noi {\bf Proof.} The results follow immediately from Propositions 6.1,
6.2 and 7.3-7.5 by the change of variables (\ref{1.19e}). \par\nobreak \hfill $\sq$\par

\par \vskip 1 truecm

\noi {\large\bf  Acknowledgements.} One of us (J. G.) is grateful to
Professor Takeshi Wada for enlightening conversations at an early stage
of this work.

\end{document}